\documentclass[a4paper,11pt]{amsart}
\usepackage{amssymb,amsfonts,amsxtra,amscd}
\usepackage[all]{xy}
\usepackage{fullpage}
\theoremstyle{plain}
\newtheorem{theorem}{Theorem}[section]
\newtheorem{lemma}[theorem]{Lemma}
\newtheorem{cor}[theorem]{Corollary}
\newtheorem{prop}[theorem]{Proposition}
\theoremstyle{definition}
\newtheorem{defi}[theorem]{Definition}
\newtheorem{example}[theorem]{Example}
\theoremstyle{remark}
\newtheorem{rem}[theorem]{Remark}
\numberwithin{equation}{section}
\newcommand{\ci}{\ensuremath{C_\infty}}
\newcommand{\ai}{\ensuremath{A_\infty}}
\newcommand{\li}{\ensuremath{L_\infty}}
\newcommand{\csalg}[1]{\ensuremath{\widehat{S}\Sigma #1^*}}
\newcommand{\ctalg}[1]{\ensuremath{\widehat{T}\Sigma #1^*}}
\newcommand{\clalg}[1]{\ensuremath{\widehat{L}\Sigma #1^*}}
\newcommand{\col}[1]{\ensuremath{\prod_{i=1}^\infty (\Sigma {#1}^*)^{\cotimes i}}}
\newcommand{\deof}[2][]{\ensuremath{\Omega^1_\mathrm{#1}(#2)}}
\newcommand{\drof}[2][]{\ensuremath{DR^1_\mathrm{#1}(#2)}}
\newcommand{\drnof}[2][]{\ensuremath{\overline{DR}^1_\mathrm{#1}(#2)}}
\newcommand{\drzf}[2][]{\ensuremath{DR^0_\mathrm{#1}(#2)}}
\newcommand{\drtf}[2][]{\ensuremath{DR^2_\mathrm{#1}(#2)}}
\newcommand{\de}[2][]{\ensuremath{\Omega^\bullet_\mathrm{#1}(#2)}}
\newcommand{\dr}[2][]{\ensuremath{DR^\bullet_\mathrm{#1}(#2)}}
\newcommand{\mstr}[4][]{\ensuremath{\mathcal{#1#2}_{#3}(#4)}}
\newcommand{\mext}[3][]{\ensuremath{\mathcal{#1E}_{#2}(#3)}}
\newcommand{\mmor}[4][]{\ensuremath{\mathcal{#1M}_{#2}(#3 ; #4)}}
\newcommand{\mmext}[5][]{\ensuremath{\mathcal{#1EM}_{#2}(#3 : #4 , #5)}}
\newcommand{\Cbr}[2][\bullet]{\ensuremath{C_{#1}^\mathrm{Bar}(#2)}}
\newcommand{\Choch}[3][\bullet]{\ensuremath{C_{#1}^\mathrm{Hoch}(#2,#3)}}
\newcommand{\Hhoch}[3][\bullet]{\ensuremath{H_{#1}^\mathrm{Hoch}(#2,#3)}}
\newcommand{\clac}[3][\bullet]{\ensuremath{C_\mathrm{CE}^{#1}(#2,#3)}}
\newcommand{\hlac}[3][\bullet]{\ensuremath{H_\mathrm{CE}^{#1}(#2,#3)}}
\newcommand{\choch}[3][\bullet]{\ensuremath{C^{#1}_\mathrm{Hoch}(#2,#3)}}
\newcommand{\hhoch}[3][\bullet]{\ensuremath{H^{#1}_\mathrm{Hoch}(#2,#3)}}
\newcommand{\caq}[3][\bullet]{\ensuremath{C^{#1}_\mathrm{Harr}(#2,#3)}}
\newcommand{\haq}[3][\bullet]{\ensuremath{H^{#1}_\mathrm{Harr}(#2,#3)}}
\newcommand{\cbr}[2][\bullet]{\ensuremath{C^{#1}_\mathrm{Bar}(#2)}}
\newcommand{\hbr}[2][\bullet]{\ensuremath{H^{#1}_\mathrm{Bar}(#2)}}
\newcommand{\cclac}[2][\bullet]{\ensuremath{CC^{#1}_\mathrm{CE}(#2)}}
\newcommand{\hclac}[2][\bullet]{\ensuremath{HC^{#1}_\mathrm{CE}(#2)}}
\newcommand{\cchoch}[2][\bullet]{\ensuremath{CC^{#1}_\mathrm{Hoch}(#2)}}
\newcommand{\hchoch}[2][\bullet]{\ensuremath{HC^{#1}_\mathrm{Hoch}(#2)}}
\newcommand{\ccaq}[2][\bullet]{\ensuremath{CC^{#1}_\mathrm{Harr}(#2)}}
\newcommand{\hcaq}[2][\bullet]{\ensuremath{HC^{#1}_\mathrm{Harr}(#2)}}
\newcommand{\ctsygan}[2][\bullet\bullet]{\ensuremath{CC^{#1}_\mathrm{Tsygan}(#2)}}

\newcommand{\echoch}[4][\bullet]{\ensuremath{C[#4]^{#1}_\mathrm{Hoch}(#2,#3)}}
\newcommand{\ehhoch}[4][\bullet]{\ensuremath{H[#4]^{#1}_\mathrm{Hoch}(#2,#3)}}
\newcommand{\ecchoch}[3][\bullet]{\ensuremath{CC[#3]^{#1}_\mathrm{Hoch}(#2)}}
\newcommand{\ehchoch}[3][\bullet]{\ensuremath{HC[#3]^{#1}_\mathrm{Hoch}(#2)}}
\newcommand{\thiso}[1][]{\ensuremath{\Theta_\mathrm{#1}}}
\newcommand{\phiso}[1][]{\ensuremath{\Phi_\mathrm{#1}}}
\newcommand{\upiso}[1][]{\ensuremath{\Upsilon_\mathrm{#1}}}
\newcommand{\zetiso}[1][]{\ensuremath{\zeta_\mathrm{#1}}}
\newcommand{\cotimes}{\ensuremath{\hat{\otimes}}}
\newcommand{\gf}{\ensuremath{\mathbb{K}}}
\newcommand{\fpsa}{\ensuremath{\gf\langle\langle\tau,\boldsymbol{t}\rangle\rangle}}
\newcommand{\pvect}{\ensuremath{\mathcal{P}Mod_\gf}}
\newcommand{\inlim}[2]{\ensuremath{\varprojlim_{#1} #2}}
\newcommand{\dilim}[2]{\ensuremath{\varinjlim_{#1} #2}}
\newcommand{\innprod}[1]{\ensuremath{\langle -,- \rangle:#1 \otimes #1 \to \gf}}
\newcommand{\noproof}{\begin{flushright} \ensuremath{\square} \end{flushright}}
\DeclareMathOperator{\Map}{Map}
\DeclareMathOperator{\Hom}{Hom}
\DeclareMathOperator{\Der}{Der}
\DeclareMathOperator{\Coder}{Coder}
\DeclareMathOperator{\Aut}{Aut}
\DeclareMathOperator{\Obs}{Obs}
\DeclareMathOperator{\obs}{obs}
\DeclareMathOperator{\ord}{ord}
\DeclareMathOperator{\ad}{ad}
\DeclareMathOperator{\id}{id}

\begin{document}
\begin{abstract}
We study cohomology theories of strongly homotopy algebras, namely $A_\infty, C_\infty$ and $L_\infty$-algebras and establish the Hodge decomposition of Hochschild and cyclic cohomology of $C_\infty$-algebras thus generalising previous work by Loday and Gerstenhaber-Schack. These results are then used to show that a $C_\infty$-algebra with an invariant inner product on its cohomology can be uniquely extended to a symplectic $C_\infty$-algebra (an $\infty$-generalisation of a commutative Frobenius algebra introduced by Kontsevich). As another application, we show that the `string topology' operations (the loop product, the loop bracket and the string bracket) are homotopy invariant and can be defined on the homology or equivariant homology of an arbitrary Poincar\'e duality space.
\end{abstract}
\title{Homotopy algebras and noncommutative geometry}
\author{Alastair Hamilton \and Andrey Lazarev}
\address{Mathematics Department, Bristol University, Bristol, England. BS8 1TW.} \email{a.e.hamilton@bristol.ac.uk \and a.lazarev@bristol.ac.uk}

\thanks{The authors were partially supported by the EPSRC grants No. GR/S07148/01 and GR/R84276/01.}
\maketitle
\tableofcontents

\section{Introduction}

An $A_\infty$-algebra is a generalisation of an associative algebra introduced in \cite{staha2} for the purposes of studying $H$-spaces. It was originally defined via a system of higher multiplication maps satisfying a series of complicated relations. One can similarly define $C_\infty$ and $L_\infty$-algebras as $\infty$-generalisations of commutative and Lie algebras respectively.

More recently, $L_\infty$ and $A_\infty$-algebras have found applications in mathematical physics, particularly in string field theory and the theory of topological $\Sigma$-models, cf. \cite{zwiebach} and \cite{azksch}. In addition, $\infty$-algebras with an invariant inner product were introduced in \cite{kontfd} and \cite{kontsg} and were shown to have a close relation with graph homology and therefore to the intersection theory on the moduli spaces of complex curves and invariants of differentiable manifolds.  A short, informal introduction to graph homology 
is contained in \cite{voronov}, a more substantial account is in \cite{vogtmann}.

In this paper we give a detailed analysis of the cohomology theories associated to $A_\infty, L_\infty$ and $C_\infty$-algebras. It is likely that much of our work could be extended to algebras over a Koszul quadratic operad, however we see no advantage in working in this more general context since our most interesting applications are concerned with $C_\infty$-algebras and do not generalise to other operads. 

One of the important ideas that we wish to advocate is working with a geometric definition of an $\infty$-algebra, see Definition \ref{def_infstr}, as a formal supermanifold together with a homological vector field. This idea is not new, cf. for example \cite{azksch} or \cite{lazmod}, however we feel that it has not been used to its full potential. Apart from the obvious advantage of being concise, this definition allows one to engage the apparatus of noncommutative differential geometry which could be quite beneficial as we hope to demonstrate.

This paper contains two main sets of results. The first one is concerned with the Hodge decomposition of the Hochschild and cyclic Hochschild cohomology of a \ci-algebra and generalises the work of previous authors; cf. \cite{loday}, \cite{gersch} and \cite{natsch}. Our geometric approach allowed us to considerably simplify the combinatorics present in the cited sources while working in the more general context of $C_\infty$-algebras. We actually get something new even for the usual (strictly commutative) algebras, namely the Hodge decomposition of the Hochschild, cyclic Hochschild and bar cohomology in the \emph{nonunital} case.

The second set of results concerns $\infty$-algebras with an invariant inner product or \emph{symplectic} $\infty$-algebras. (It should be noted that some authors use the term `cyclic $\infty$-algebra', cf.  \cite{markl}, for example. We feel that the term `symplectic $\infty$-algebra' is most appropriate since the structure in question is simply a homological vector field preserving a given symplectic form.)  
Our main theorem here states that a unital $C_\infty$-algebra with the structure of a Frobenius algebra on its cohomology is weakly equivalent to a symplectic $C_\infty$-algebra which is uniquely determined up to homotopy. One of the consequences of that is that the cochain algebra of a rational Poincar\'e duality space is weakly equivalent to a symplectic $C_\infty$-algebra. We then use this fact to define `string topology' operations on the homology and equivariant homology of such spaces. Another consequence is that the graph homology class associated to a symplectic $C_\infty$-algebra only depends on the (homotopy class of) the underlying $C_\infty$-algebra together with the invariant inner product on its cohomology.

It would be interesting to find a physical explanation of the above results. It seems that the BRST/Batalin-Vilkovisky formalism (cf. \cite{batvi}) in quantum field theory is especially relevant in this connection. See also \cite{Stasheff}, \cite{stasheff}, \cite{schwarz} and \cite{Schwarz} for an overview of this method and its applications.

To prove the main theorem we develop an obstruction theory for lifting the  symplectic $C_n$-structures and $C_n$-morphisms which is of independent interest. The proof is constructed by combining the obstruction theory with the Hodge decomposition for cyclic cohomology.

The paper is organised as follows. Sections \ref{sec_ncgeom} and \ref{relations} introduce the basics of formal noncommutative geometry in the commutative, associative and Lie worlds. This is largely a detailed exposition of a part of Kontsevich's paper \cite{kontfd}. The language of topological modules and topological algebras is used here and throughout the paper and we collected the necessary facts in Appendix \ref{app_todual}. Section \ref{prerequisites} deals with the definitions and basic properties of $\infty$-algebras, emphasising the geometrical viewpoint. Here we consider $\mathbb{Z}$-graded $\infty$-algebras; however all the results carry over with obvious modifications to the $\mathbb{Z}/2$-graded framework used in the cited work of Kontsevich. In section \ref{minimal} we prove the analogue of Kadeishvili's minimality theorem for $C_\infty$-algebras.

In section \ref{sec_iachom} Chevalley-Eilenberg, Hochschild and Harrison cohomology theories for $\infty$-algebras are defined along with their cyclic counterparts and section \ref{sec_obstrc} is devoted to the obstruction theory for lifting the $C_n$-structures and $C_n$-morphisms. The Hodge decomposition for $C_\infty$-algebras is established in sections \ref{sec_hdgdec} and \ref{cyclichodge}; these two sections are the technical heart of the paper.

Section \ref{symplectic} discusses formal noncommutative \emph{symplectic} geometry; again main ideas here are due to Kontsevich, \cite{kontfd}; also in this section the notion of a \emph{symplectic $\infty$-algebra} is introduced. In section \ref{examples} we list some examples of symplectic $\infty$-algebras, particularly the Moore algebras which were introduced under this name in \cite{lazmod} and further studied in \cite{hamilt}. A version of these $A_\infty$-algebras leads to the Mumford-Morita-Miller classes in the cohomology spaces of complex algebraic curves, cf. \cite{kontfd} and \cite{igusa}.  Section \ref{sec_symobs} deals with the obstruction theory in the symplectic context. In section \ref{correspondence} we prove the theorem mentioned above on the relationship between symplectic and nonsymplectic  
$C_\infty$-algebras. Finally, section \ref{strings} discusses applications to string topology.

{\it Acknowledgement.} We are grateful to Jim Stasheff for many useful discussions and comments.

\subsection{Notation and conventions} \label{sec_notcon}

Throughout the paper our ground ring $\gf$ will be an evenly graded commutative ring containing the 
field $\mathbb{Q}$. For most applications it will be enough to assume that $\gf$ is a field or the rational numbers, however, our approach is designed to accommodate \emph{deformation theory} (which is not considered in the present paper) and to do this we need to allow ground rings which are not necessarily fields. \gf-algebras and \gf-modules will simply be called algebras and modules. We will assume that all of our \gf-modules are obtained from $k$-vector spaces by extension of scalars where $k=\mathbb{Q}$ or any other subfield of $\gf$. All of our tensors will be taken over the ground ring {\gf} unless stated otherwise.                                                             

Given a graded module V we define the tensor algebra $TV$ by
\[ TV:=\gf \oplus V \oplus V^{\otimes 2} \oplus \ldots \oplus V^{\otimes n} \oplus \ldots. \]
We define the symmetric algebra $SV$ as the quotient of $TV$ by the relation $x \otimes y = (-1)^{|x||y|}y \otimes x$. Finally we define the free Lie algebra $LV$ as the Lie subalgebra of the commutator algebra $TV$ which consists of linear combinations of Lie monomials.

We will use $~\widehat{ \ }~$ to denote completion. Given a profinite graded module $V$ we can define the completed versions of $SV$, $TV$ and $LV$. For instance the completed tensor algebra would be
\[ \widehat{T}V :=\prod_{i=0}^{\infty} V^{\cotimes i} \]
where $\cotimes$ denotes the completed tensor product. The completed symmetric algebra $\widehat{S}V$ is then the quotient of $\widehat{T}V$ by the usual relations and $\widehat{L}V$ is the Lie subalgebra of $\widehat{T}V$ consisting of all convergent (possibly uncountably) infinite linear combinations of Lie monomials. Further details on formal objects and constructions can be found in Appendix \ref{app_todual}. We will assume that all the formal commutative and associative $\gf$-algebras considered in the main text are augmented.

Given a \emph{formal} graded commutative, associative or Lie algebra $X$, the Lie algebra consisting of all \emph{continuous} derivations
\[ \xi : X \to X \]
is denoted by $\Der(X)$. In order to emphasise our use of geometrical ideas in this paper we will use the term `vector field' synonymously with `continuous derivation'. The group consisting of all invertible \emph{continuous} (commutative, associative or Lie) algebra homomorphisms
\[ \phi: X \to X \]
will be denoted by $\Aut(X)$. Again, in order to emphasise the geometrical approach we will call an `invertible continuous homomorphism' of formal graded commutative, associative or Lie algebras a `diffeomorphism'.

On many occasions we will want to deal with commutative, associative and Lie algebras simultaneously. When we do so we will often abuse the terminology; for instance calling a bimodule a `module' or calling a Lie algebra an `algebra'.

Given a formal graded associative algebra $X$, the module of commutators $[X,X]$ is defined as the module consisting of all convergent (possibly uncountably) infinite linear combinations of elements of the form,
\[ [x,y]:=x \cdot y - (-1)^{|x||y|}y \cdot x; \quad x,y \in X \]
We will often denote the module of commutators by $[-,-]$ when the context makes it clear what $X$ is.

Given a profinite graded module $V$ we can place a grading on $\widehat{T}V$ which is different from the grading which is naturally inherited from $V$. We say an element $x \in \widehat{T}V$ has homogeneous \emph{order} $n$ if $x \in V^{\cotimes n}$. This also defines a grading by order on $\widehat{S}V$ and $\widehat{L}V$ since they are quotients and submodules respectively of $\widehat{T}V$.

We say a continuous endomorphism (linear map) $f:\widehat{T}V \to \widehat{T}V$ has homogeneous \emph{order} $n$ if it takes any element of order $i$ to an element of order $i+n-1$. Any vector field $\xi \in \Der(\widehat{T}V)$ could be written in the form
\begin{equation} \label{eqn_vecfrm}
\xi = \xi_0 + \xi_1 + \xi_2 + \ldots + \xi_n + \ldots,
\end{equation}
where $\xi_i$ is a vector field of order $i$. We say that $\xi$ \emph{vanishes at zero} if $\xi_0 = 0$.

Likewise any diffeomorphism $\phi \in \Aut(\widehat{T}V)$ could be written in the form
\[ \phi = \phi_1 + \phi_2 + \phi_3 + \ldots + \phi_n + \ldots, \]
where $\phi_i$ is an endomorphism of order $i$. We call $\phi$ a \emph{pointed diffeomorphism} if $\phi_1 = \id$. Similarly, we could make the same definitions and observations if we were to replace $\widehat{T}V$ with $\widehat{S}V$ or $\widehat{L}V$ in the above.

We will denote the symmetric group on $n$ letters by $S_n$ and the cyclic group of order $n$ by $Z_n$. Given a module $M$ over a group $G$ the module of coinvariants will be denoted by $M_G$ and the module of invariants by $M^G$.

Given a graded profinite module $V$, its completed tensor power $V^{\cotimes n}$ has a continuous action of the cyclic group $Z_n$ on it which is the restriction of the canonical action of $S_n$ to the subgroup $Z_n \cong \langle (n\,n-1\ldots 2\,1) \rangle \subset S_n$. If we define the algebra $\Lambda$ as
\[ \Lambda:= \prod_{n=1}^\infty\mathbb{Z}[Z_n] \]
then the action of each $Z_n$ on $V^{\cotimes n}$ for $n\geq 1$ gives $\prod_{n=1}^\infty V^{\cotimes n}$ the structure of a left $\Lambda$-module.

Let $z_n$ denote the generator of $Z_n$ corresponding to the cycle $(n\,n-1\ldots 2\,1)$. We define $z \in \Lambda$ by the formula,
\[ z:=\sum_{n=1}^\infty z_n. \]
Let the norm operator $N_n \in \mathbb{Z}[Z_n]$ be the element given by the formula,
\[ N_n:=1+z_n+z_n^2 +\ldots + z_n^{n-1}. \]
We define $N \in \Lambda$ by the formula,
\[ N:= \sum_{n=1}^\infty N_n. \]
The operators $z$ and $N$ will be used regularly throughout the paper.

Our convention will be to always work with cohomologically graded objects and we consequently define the suspension of a graded module $V$ as $\Sigma V$ where $\Sigma V^i:=V^{i+1}$. We define the desuspension of $V$ as $\Sigma^{-1}V$ where $\Sigma^{-1}V^{i}:= V^{i-1}$. The term `differential graded algebra' will be abbreviated as `DGA'.

Given a graded module $V$, we denote the graded \gf-linear dual
\[ \Hom_\gf(V,\gf) \]
by $V^*$. If $V$ is free with basis $\{v_\alpha, \alpha \in I\}$, then we denote the dual basis of $V^*$ by
\[ \{v_\alpha^*, \alpha \in I\}. \]
For the sake of clarity, when we write $\Sigma V^*$ we mean the graded module $\Hom_\gf(\Sigma V,\gf)$.
\section{Formal Noncommutative Geometry} \label{sec_ncgeom}

In this section we will collect the definitions and fundamental facts about formal commutative and noncommutative geometry. Although much of this work can be done in the generality of working with an algebra over an operad, our applications will be concerned with the three particular theories of (non)commutative geometry corresponding to commutative, associative and Lie algebras and we will describe each theory in its own detail.

The starting point is to define the (non)commutative 1-forms in each of the three settings. From this we can construct a differential envelope of our commutative, associative or Lie algebra. We can then construct the de Rham complex from the differential envelope and introduce the contraction and Lie operators. We conclude the section with some basic facts from noncommutative geometry. The theory we consider here will be the formal version of the theory as described by Kontsevich in \cite{kontsg}. The reader can refer to Appendix \ref{app_todual} for a discussion of the formal objects and constructions that we use.

\begin{defi} \label{def_nconfm}
\
\begin{enumerate}
\item[(a)]
Let $A$ be a formal graded commutative algebra. Consider the module $A \cotimes (A/\gf)$ and write $x \cotimes y$ as $x\cotimes dy$. The module of commutative 1-forms \deof[Com]{A} is defined as the quotient of $A \cotimes (A/\gf)$ by the relations
\[ x \cotimes d(yz)=(-1)^{|z|(|y|+|x|)}zx\cotimes dy + xy\cotimes dz. \]
This is a formal left $A$-module via the action
\[ a \cdot x\cotimes dy:=ax\cotimes dy. \]
\item[(b)]
Let $A$ be a formal graded associative algebra. The module of noncommutative 1-forms \deof[Ass]{A} is defined as
\[ \deof[Ass]{A}:= A \cotimes (A/\gf). \]
Let us write $x\cotimes y$ as $x\cotimes dy$. \deof[Ass]{A} has the structure of a formal $A$-bimodule via the actions
\begin{align*}
& a \cdot x\cotimes dy:= ax \cotimes dy, \\
& x \cotimes dy \cdot a:= x \cotimes d(ya) - xy \cotimes da. 
\end{align*}
\item[(c)]
Let $g$ be a formal graded Lie algebra. Consider the module $\widehat{\mathcal{U}}(g) \cotimes g$ and write $x\cotimes y=x\cotimes dy$. Then the module of Lie 1-forms \deof[Lie]{g} is defined as the quotient of $\widehat{\mathcal{U}}(g) \cotimes g$ by the relations
\[ x \cotimes d[y,z]=xy \cotimes dz -(-1)^{|z||y|}xz \cotimes dy. \]
This is a formal left $g$-module via the action
\[ x \cdot (y \cotimes dz):= xy \cotimes dz. \]
\end{enumerate}
\end{defi}

\begin{rem}
When formulating theorems and proofs it will often be our policy to discuss these three theories simultaneously where this is practical. When we do this we will omit the subscripts $\mathrm{Com}$, $\mathrm{Ass}$ and $\mathrm{Lie}$. It should be understood that the reader should choose the appropriate subscript/construction depending on whether they wish to work with commutative, associative or Lie algebras.
\end{rem}

Let X be either a formal graded commutative, associative or Lie algebra. The following proposition shows that the module of 1-forms \deof{X} could be introduced as the unique $X$-module representing the functor $\Der^0(X,-)$ which sends a formal graded $X$-module $M$ to the module consisting of all continuous derivations $\xi:X \to M$ of degree zero. This proposition is the formal analogue of Lemma 5.5 of \cite{ginzsg}:

\begin{prop}
Let $X$ be either a formal graded commutative, associative or Lie algebra, then the map $d:X \to \Omega^1(X)$ is a derivation of degree zero. Given any formal graded $X$-module $M$ there is an isomorphism which is natural in both variables:
\begin{displaymath}
\begin{array}{ccc}
\Der^0(X,M) & \cong & \Hom^0_{X\mathrm{-module}}(\Omega^1(X),M), \\
\partial & \mapsto & [dx \mapsto \partial x]. \\
\end{array}
\end{displaymath}
\end{prop}
\noproof

This proposition could be summarised by the following diagram:
\[ \xymatrix{X \ar^{d}[rr] \ar_{\partial}[rd] && \Omega^1(X) \ar@{-->}^{\exists ! \phi}[ld] \\ & M} \]

\begin{rem}
Note that if $X$ is an associative algebra then the word `$X$-module' should be replaced with the word `$X$-bimodule'.
\end{rem}

Now we want to extend the 1-forms $\Omega^1(X)$ to the module of forms $\Omega^*(X)$ by forming a differential envelope of $X$.

\begin{defi} \label{def_ncfrms}
\
\begin{enumerate}
\item[(a)]
Let $A$ be a formal graded commutative algebra. The module of commutative forms \de[Com]{A} is defined as
\[ \de[Com]{A}:= \widehat{S}_A(\Sigma^{-1}\deof[Com]{A})=A \times \prod_{i=1}^\infty (\underbrace{\Sigma^{-1}\deof[Com]{A}\underset{A}{\cotimes}\ldots\underset{A}{\cotimes}\Sigma^{-1}\deof[Com]{A}}_{i \text{ factors}})_{S_i}. \]
Since \deof[Com]{A} is a module over the commutative algebra $A$, \de[Com]{A} has the structure of a formal graded commutative algebra whose multiplication is the standard multiplication on the completed symmetric algebra $\widehat{S}_A(\Sigma^{-1}\deof[Com]{A})$. The map $d:A \to \deof[Com]{A}$ lifts uniquely to give \de[Com]{A} the structure of a formal DGA.
\item[(b)]
Let $A$ be a formal graded associative algebra. The module of noncommutative forms \de[Ass]{A} is defined as
\[ \de[Ass]{A}:=\widehat{T}_A(\Sigma^{-1}\deof[Ass]{A})=A \times \prod_{i=1}^\infty \underbrace{\Sigma^{-1}\deof[Ass]{A} \underset{A}{\cotimes} \ldots \underset{A}{\cotimes} \Sigma^{-1}\deof[Ass]{A}}_{i \text{ factors}}. \]
Since \deof[Ass]{A} is an $A$-bimodule, \de[Ass]{A} has the structure of a formal associative algebra whose multiplication is the standard associative multiplication on the tensor algebra $\widehat{T}_A(\Sigma^{-1}\deof[Ass]{A})$. The map $d:A \to \deof[Ass]{A}$ lifts uniquely to give \de[Ass]{A} the structure of a formal DGA.
\item[(c)]
Let $g$ be a formal graded Lie algebra. The module of Lie forms \de[Lie]{g} is defined as
\[ \de[Lie]{g}:=g\ltimes \widehat{L}(\Sigma^{-1}\deof[Lie]{g}) \]
where the action of $g$ on $\widehat{L}(\Sigma^{-1}\deof[Lie]{g})$ is the restriction of the standard action of $g$ on $\widehat{T}(\Sigma^{-1}\deof[Lie]{g})$ to the Lie subalgebra of Lie monomials. The map $d:g \to \deof[Lie]{g}$ lifts uniquely to give \de[Lie]{g} the structure of a formal DGLA.
\end{enumerate}
\end{defi}

\begin{rem}
Since \deof{X} is a formal $X$-module, $\Sigma^{-1}\deof{X}$ is a formal $X$-module via the action
\[ x \cdot \Sigma^{-1}y:=(-1)^{|x|}\Sigma^{-1}(x \cdot y) \quad , \quad x \in X, y \in \deof{X}. \]
The additional signs appear because of the Koszul sign rule. This results in the map $d:\de{X} \to \de{X}$ being a graded derivation of degree one.
\end{rem}

\begin{rem}
The module of forms \de{X} inherits a grading from the graded modules $X$ and \deof{X}. In fact it has a bigrading. An element $x \in X$ has bidegree $(0,|x|)$ and an element $x \cotimes dy \in \deof{X}$ has bidegree $(1,|x|+|y|)$. This implicitly defines a bigrading on the whole of \de{X}. The map $d:\de{X} \to \de{X}$ has bidegree $(1,0)$ in this bigrading. The natural grading on \de{X} inherited from the graded modules $X$ and \deof{X} coincides with the total grading of this bigrading.
\end{rem}

The following proposition says that \de{X} could be uniquely defined as the differential envelope of $X$. This proposition is proved in the operadic and nonformal context in Proposition 5.6 of \cite{ginzsg}.

\begin{prop} \label{prop_difenv}
Let $X$ be either a formal graded commutative, associative or Lie algebra and let $M$ be either a formal differential graded commutative, associative or Lie algebra respectively. There is a natural adjunction isomorphism:
\[ \Hom_\mathrm{Alg}(X,M) \cong \Hom_\mathrm{DGA}(\de{X},M). \]
The adjunction isomorphism is given by the following diagram:
\begin{displaymath}
\xymatrix{X \ar@<-0.25ex>@{^{(}->}[rr] \ar^{\phi}[rd] && \de{X} \ar@{-->}^{\exists !\psi}[ld] \\ & M}
\end{displaymath}
\end{prop}
\noproof

\begin{rem} \label{rem_defunc}
It follows from Proposition \ref{prop_difenv} that $\Omega^\bullet$ is functorial. Suppose that both $X$ and $Y$ are either formal graded commutative, associative or Lie algebras. By Proposition \ref{prop_difenv}, any continuous algebra homomorphism $\phi:X \to Y$ lifts uniquely to a continuous differential graded algebra homomorphism $\phi^*:\de{X} \to \de{Y}$.
\end{rem}

Next we want to introduce contraction and Lie operators onto the module of forms \de{X}.

\begin{defi} \label{def_ncoper}
Let $X$ be either a formal graded commutative, associative or Lie algebra and let $\xi:X \to X$ be a vector field:
\begin{enumerate}
\item[(i)]
We can define a vector field $L_\xi: \de{X} \to \de{X}$ of bidegree $(0,|\xi|)$, called the Lie derivative, by the formula;
\begin{displaymath}
\begin{array}{lc}
L_\xi(x):=\xi(x), & x \in X; \\
L_\xi(dx):=(-1)^{|\xi|}d(\xi(x)), & x \in X. \\
\end{array}
\end{displaymath}
\item[(ii)]
We can define a vector field $i_\xi:\de{X} \to \de{X}$ of bidegree $(-1,|\xi|)$, called the contraction operator, by the formula;
\begin{displaymath}
\begin{array}{lc}
i_\xi(x):=0, & x \in X; \\
i_\xi(dx):=\xi(x), & x \in X. \\
\end{array}
\end{displaymath}
\end{enumerate}
\end{defi}

These operators satisfy certain important identities which are summarised by the following lemma.

\begin{lemma} \label{lem_schf}
Let both $X$ and $Y$ be either formal graded commutative, associative or Lie algebras, let $\xi:X \to X$ and $\gamma:X \to X$ be vector fields and let $\phi:X \to Y$ be a diffeomorphism, then we have the following identities:
\begin{displaymath}
\begin{array}{rl}
\textnormal{(i)} & L_\xi=[i_\xi,d]. \\
\textnormal{(ii)} & [L_\xi,i_\gamma]=i_{[\xi,\gamma]}. \\
\textnormal{(iii)} & L_{[\xi,\gamma]}=[L_\xi,L_\gamma]. \\
\textnormal{(iv)} & [i_\xi,i_\gamma]=0. \\
\textnormal{(v)} & [L_\xi,d]=0. \\
\textnormal{(vi)} & L_{\phi\xi\phi^{-1}}=\phi^* L_\xi \phi^{*-1}. \\
\textnormal{(vii)} & i_{\phi\xi\phi^{-1}}=\phi^* i_\xi \phi^{*-1}. \\
\end{array}
\end{displaymath}
\end{lemma}

\begin{proof}
Since \de{X} is generated by elements in $X$ and $d(X)$, we only need to establish the identities on these generators. Let us prove (i), the  Cartan homotopy formula. It should then be clear to the reader how to obtain the other identities.
\begin{align*}
& i_\xi d(x)=\xi(x), \quad di_\xi(x)=0; \\
& L_\xi(x)=\xi(x)=[i_\xi,d](x). 
\end{align*}
Similarly,
\begin{align*}
& i_\xi d(dx)=0, \quad di_\xi(dx)=d(\xi(x)); \\
& L_\xi(dx)=(-1)^{|\xi|}d(\xi(x))=[i_\xi,d](dx). 
\end{align*}
\end{proof}

We are now in a position to define the de Rham complex which will be the fundamental object in our applications of formal noncommutative geometry. For the purposes of working with commutative, associative or Lie algebras, one could use the definitions introduced by Kontsevich in \cite{kontsg}. These definitions were extended to the framework of operads by Getzler and Kapranov in \cite{getkap}. The functor $F: \mathcal{F}Alg \to \pvect$ (see Appendix \ref{app_todual} for definitions of the categories) is defined by the formula;
\[ F(X):=\widehat{S}^{\,2}(X)/(x\cotimes\mu(y,z)=\mu(x,y)\cotimes z), \]
where $X$ is a formal commutative, associative or Lie algebra and $\mu$ is the multiplication or Lie bracket respectively. This functor was originally introduced by Kontsevich in \cite{kontsg}. The de Rham complex is then defined as the result of applying the functor $F$ to the differential envelope \de{X}. The differential $d:\de{X} \to \de{X}$ induces the differential on the de Rham complex by applying the Leibniz rule, although in general the algebra structure will be lost and we will just have a complex.

\begin{defi} \label{def_derham}
\
\begin{enumerate}
\item[(a)]
Let $A$ be a formal graded commutative algebra. The de Rham complex \dr[Com]{A} is defined as
\[ \dr[Com]{A}:=\de[Com]{A}. \]
The differential is the same as the differential defined in Definition \ref{def_ncfrms}.
\item[(b)]
Let $A$ be a formal graded associative algebra. The de Rham complex \dr[Ass]{A} is defined as
\[ \dr[Ass]{A}:=\de[Ass]{A}/[\de[Ass]{A},\de[Ass]{A}]. \]
The differential is induced by the differential defined in Definition \ref{def_ncfrms}.
\item[(c)]
Let $g$ be a formal graded Lie algebra. The de Rham complex \dr[Lie]{g} is defined as the quotient of $\de[Lie]{g} \cotimes \de[Lie]{g}$ by the relations
\begin{displaymath}
\begin{array}{ll}
x \cotimes y = (-1)^{|x||y|}y \cotimes x; & x,y \in \de[Lie]{g}, \\
x \cotimes [y,z] = [x,y] \cotimes z; & x,y,z \in \de[Lie]{g}. \\
\end{array}
\end{displaymath}
The differential $d:\dr[Lie]{g} \to \dr[Lie]{g}$ is induced by the differential $d:\de[Lie]{g} \to \de[Lie]{g}$ by specifying
\[ d(x \cotimes y):=dx \cotimes y + (-1)^{|x|}x \cotimes dy; \quad x,y \in \de[Lie]{g}. \]
\end{enumerate}
\end{defi}

\begin{rem}
In the commutative and associative cases the definition in terms of the functor $F$ can be simplified, resulting in the definitions given above. The identification of $F(\de{X})$ with \dr{X} is given by
\begin{displaymath}
\begin{array}{ccc}
F(\de{X}) & \cong & \dr{X}, \\
x \cotimes y & \mapsto & x \cdot y. \\
\end{array}
\end{displaymath}
In the commutative case we see that \de[Com]{X} and \dr[Com]{X} actually coincide. In the Lie case we must resort to the slightly awkward definition of the de Rham complex given by the functor $F$.
\end{rem}

\begin{rem} \label{rem_ordgra}
Suppose we are given a profinite graded module $W$ and consider the commutative, associative or Lie algebras $\widehat{S}W$, $\widehat{T}W$ and $\widehat{L}W$. Recall from section \ref{sec_notcon} that these modules have a grading on them which we called the grading by order. An element $x \in \widehat{T}W$ has homogeneous order $n$ if $x \in W^{\cotimes n}$. We extend the grading by order on $\widehat{T}W$ to the whole of \dr[Ass]{\widehat{T}W} by stipulating that a 1-form $xdy$ has order $\ord(x)+\ord(y)-1$. This suffices to completely determine the grading by order on \dr[Ass]{\widehat{T}W}. For instance an $n$-form
\[ x_0 dx_1 \ldots dx_n \]
has order $\ord(x_0) + \ldots + \ord(x_n) - n$. We obtain a grading by order on \dr[Com]{\widehat{S}W} and \dr[Lie]{\widehat{L}W} in a similar and obvious manner.
\end{rem}

\begin{rem} \label{rem_drfunc}
It follows from Remark \ref{rem_defunc} that the de Rham complex construction is functorial. Suppose that both $X$ and $Y$ are either formal graded commutative, associative or Lie algebras and that $\phi:X \to Y$ is a continuous algebra homomorphism. The continuous differential graded algebra homomorphism $\phi^*:\de{X} \to \de{Y}$ induces a map between the de Rham complexes which, by an abuse of notation, we also denote by $\phi^*$;
\begin{equation} \label{eqn_drfunc}
\phi^*: \dr{X} \to \dr{Y}.
\end{equation}
For instance, in the Lie case the map $\phi^*:\dr[Lie]{g} \to \dr[Lie]{h}$ is given by;
\[ \phi^*(x \cotimes y):= \phi^*(x) \cotimes \phi^*(y); \quad x,y \in \de[Lie]{g}. \]
\end{rem}

The de Rham complex inherits a bigrading derived from the bigrading on \de{X}. Given a vector field $\xi:X \to X$ we can define Lie and contraction operators $L_\xi,i_\xi:\dr{X} \to \dr{X}$ of bidegrees $(0,|\xi|)$ and $(-1,|\xi|)$ respectively. These are induced by the Lie and contraction operators defined in Definition \ref{def_ncoper} and by an abuse of notation, we denote them by the same letters. It should be clear how $L_\xi$ and $i_\xi$ are defined except perhaps in the Lie case. In this case $L_\xi:\dr[Lie]{g} \to \dr[Lie]{g}$ is defined by specifying
\[ L_\xi(x \cotimes y):= L_\xi(x) \cotimes y + (-1)^{|x|}x \cotimes L_\xi(y); \quad x,y \in \de[Lie]{g}. \]
The operator $i_\xi:\dr[Lie]{g} \to \dr[Lie]{g}$ is defined in the same way. It is obvious that the operators $L_\xi,i_\xi:\dr{X} \to \dr{X}$ satisfy all the identities of Lemma \ref{lem_schf}.

We denote the $i$th fold suspension of the component of \dr{X} of bidegree $(i,\bullet)$ by $DR^i(X)$ and call these the $i$-forms. This definition gives us the identity;
\[ \dr{X} = DR^0(X) \times \Sigma^{-1}DR^1(X) \times \Sigma^{-2}DR^2(X) \times \ldots \times \Sigma^{-n}DR^n(X) \times \ldots. \]
The following lemma describing $DR^1(X)$ could be found in \cite{ginzsg}. There is also a description of $DR^0(X)$ in Proposition 4.9 of \cite{getkap} and in \cite{ginzsg}.

\begin{lemma} \label{lem_ofisom}
Let $W$ be a profinite graded module. We have the following isomorphisms of graded modules:
\begin{displaymath}
\begin{array}{rc}
\textnormal{(a)} & \thiso[Com] : W \cotimes \widehat{S}W \to \drof[Com]{\widehat{S}W}, \\
& \thiso[Com](x \cotimes y) := dx \cdot y. \\
\textnormal{(b)} & \thiso[Ass] : W \cotimes \widehat{T}W \to \drof[Ass]{\widehat{T}W}, \\
& \thiso[Ass](x \cotimes y) := dx \cdot y. \\
\textnormal{(c)} & \thiso[Lie] : W \cotimes \widehat{L}W \to \drof[Lie]{\widehat{L}W}, \\
& \thiso[Lie](x \cotimes y) := dx \cotimes y. \\
\end{array}
\end{displaymath}
\end{lemma}

\begin{proof}
Let us prove (c) since it is the most abstract. It should then be clear how to prove the other maps are isomorphisms. The 1-forms \drof[Lie]{\widehat{L}W} are generated by representatives of the form $x \cotimes y$ where $x \in \deof[Lie]{\widehat{L}W}$ and $y \in \widehat{L}W$. Using the Leibniz rule
\[ d[a,b]= (-1)^{|a|}[a,db] -(-1)^{(|a|+1)|b|}[b,da]; \quad a,b \in \widehat{L}W \]
we conclude that \deof[Lie]{\widehat{L}W} is generated by elements of the form $u \cdot dw\cotimes y$ where $u \in \widehat{\mathcal{U}}(\widehat{L}W)$, $w \in W$ and $y \in \widehat{L}W$. Now using the relation we imposed on \dr[Lie]{\widehat{L}W} in Definition \ref{def_derham},
\[ [a,b] \cotimes c = a \cotimes [b,c]; \quad a,b,c \in \de[Lie]{\widehat{L}W}; \]
we assert that \drof[Lie]{\widehat{L}W} is generated by representatives of the form $dw \cotimes y$ where $w \in W$ and $y \in \widehat{L}W$. It is now clear that \thiso[Lie] is an isomorphism.
\end{proof}

\begin{rem}
Note that since $W \cotimes \widehat{T}W = \prod_{i=1}^\infty W^{\cotimes i}$ we may also regard \thiso[Ass] as a map
\[ \thiso[Ass] : \prod_{i=1}^\infty W^{\cotimes i} \to \drof[Ass]{\widehat{T}W}. \]
\end{rem}

Given either a formal commutative, associative or Lie algebra $X$, we define the de Rham cohomology of $X$ to be the cohomology of the complex \dr{X}. We define $H^i(\dr{X})$ as
\[ H^i(\dr{X}):=\frac{\{ x \in DR^i(X):dx=0 \}}{d(DR^{i-1}(X))}. \]
The following result is the (non)commutative analogue of the Poincar\'e lemma and can be found in \cite{kontsg}.

\begin{lemma} \label{lem_poinca}
Let $W$ be a profinite graded module:
\begin{enumerate}
\item[(a)]
\begin{displaymath}
\begin{array}{lcll}
H^i(\dr[Com]{\widehat{S}W}) & = & 0 & i \geq 1, \\
H^0(\dr[Com]{\widehat{S}W}) & = & \gf. \\
\end{array}
\end{displaymath}
\item[(b)]
\begin{displaymath}
\begin{array}{lcll}
H^i(\dr[Ass]{\widehat{T}W}) & = & 0 & i \geq 1, \\
H^0(\dr[Ass]{\widehat{T}W}) & = & \gf. \\
\end{array}
\end{displaymath}
\item[(c)]
\[ H^\bullet(\dr[Lie]{\widehat{L}W}) = 0. \]
\end{enumerate}
\end{lemma}

\begin{proof}
Let us prove (b). It will then be clear that the same argument will work to establish the other cases. Choose a topological basis $\{x_i\}_{i \in I}$ for the profinite module $W$ and consider Euler's vector field $\xi:\widehat{T}W \to \widehat{T}W$ of degree zero and order one:
\[ \xi:= \sum_{i \in I}x_i\partial_{x_i}. \]
Since we are working in characteristic zero, $L_\xi:\dr[Ass]{\widehat{T}W} \to \dr[Ass]{\widehat{T}W}$ establishes a one-to-one correspondence on every component except on the component
\[ \gf \subset \widehat{T}W/[\widehat{T}W,\widehat{T}W]=\drzf[Ass]{\widehat{T}W} \]
of bidegree $(0,0)$ and order $0$ which it maps to zero, however, by the  Cartan homotopy formula of Lemma \ref{lem_schf} part (i) we know that $L_\xi$ is nullhomotopic whence the result.
\end{proof}

\section{Relations between Commutative, Associative and Lie Geometries}\label{relations}

In this section we will discuss some aspects of noncommutative geometry which describe the relationships between the commutative, associative and Lie versions of the theory described in section \ref{sec_ncgeom}. We will only concern ourselves with describing this relationship for the pro-free algebras $\widehat{S}W$, $\widehat{T}W$ and $\widehat{L}W$. One of the main purposes of this section will be to lay the groundwork for sections \ref{sec_iachom}, \ref{sec_hdgdec}, \ref{cyclichodge} and \ref{symplectic} and also to expand upon a theorem of Kontsevich \cite{kontsg}. This theorem will allow us to consider the closed commutative, associative or Lie 2-forms of the pro-free algebras $\widehat{S}U^*$, $\widehat{T}U^*$ or $\widehat{L}U^*$ respectively, as linear maps $\alpha:TU \to \gf$.

We first of all establish some fundamental identities for the operators $N$ and $z$ defined in section \ref{sec_notcon}. We then introduce a series of maps which we use to establish some relationships between the complexes \dr[Com]{\widehat{S}W}, \dr[Ass]{\widehat{T}W} and \dr[Lie]{\widehat{L}W}. We use these relationships to prove Theorem \ref{thm_tfmaps} and we use this theorem to interpret the closed 2-forms as linear maps. In particular we identify the 2-forms giving rise to symmetric bilinear forms.

\begin{lemma} \label{lem_nrmdif}
Let $W$ be a graded profinite module:
\begin{enumerate}
\item[(i)]
The following identity holds:
\begin{equation} \label{eqn_nrmdifdummy}
\widehat{T}W/[\widehat{T}W,\widehat{T}W] = \prod_{i=0}^\infty (W^{\cotimes i})_{Z_i}.
\end{equation}
\item[(ii)]
The following diagram commutes:
\begin{displaymath}
\xymatrix{ W \cotimes \widehat{T}W \ar^{\thiso[Ass]}_{\cong}[r] \ar@{=}[dd] & \drof[Ass]{\widehat{T}W} & \drzf[Ass]{\widehat{T}W} \ar_{d}[l] \ar@{=}[d] \\ & \widehat{T}W \ar@{->>}[r] & \widehat{T}W/[\widehat{T}W,\widehat{T}W] \\ \prod_{i=1}^\infty W^{\cotimes i} & \prod_{i=1}^\infty W^{\cotimes i} \ar_{N}[l] \ar@{^{(}->}[u] \ar@{->>}[r] & \prod_{i=1}^\infty (W^{\cotimes i})_{Z_i} \ar@{^{(}->}[u] \\ }
\end{displaymath}
\end{enumerate}
\end{lemma}

\begin{proof}
Let $x:=x_1\cotimes\ldots\cotimes x_n \in W^{\cotimes n} \subset \widehat{T}W$:
\begin{enumerate}
\item[(i)]
For $1\leq i \leq n-1$ we calculate;
\begin{equation} \label{eqn_cyccom}
\begin{split}
(1-z_n^i)\cdot x & = x_1\cotimes\ldots\cotimes x_n - (-1)^{(|x_1|+\ldots+|x_i|)(|x_{i+1}|+\ldots+|x_n|)}x_{i+1}\cotimes\ldots\cotimes x_n\cotimes x_1\cotimes \ldots x_i, \\
& = [x_1\cotimes\ldots\cotimes x_i,x_{i+1}\cotimes\ldots\cotimes x_n]. \\
\end{split}
\end{equation}
This calculation suffices to establish equation \eqref{eqn_nrmdifdummy}.

\item[(ii)]
\begin{equation} \label{eqn_nrmdif}
\begin{split}
d(x) &= \sum_{i=1}^n (-1)^{|x_1|+\ldots+|x_{i-1}|} x_1\cotimes\ldots\cotimes x_{i-1}\cdot dx_i \cdot x_{i+1}\cotimes\ldots\cotimes x_n, \\
& = \sum_{i=1}^n (-1)^{(|x_1|+\ldots+|x_{i-1}|)(|x_i|+\ldots+|x_n|)} dx_i \cdot x_{i+1}\cotimes\ldots\cotimes x_n \cotimes x_1 \cotimes \ldots \cotimes x_{i-1} \mod [-,-]. \\
N_n\cdot x & = \sum_{i=0}^{n-1} (-1)^{(|x_1|+\ldots+|x_i|)(|x_{i+1}|+\ldots+|x_n|)}x_{i+1}\cotimes\ldots\cotimes x_n\cotimes x_1 \cotimes \ldots \cotimes x_i. \\
\end{split}
\end{equation}
From this calculation we see that $\thiso[Ass](N_n \cdot x) = d(x)$ and therefore the diagram is commutative.
\end{enumerate}
\end{proof}

Our next task will be to identify the closed 2-forms of \dr[Com]{\widehat{S}U^*}, \dr[Ass]{\widehat{T}U^*} and \dr[Lie]{\widehat{L}U^*} with a submodule of the module consisting of all linear maps $\alpha:TU \to \gf$. In order to do this we will need to introduce maps between the de Rham complexes and establish some basic properties of them.

Let $W$ be a graded profinite module. The completed tensor algebra $\widehat{T}W$ is a Lie algebra under the commutator and the differential envelope \de[Ass]{\widehat{T}W} is a differential graded Lie algebra under the commutator. There is a series of canonical Lie algebra inclusions;
\[ \widehat{L}W \hookrightarrow \widehat{T}W \hookrightarrow \de[Ass]{\widehat{T}W}. \]
By Proposition \ref{prop_difenv} the composition of these inclusions extends to a map of differential graded Lie algebras
\[ l':\de[Lie]{\widehat{L}W} \to \de[Ass]{\widehat{T}W}. \]
One can easily check that the following formula holds for all $a,b,c \in \de[Ass]{\widehat{T}W}$:
\begin{equation} \label{eqn_modcom}
a[b,c]=[a,b]c \mod [-,-].
\end{equation}
This allows us to define a map $l: \dr[Lie]{\widehat{L}W} \to \dr[Ass]{\widehat{T}W}$ which is induced by the map $l'$ and given by
\begin{equation} \label{eqn_liasmp}
l(x\cotimes y):=l'(x)l'(y); \quad x,y \in \de[Lie]{\widehat{L}W}.
\end{equation}
It follows from equation \eqref{eqn_modcom} that this map is well defined modulo the relations of Definition \ref{def_derham} part (c).

We can define another map $p:\dr[Ass]{\widehat{T}W} \to \dr[Com]{\widehat{S}W}$ as follows: The commutative algebra $\widehat{S}W$ is obviously an associative algebra and the canonical projection $\widehat{T}W \twoheadrightarrow \widehat{S}W$ is a map of associative algebras. By Proposition \ref{prop_difenv} this map extends to a map of DGAs
\[ p':\de[Ass]{\widehat{T}W} \to \de[Com]{\widehat{S}W}. \]
Clearly this map is zero on commutators and hence induces a map
\begin{equation} \label{eqn_ascomp}
p: \dr[Ass]{\widehat{T}W} \to \dr[Com]{\widehat{S}W}.
\end{equation}

Next we need to introduce a map $i:\widehat{S}W \to \prod_{i=0}^\infty (W^{\cotimes i})^{S_i}$. This is defined as the canonical map identifying coinvariants with invariants and is given by
\begin{equation} \label{eqn_symiso}
i(x_1\cotimes\ldots\cotimes x_n):= \sum_{\sigma \in S_n} \sigma\cdot x_1\cotimes\ldots\cotimes x_n.
\end{equation}

Finally, we can describe a map $j:\drof[Com]{\widehat{S}W} \to \drof[Ass]{\widehat{T}W}$. This is defined by the following commutative diagram which uses the isomorphisms defined by Lemma \ref{lem_ofisom} and the map $i$ defined above by equation \eqref{eqn_symiso}:
\begin{equation} \label{fig_caofmp}
\xymatrix{ W \cotimes \widehat{S}W \ar^{\thiso[Com]}_{\cong}[r] \ar^{1 \cotimes i}[d] & \drof[Com]{\widehat{S}W} \ar^{j}[d] \\ W \cotimes \widehat{T}W \ar^{\thiso[Ass]}_{\cong}[r] & \drof[Ass]{\widehat{T}W} \\ }
\end{equation}

The map $l:\dr[Lie]{\widehat{L}W} \to \dr[Ass]{\widehat{T}W}$ defined by equation \eqref{eqn_liasmp} can be lifted to a unique map
\[ \frac{\drof[Lie]{\widehat{L}W}}{d(\drzf[Lie]{\widehat{L}W})} \to \frac{\drof[Ass]{\widehat{T}W}}{d(\drzf[Ass]{\widehat{T}W})}. \]
Similarly, the map $p:\drof[Ass]{\widehat{T}W}\to\drof[Com]{\widehat{S}W}$  of \eqref{eqn_ascomp} can be lifted to a unique map
\[ \frac{\drof[Ass]{\widehat{T}W}}{d(\drzf[Ass]{\widehat{T}W})} \to \frac{\drof[Com]{\widehat{S}W}}{d(\drzf[Com]{\widehat{S}W})}. \]
Abusing notation, we will denote the induced maps by the same symbols $l$ and $p$.

\begin{lemma} \label{lem_clsdtf}
Let $W$ be a graded profinite module:
\begin{enumerate}
\item[(i)]
The map $l:\drof[Lie]{\widehat{L}W} \to \drof[Ass]{\widehat{T}W}$ is injective.
\item[(ii)]
The map
\[ l:\frac{\drof[Lie]{\widehat{L}W}}{d(\drzf[Lie]{\widehat{L}W})} \to \frac{\drof[Ass]{\widehat{T}W}}{d(\drzf[Ass]{\widehat{T}W})} \]
is also injective.  
\item[(iii)]
The map $j:\drof[Com]{\widehat{S}W}\to\drof[Ass]{\widehat{T}W}$ defined by diagram \eqref{fig_caofmp} is injective and the map $p:\drof[Ass]{\widehat{T}W}\to\drof[Com]{\widehat{S}W}$  is surjective.
\item[(iv)]
The map $j:\drof[Com]{\widehat{S}W}\to\drof[Ass]{\widehat{T}W}$ can be lifted to a unique map
\[ \frac{\drof[Com]{\widehat{S}W}}{d(\drzf[Com]{\widehat{S}W})} \to \frac{\drof[Ass]{\widehat{T}W}}{d(\drzf[Ass]{\widehat{T}W})} \]
which will be denoted by the same letter $j$ by the customary abuse of notation.
\item[(v)]
The map
\[ j:\frac{\drof[Com]{\widehat{S}W}}{d(\drzf[Com]{\widehat{S}W})} \to \frac{\drof[Ass]{\widehat{T}W}}{d(\drzf[Ass]{\widehat{T}W})} \]
is injective and the map
\[ p:\frac{\drof[Ass]{\widehat{T}W}}{d(\drzf[Ass]{\widehat{T}W})} \to \frac{\drof[Com]{\widehat{S}W}}{d(\drzf[Com]{\widehat{S}W})} \]
is surjective.
\end{enumerate}
All told, we have the following diagram of injections and surjections:
\begin{displaymath}
\xymatrix{ \drof[Com]{\widehat{S}W} \ar@{->>}[r] \ar@<-1ex>@{^{(}->}_{j}[d] & \frac{\drof[Com]{\widehat{S}W}}{d(\drzf[Com]{\widehat{S}W})} \ar@<-1ex>@{^{(}->}_{j}[d] \\  \drof[Ass]{\widehat{T}W} \ar@{->>}[r] \ar@<-1ex>@{->>}_{p}[u] & \frac{\drof[Ass]{\widehat{T}W}}{d(\drzf[Ass]{\widehat{T}W})} \ar@<-1ex>@{->>}_{p}[u] \\ \drof[Lie]{\widehat{L}W} \ar@{->>}[r] \ar@{^{(}->}^{l}[u] & \frac{\drof[Lie]{\widehat{L}W}}{d(\drzf[Lie]{\widehat{L}W})} \ar@{^{(}->}^{l}[u] \\ }
\end{displaymath}
\end{lemma}

\begin{proof}
From the identities;
\begin{enumerate}
\item[(a)] $\prod_{i=1}^\infty W^{\cotimes i} = W \cotimes \widehat{T}W$,
\item[(b)] $\drzf[Ass]{\widehat{T}W} = \widehat{T}W/[\widehat{T}W,\widehat{T}W]$,
\item[(c)] $\drzf[Lie]{\widehat{L}W} = \widehat{L}W \cotimes \widehat{L}W / ([x,y]\cotimes z = x \cotimes [y,z])$;
\end{enumerate}
it follows that there are maps;
\[ \pi_\mathrm{Ass} : W \cotimes \widehat{T}W \to \drzf[Ass]{\widehat{T}W}, \]
\[ \pi_\mathrm{Lie} : W \cotimes \widehat{L}W \to \drzf[Lie]{\widehat{L}W}; \]
which are just the canonically defined projection maps. Let $u:\widehat{L}W \to \widehat{T}W$ be the canonical inclusion of the Lie subalgebra $\widehat{L}W$ into $\widehat{T}W$. These maps fit together with the maps \thiso[Ass] and \thiso[Lie] defined by Lemma \ref{lem_ofisom} and the map $l$ defined by equation \eqref{eqn_liasmp} to form the following commutative diagram:
\begin{equation} \label{fig_clsdtf}
\xymatrix{ \drof[Ass]{\widehat{T}W} & W \cotimes \widehat{T}W \ar^{\pi_\mathrm{Ass}}[r] \ar_{\thiso[Ass]}^{\cong}[l] & \drzf[Ass]{\widehat{T}W} \\ \drof[Lie]{\widehat{L}W} \ar^{l}[u] & W \cotimes \widehat{L}W \ar^{1 \cotimes u}[u] \ar^{\pi_\mathrm{Lie}}[r] \ar_{\thiso[Lie]}^{\cong}[l] & \drzf[Lie]{\widehat{L}W} \ar^{l}[u] \\ }
\end{equation}
\begin{enumerate}
\item[(i)]
Since $u$ is injective we conclude from diagram \eqref{fig_clsdtf} that the map
\[ l:\drof[Lie]{\widehat{L}W} \to \drof[Ass]{\widehat{T}W} \]
is injective.
\item[(ii)]
Let $x \in \drof[Lie]{\widehat{L}W}$ be a 1-form of order $n$ and suppose that $l(x) \in d(\drzf[Ass]{\widehat{T}W})$. By Lemma \ref{lem_nrmdif} part (ii), \thiso[Ass] satisfies the identity,
\begin{equation} \label{eqn_clsdtfdummyb}
\thiso[Ass] \circ N = d \circ \pi_\mathrm{Ass}
\end{equation}
and therefore identifies $d(\drzf[Ass]{\widehat{T}W})$ with the module of invariants $\prod_{i=1}^\infty (W^{\cotimes i})^{Z_i}$. Using the commutativity of diagram \eqref{fig_clsdtf} we perform the following calculation:
\begin{displaymath}
\begin{split}
\thiso[Ass]^{-1}l(x) & = \frac{1}{n+1} N_{n+1} \cdot \thiso[Ass]^{-1}l(x), \\
l(x) & = \frac{1}{n+1} \thiso[Ass] N_{n+1} \thiso[Ass]^{-1}l(x),  \\
& = \frac{1}{n+1} ld\pi_\mathrm{Lie}\thiso[Lie]^{-1}(x). \\
\end{split}
\end{displaymath}
From part (i) of this lemma we conclude that $x=\frac{1}{n+1}d\pi_\mathrm{Lie}\thiso[Lie]^{-1}(x)$ and hence the map
\[ l: \frac{\drof[Lie]{\widehat{L}W}}{d(\drzf[Lie]{\widehat{L}W})} \to \frac{\drof[Ass]{\widehat{T}W}}{d(\drzf[Ass]{\widehat{T}W})} \]
is injective.
\item[(iii)]
One can easily verify the following equation from which the assertion follows immediately:
\begin{equation} \label{eqn_clsdtfdummya}
pj(dx_0 \cdot x_1 \cotimes \ldots \cotimes x_n) = n!dx_0 \cdot x_1 \cotimes \ldots \cotimes x_n.
\end{equation}
\item[(iv)]
Let $x:=x_1 \cotimes \ldots \cotimes x_{n+1} \in \widehat{T}W$ and let $\bar{x}$ denote its image in the quotient $\widehat{S}W =\drzf[Com]{\widehat{S}W}$, then
\begin{displaymath}
\begin{split}
\thiso[Ass]^{-1}jd(\bar{x}) & = \sum_{\sigma \in S_n}\sum_{\tau \in Z_{n+1}} \sigma\tau \cdot x, \\
& = \sum_{\sigma \in S_{n+1}} \sigma \cdot x, \\
& = \sum_{\tau \in Z_{n+1}}\sum_{\sigma \in S_n} \tau\sigma \cdot x. \\
\end{split}
\end{displaymath}
Using equation \eqref{eqn_clsdtfdummyb} we obtain;
\[ jd(\bar{x}) = d\pi_\mathrm{Ass}\left(\sum_{\sigma \in S_n} \sigma \cdot x\right), \]
hence the map defined by diagram \eqref{fig_caofmp} can be lifted to a map
\[ j:\frac{\drof[Com]{\widehat{S}W}}{d(\drzf[Com]{\widehat{S}W})} \to \frac{\drof[Ass]{\widehat{T}W}}{d(\drzf[Ass]{\widehat{T}W})}. \]
\item[(v)]
This is an immediate consequence of equation \eqref{eqn_clsdtfdummya}.
\end{enumerate}
\end{proof}

\begin{rem} \label{rem_vecseq}
Let $W$ be a graded profinite module. We will now describe a sequence of Lie algebra homomorphisms;
\[ \Der(\widehat{L}W) \hookrightarrow \Der(\widehat{T}W) \twoheadrightarrow \Der(\widehat{S}W). \]

Any vector field $\xi \in \Der(\widehat{T}W)$ is specified uniquely by its restriction $\xi:W \to \widehat{T}W$ (cf. Proposition \ref{prop_dergen}) and likewise for $\widehat{L}W$ and $\widehat{S}W$. It follows that any vector field $\xi \in \Der(\widehat{L}W)$ can be extended uniquely to a vector field
\[ \xi: \widehat{T}W \to \widehat{T}W. \]
In fact, this vector field is a Hopf algebra derivation (i.e. a derivation and a coderivation) when $\widehat{T}W$ is equipped with the cocommutative comultiplication (\emph{shuffle coproduct}) defined in Remark \ref{rem_ccmult}. Furthermore all continuous Hopf algebra derivations are obtained from $\Der(\widehat{L}W)$ in this way. This is because the Lie subalgebra of primitive elements of this Hopf algebra coincides with the Lie subalgebra $\widehat{L}W$. This establishes a one-to-one correspondence between $\Der(\widehat{L}W)$ and continuous Hopf algebra derivations on $\widehat{T}W$.

Any map $\xi:W \to \widehat{T}W$ gives rise to a map $\tilde{\xi}:W \to \widehat{S}W$ simply by composing it with the canonical projection $\widehat{T}W \to \widehat{S}W$. It follows that any vector field $\xi \in \Der(\widehat{T}W)$ can be lifted to a unique vector field
\[ \tilde{\xi}: \widehat{S}W \to \widehat{S}W \]
and that all the vector fields in $\Der(\widehat{S}W)$ are obtained from $\Der(\widehat{T}W)$ in this way.
\end{rem}

We will now formulate a lemma, the proof of which is a simple check, which describes how this sequence of Lie algebra homomorphisms interact with our maps $l$ and $p$:

\begin{lemma} \label{lem_plcomm}
Let $W$ be a graded profinite module:
\begin{enumerate}
\item[(i)]
Given any vector field $\xi:\widehat{L}W \to \widehat{L}W$, the following diagrams commute:
\begin{displaymath}
\xymatrix{ \dr[Ass]{\widehat{T}W} \ar^{i_\xi}[r] & \dr[Ass]{\widehat{T}W} \\ \dr[Lie]{\widehat{L}W} \ar^{i_\xi}[r] \ar@{^{(}->}^{l}[u] & \dr[Lie]{\widehat{L}W} \ar@{^{(}->}^{l}[u] \\ }
\qquad
\xymatrix{ \dr[Ass]{\widehat{T}W} \ar^{L_\xi}[r] & \dr[Ass]{\widehat{T}W} \\ \dr[Lie]{\widehat{L}W} \ar^{L_\xi}[r] \ar@{^{(}->}^{l}[u] & \dr[Lie]{\widehat{L}W} \ar@{^{(}->}^{l}[u] \\ }
\end{displaymath}
\item[(ii)]
Given any vector field $\xi:\widehat{T}W \to \widehat{T}W$, the following diagrams commute:
\begin{displaymath}
\xymatrix{ \dr[Com]{\widehat{S}W} \ar^{i_{\tilde{\xi}}}[r] & \dr[Com]{\widehat{S}W} \\ \dr[Ass]{\widehat{T}W} \ar^{i_\xi}[r] \ar@{->>}^{p}[u] & \dr[Ass]{\widehat{T}W} \ar@{->>}^{p}[u] \\ }
\qquad
\xymatrix{ \dr[Com]{\widehat{S}W} \ar^{L_{\tilde{\xi}}}[r] & \dr[Com]{\widehat{S}W} \\ \dr[Ass]{\widehat{T}W} \ar^{L_\xi}[r] \ar@{->>}^{p}[u] & \dr[Ass]{\widehat{T}W} \ar@{->>}^{p}[u] \\ }
\end{displaymath}
\end{enumerate}
\end{lemma}
\noproof

Using Lemma \ref{lem_clsdtf} we can prove the following interesting theorem. This result was proved in \cite{kontsg} for closed associative 2-forms and stated without a proof for closed Lie 2-forms. We will describe the analogous result in the commutative case as well. Given a graded profinite module $W$ we can consider the module of commutators
\[ [W,\textstyle{\prod_{i=0}^\infty (W^{\cotimes i})^{S_i}}] \subset [\widehat{T}W,\widehat{T}W] \subset \widehat{T}W. \]

\begin{theorem} \label{thm_tfmaps}
Let $W$ be a graded profinite module. There are isomorphisms:
\begin{enumerate}
\item[(i)]
\[ \zetiso[Ass]:\{ \omega \in \drtf[Ass]{\widehat{T}W} : d\omega = 0 \} \to [\widehat{T}W,\widehat{T}W]. \]
\item[(ii)]
\[ \zetiso[Com]:\{ \omega \in \drtf[Com]{\widehat{S}W} : d\omega = 0 \} \to [W,\textstyle{\prod_{i=0}^\infty (W^{\cotimes i})^{S_i}}]. \]
\item[(iii)]
\[ \zetiso[Lie]:\{ \omega \in \drtf[Lie]{\widehat{L}W} : d\omega = 0 \} \to [\widehat{L}W,\widehat{L}W]. \]
\end{enumerate}
\end{theorem}

\begin{proof}
Firstly if $X$ is one of the three algebras $\widehat{S}W$, $\widehat{T}W$ or $\widehat{L}W$ then by Lemma \ref{lem_poinca} there is an isomorphism of modules,
\[ f: \drof{X}/d(\drzf{X}) \to \{ \omega \in \drtf{X} : d\omega = 0 \}; \]
where $f$ is defined by the formula $f(x):=dx$. We will denote this isomorphism by $f_\mathrm{Com}$, $f_\mathrm{Ass}$ or $f_\mathrm{Lie}$ respectively.
\begin{enumerate}
\item[(i)]
First of all note that by Lemma \ref{lem_nrmdif}, the isomorphism \thiso[Ass] defined in Lemma \ref{lem_ofisom} induces a map,
\[ \bar{\Theta}_\mathrm{Ass}: \frac{\prod_{i=1}^\infty W^{\cotimes i}}{N\cdot\prod_{i=1}^\infty W^{\cotimes i}} \to \frac{\drof[Ass]{\widehat{T}W}}{d(\drzf[Ass]{\widehat{T}W})} \]
which is also an isomorphism.

Equation \eqref{eqn_cyccom} gives us the identity
\[ (1-z)\cdot\prod_{i=1}^\infty W^{\cotimes i} = [\widehat{T}W,\widehat{T}W]. \]
This allows us to define a module isomorphism
\[ g: \frac{\drof[Ass]{\widehat{T}W}}{d(\drzf[Ass]{\widehat{T}W})} \to [\widehat{T}W,\widehat{T}W] \]
by the formula $g(x):=(1-z)\cdot\bar{\Theta}_\mathrm{Ass}^{-1}(x)$. We then define the isomorphism \zetiso[Ass] by the formula
\[ \zetiso[Ass]:= g \circ f_\mathrm{Ass}^{-1}. \]
\item[(ii)]
By Lemma \ref{lem_clsdtf} there is an injection
\[ j: \frac{\drof[Com]{\widehat{S}W}}{d(\drzf[Com]{\widehat{S}W})} \hookrightarrow \frac{\drof[Ass]{\widehat{T}W}}{d(\drzf[Ass]{\widehat{T}W})}. \]
Composing $j$ with $g$ yields an injection $g\circ j$. It requires a simple check using the definition of $j$ and equation \eqref{eqn_cyccom} to see that the image of $g\circ j$ is $[W,\prod_{i=1}^\infty (W^{\cotimes i})^{S_i}]$. It follows that the map
\[ \zetiso[Com]:\{ \omega \in \drtf[Com]{\widehat{S}W} : d\omega = 0 \} \to [W,\textstyle{\prod_{i=0}^\infty (W^{\cotimes i})^{S_i}}] \]
defined by the formula $\zetiso[Com](x):=gjf_\mathrm{Com}^{-1}(x)$ is an isomorphism.
\item[(iii)]
By Lemma \ref{lem_clsdtf} there is an injection
\[l: \frac{\drof[Lie]{\widehat{L}W}}{d(\drzf[Lie]{\widehat{L}W})} \hookrightarrow \frac{\drof[Ass]{\widehat{T}W}}{d(\drzf[Ass]{\widehat{T}W})}. \]
Composing $l$ with $g$ yields an injection $g\circ l$. Again it requires a simple check using the definitions and equation \eqref{eqn_cyccom} to see that the image of $g \circ l$ is $[W,\widehat{L}W]$, however, using the Jacobi identity we obtain the equality
\[ [W,\widehat{L}W] = [\widehat{L}W,\widehat{L}W]. \]
It follows that the map
\[ \zetiso[Lie]:\{ \omega \in \drtf[Lie]{\widehat{L}W} : d\omega = 0 \} \to [\widehat{L}W,\widehat{L}W] \]
defined by the formula $\zetiso[Lie](x):=glf_\mathrm{Lie}^{-1}(x)$ is an isomorphism.
\end{enumerate}
\end{proof}

\begin{rem}
It follows from the definitions that the maps \zetiso[Com], \zetiso[Ass] and \zetiso[Lie] satisfy the following relations:
\begin{equation} \label{eqn_zetrel}
\begin{split}
\zetiso[Ass]\circ d & = (1-z) \circ \thiso[Ass]^{-1}, \\
\zetiso[Com] \circ d & = \zetiso[Ass] \circ d \circ j, \\
\zetiso[Lie] & = \zetiso[Ass] \circ l. \\
\end{split}
\end{equation}
\end{rem}

If we set $W:=U^*$ in Theorem \ref{thm_tfmaps}, where $U$ is a free graded module, then we can interpret the closed 2-forms as linear maps $(TU)^* \to \gf$. We now want to concentrate on 2-forms of order zero and determine what linear maps in $(TU)^*$ they give rise to under the maps \zetiso[Com], \zetiso[Ass] and \zetiso[Lie]. 2-forms of order zero are naturally closed since it is easy to see that any 2-form $\omega$ \emph{of order zero} has the following form:
\begin{displaymath}
\begin{array}{lcll}
\omega \in \drtf[Com]{\widehat{S}U^*} & \Rightarrow & \omega = \sum_i dx_i\cdot dy_i; & x_i,y_i \in U^*; \\
\omega \in \drtf[Ass]{\widehat{T}U^*} & \Rightarrow & \omega = \sum_i dx_i\cdot dy_i; & x_i,y_i \in U^*; \\
\omega \in \drtf[Lie]{\widehat{L}U^*} & \Rightarrow & \omega = \sum_i dx_i\cotimes dy_i; & x_i,y_i \in U^*. \\
\end{array}
\end{displaymath}

We have the following lemma which describes the maps that two forms of order zero give rise to:

\begin{lemma} \label{lem_tfmbmp}
Let $U$ be a free graded module. Every map in the following commutative diagram is an isomorphism:
\begin{displaymath}
\xymatrix{ \{ \omega \in \drtf[Com]{\widehat{S}U^*} : \ord(\omega) = 0 \} \ar_{j}[d] \ar^{\zetiso[Com]}[rd] \\ \{ \omega \in \drtf[Ass]{\widehat{T}U^*} : \ord(\omega) = 0 \} \ar^-{\zetiso[Ass]}[r] & (\Lambda^2 U)^* \\ \{ \omega \in \drtf[Lie]{\widehat{L}U^*} : \ord(\omega) = 0 \} \ar^{l}[u] \ar^{\zetiso[Lie]}[ru] \\ }
\end{displaymath}
\end{lemma}

\begin{proof}
This follows from Theorem \ref{thm_tfmaps} and the identity
\[ [U^*,U^*] = (1-z)\cdot U^*\cotimes U^* = (\Lambda^2 U)^*. \]
\end{proof}

This means there is a one-to-one correspondence between the module of 2-forms (commutative, associative or Lie) of order zero and the module consisting of all skew-symmetric bilinear forms
\[ \innprod{U}. \]
We will now recall the definition of a nondegenerate 2-form (cf. \cite{ginzsg}) and a nondegenerate bilinear form:

\begin{defi} \label{def_nondeg}
\
\begin{enumerate}
\item[(i)]
Let $X$ be either a formal graded commutative, associative or Lie algebra and let $\omega \in \drtf{X}$ be a 2-form. We say $\omega$ is nondegenerate if the following map is a bijection;
\begin{displaymath}
\begin{array}{ccc}
\Der(X) & \to & \drof{X}, \\
\xi & \mapsto & i_\xi(\omega). \\
\end{array}
\end{displaymath}
\item[(ii)]
Let $U$ be a free graded module of finite rank and let \innprod{U} be a bilinear form. We say that $\langle -,- \rangle$ is nondegenerate if the following map is a bijection;
\begin{displaymath}
\begin{array}{ccc}
U & \to & U^*, \\
u & \mapsto & [ x \mapsto \langle u,x \rangle ]. \\
\end{array}
\end{displaymath}
If in addition $\langle -,- \rangle$ is symmetric then we will call it an \emph{inner product} on $U$.
\end{enumerate}
\end{defi}

We have the following proposition relating the two notions:

\begin{prop} \label{prop_nondeg}
Let $U$ be a free graded module of finite rank and let $X$ be one of the three algebras $\widehat{S}U^*$, $\widehat{T}U^*$ or $\widehat{L}U^*$. Let $\omega \in \drtf{X}$ be a 2-form of order zero. Let $\langle -,- \rangle:=\zeta(\omega)$ be the skew-symmetric bilinear form corresponding to the 2-form $\omega$, then $\omega$ is nondegenerate if and only if $\langle -,- \rangle$ is nondegenerate.
\end{prop}

\begin{proof}
We will treat the three cases $X=\widehat{S}U^*$, $X=\widehat{T}U^*$ and $X=\widehat{L}U^*$ simultaneously. Let $x_1,\ldots,x_n$ be a basis of the free module $U$. There are coefficients $a_{ij} \in \gf$ such that
\[ \omega = \sum_{1\leq i,j \leq n} a_{ij}dx_i^* dx_j^*. \]
It follows from the definition of $\zeta$ that
\[ \langle a,b \rangle = \sum_{1 \leq i,j \leq n} a_{ij}\left[(-1)^{|x_i|(|x_j|+1)}x_j^*(a)x_i^*(b) - (-1)^{|x_i|}x_i^*(a)x_j^*(b)\right]. \]

Let us define the map $\Phi: \Der(X) \to \drof{X}$ by the formula
\[ \Phi(\xi):= i_\xi(\omega). \]
Let us also define a map $D:U \to U^*$ by the formula
\[ D(u):= [x \mapsto \langle u,x \rangle]. \]

We calculate $\Phi$ as follows: Let $\xi \in \Der(X)$, then
\begin{displaymath}
\begin{split}
\Phi(\xi) = i_\xi(\omega) & = \sum_{1 \leq i,j \leq n} a_{ij}[\xi(x_i^*)\cdot dx_j^* + (-1)^{(|x_i|+1)(|\xi|+1)}dx_i^*\cdot \xi(x_j^*)], \\
& =  \sum_{1 \leq i,j \leq n} a_{ij}[(-1)^{(|x_j|+1)(|x_i|+|\xi|)}dx_j^*\cdot \xi(x_i^*) + (-1)^{(|x_i|+1)(|\xi|+1)}dx_i^*\cdot \xi(x_j^*)]. \\
\end{split}
\end{displaymath}
We calculate $D$ as follows: For all $a \in U$, and $1\leq k \leq n$;
\begin{displaymath}
\begin{split}
D(x_k)[a] = \langle x_k,a \rangle & = \sum_{1 \leq i,j \leq n} a_{ij}\left[(-1)^{|x_i|(|x_j|+1)}x_j^*(x_k)x_i^*(a) - (-1)^{|x_i|}x_i^*(x_k)x_j^*(a)\right], \\
& = \sum_{1\leq i \leq n}\left[(-1)^{|x_i|(|x_k|+1)}a_{ik} - (-1)^{|x_k|}a_{ki}\right]x_i^*(a). \\
\end{split}
\end{displaymath}

Let $h:U^* \to X$ be the canonical inclusion of the submodule $U^*$ into $X$. Since any vector field $\xi \in \Der(X)$ is completely determined by its restriction $\xi\circ h$ (cf. Proposition \ref{prop_dergen}) we have the following isomorphism of graded modules:
\begin{displaymath}
\begin{array}{ccc}
\Der(X) & \cong & \Hom_\gf(U^*,X), \\
\xi & \mapsto & \xi\circ h. \\
\end{array}
\end{displaymath}
It follows from the preceding calculations that the following diagram is commutative:
\[ \xymatrix{ \drof{X} \ar^{\Theta^{-1}}_{\cong}[r] & U^* \otimes X \ar^{u \otimes x \mapsto (-1)^{|\omega|+(|u|+1)(|x|+1)}u \otimes x}[rr] && U^* \otimes X \\ \Der(X) \ar^{\Phi}[u] \ar_-{\cong}^-{\xi \mapsto \xi \circ h}[r] & \Hom_\gf(U^*,X) \ar@{=}[r] & U^{**} \otimes X \ar@{=}[r] & U \otimes X \ar_{D \otimes 1}[u] \\ } \]
We conclude that $\Phi$ is a bijection if and only if $D$ is a bijection.
\end{proof}

\begin{rem} \label{rem_canfrm}
For simplicity's sake, let us assume we are working over a field of characteristic zero which is closed under taking square roots and that $U$ is finite dimensional; then every homogeneous 2-form $\omega \in \drtf{X}$ (where $X=\widehat{S}U^*$, $\widehat{T}U^*$ or $\widehat{L}U^*$) of order zero has a canonical form.

Suppose that $\omega$ has \emph{even} degree, then there exist linearly independent vectors;
\begin{equation} \label{eqn_canfrmdummy}
p_1,\ldots,p_n;q_1,\ldots,q_n;x_1,\ldots,x_m \in U
\end{equation}
(where the $p_i$'s and the $q_i$'s are \emph{even} and the $x_i$'s are \emph{odd}) such that $\omega$ has the following form:
\[ \omega = \sum_{i=1}^n d p_i^* d q_i^* + \sum_{i=1}^m d x_i^* d x_i^*. \]
If $\omega$ is nondegenerate then \eqref{eqn_canfrmdummy} is a basis for $U$. The bilinear form $\langle -,- \rangle:=\zeta(\omega)$ is given by the formula:
\begin{displaymath}
\begin{array}{rcr}
\langle q_i,p_j \rangle = \frac{1}{2}\langle x_i,x_j \rangle & = & \delta_{ij}, \\
\langle x_i,p_j \rangle = \langle x_i,q_j \rangle = \langle p_i,p_j \rangle = \langle q_i,q_j \rangle & = & 0. \\
\end{array}
\end{displaymath}

Now suppose that $\omega$ has \emph{odd} degree, then there exist linearly independent vectors;
\begin{equation} \label{eqn_canfrmdummya}
x_1,\ldots,x_n;y_1,\ldots,y_n \in U
\end{equation}
(where the $x_i$'s have \emph{odd} degree and the $y_i$'s have \emph{even} degree) such that $\omega$ has the following form:
\[ \omega = \sum_{i=1}^n d x_i^* d y_i^*. \]
Again, if $\omega$ is nondegenerate then \eqref{eqn_canfrmdummya} is a basis for $U$. The bilinear form $\langle -,- \rangle:=\zeta(\omega)$ is given by the formula:
\begin{displaymath}
\begin{array}{rcr}
\langle x_i,y_j \rangle & = & \delta_{ij}, \\
\langle x_i,x_j \rangle = \langle y_i,y_j \rangle & = & 0. \\
\end{array}
\end{displaymath}
\end{rem}

\section{Infinity algebra Prerequisites} \label{prerequisites}

In this section we will review the definitions of three types of infinity algebra, namely {\ci}, {\ai} and {\li}. These are the (strong) homotopy generalisations of commutative, associative and Lie algebras respectively. We will also define the appropriate notion of a unital \ai and \ci-algebra. We shall utilise the duality between derivations and coderivations to define an $\infty$-structure as a homological vector field on a certain formal supermanifold. This approach will provide us with the natural framework in which to apply constructions from noncommutative geometry. 

Let $V$ be a free graded module and choose a topological basis $\boldsymbol{t}:=\{t_i\}_{i \in I}$ of $\Sigma V^*$, then $\ctalg{V}=\gf\langle\langle\boldsymbol{t}\rangle\rangle$. Recall from section \ref{sec_notcon} that we say a vector field (continuous derivation) $\xi:\ctalg{V}\to\ctalg{V}$ vanishes at zero if it has the form
\[ \xi=\sum_{i \in I}A_i(\boldsymbol{t})\partial_{t_i}, \]
where the power series $A_i(\boldsymbol{t}), i \in I$ have vanishing constant terms.

We will now recall from \cite{azksch}, \cite{getjon} and \cite{lazmod} the definition of an \li, {\ai} and \ci-structure on a free graded module $V$. An \li, {\ai} or \ci-algebra is a free graded module together with an \li, {\ai} or \ci-structure. Given an element $x$ in a graded algebra $A$ we define the derivation $\ad x:A \to A$ by the formula $\ad x(y):=[x,y]$.

\begin{defi} \label{def_infstr}
Let $V$ be a free graded module:
\begin{enumerate}
\item[(a)]
An \li-structure on $V$ is a vector field
\[ m:\csalg{V} \to \csalg{V} \]
of degree one and vanishing at zero; such that $m^2=0$.
\item[(b)]
An \ai-structure on $V$ is a vector field
\[ m:\ctalg{V} \to \ctalg{V} \]
of degree one and vanishing at zero; such that $m^2=0$.
\item[(c)]
A \ci-structure on $V$ is a vector field
\[ m:\clalg{V} \to \clalg{V} \]
of degree one; such that $m^2=0$.
\end{enumerate}
\end{defi}

\begin{rem} \label{rem_infcel}
For most applications the definition just given will suffice. However for an arbitrary graded ring $\gf$ (as opposed to a field) it might lead to homotopy noninvariant constructions. E.g. the Hochschild cohomology of two weakly equivalent $A_\infty$-algebras may not be isomorphic. In order to avoid troubles like this the definition needs to be modified as follows. In all three cases the $\infty$-structure $m$ can be represented as
\[ m=m_1+m_2+\ldots+m_n+\ldots, \]
where $m_i$ is a vector field of order $i$. The condition 
$m^2=0$ implies that $m_1^2=0$. In other words $V$ together with the self-map $m_1$ forms a complex of $\gf$-modules. We then require that this complex be \emph{cellular}, in the sense of  \cite{krizmay}. This requirement is extraneous when for example, $V$ is finitely generated and free over $\gf$ or when $\gf$ is a field. 
\end{rem}

\begin{defi} \label{def_infmor}
Let $V$ and $U$ be free graded modules:
\begin{enumerate}
\item[(a)]
Let $m$ and $m'$ be \li-structures on $V$ and $U$ respectively: An \li-morphism from $V$ to $U$ is a continuous algebra homomorphism
\[ \phi:\csalg{U} \to \csalg{V} \]
of degree zero such that $\phi \circ m'=m \circ \phi$.
\item[(b)]
Let $m$ and $m'$ be \ai-structures on $V$ and $U$ respectively: An \ai-morphism from $V$ to $U$ is a continuous algebra homomorphism
\[ \phi:\ctalg{U} \to \ctalg{V} \]
of degree zero such that $\phi \circ m'=m \circ \phi$.
\item[(c)]
Let $m$ and $m'$ be \ci-structures on $V$ and $U$ respectively: A \ci-morphism from $V$ to $U$ is a continuous algebra homomorphism
\[ \phi:\clalg{U} \to \clalg{V} \]
of degree zero such that $\phi \circ m'=m \circ \phi$.
\end{enumerate}
\end{defi}

\begin{rem} \label{rem_calseq}
We have a diagram of functors
\begin{equation} \label{fig_calseq}
\ci\mathrm{-}algebras \longrightarrow \ai\mathrm{-}algebras \longrightarrow \li\mathrm{-}algebras
\end{equation}
which depends upon the sequence of Lie algebra homomorphisms defined in Remark \ref{rem_vecseq}.

Recall that any vector field $m:\clalg{V} \to \clalg{V}$ can be uniquely extended to a continuous Hopf algebra derivation (where \ctalg{V} is equipped with the cocommutative comultiplication (\emph{shuffle coproduct}) defined in Remark \ref{rem_ccmult})
\[ m : \ctalg{V} \to \ctalg{V} \]
and that all continuous Hopf algebra derivations on \ctalg{V} are obtained from $\Der(\clalg{V})$ in this manner. It follows that any \ci-structure gives rise to an \ai-structure in this way. Similarly, any continuous Lie algebra homomorphism $\phi:\clalg{U} \to \clalg{V}$ can be uniquely extended to a continuous Hopf algebra homomorphism
\[ \phi:\ctalg{U} \to \ctalg{V} \]
and all continuous Hopf algebra homomorphisms are obtained from continuous Lie algebra homomorphisms in this way. This is because the Lie subalgebra of primitive elements of our Hopf algebra \ctalg{V} coincides with the Lie subalgebra \clalg{V}.

The category $\ai\mathrm{-}algebras$ has a subcategory whose objects consist of the \ai-algebras whose \ai-structure is also a Hopf algebra derivation and whose morphisms are the \ai-morphisms which are also Hopf algebra homomorphisms. It follows from the discussion above that this category is isomorphic to $\ci\mathrm{-}algebras$.

Recall as well that any vector field $m:\ctalg{V} \to \ctalg{V}$ can be lifted to a unique vector field
\[ \tilde{m}: \csalg{V} \to \csalg{V}. \]
It follows that any \ai-structure gives rise to a \li-structure in this manner. Similarly any continuous algebra homomorphism $\phi:\ctalg{U} \to \ctalg{V}$ can be lifted to a unique continuous homomorphism
\[ \tilde{\phi}:\csalg{U} \to \csalg{V}. \]

Observe that given a minimal \ci-algebra, the corresponding \li-algebra under the maps of diagram \eqref{fig_calseq} is a trivial algebra (i.e. its multiplication sends everything to zero).
\end{rem}

We will now define an important type of {\ai} and \ci-algebra; unital {\ai} and \ci-algebras:

\begin{defi}
Let $V$ be a free graded module:
\begin{enumerate}
\item[(i)]
\begin{enumerate}
\item
We say that an \ai-structure $m:\ctalg{V} \to \ctalg{V}$ is unital if there is a distinguished element $1 \in V$ (the unit) of degree zero which can be extended to a basis $1,\{x_i\}_{i\in I}$ of $V$ such that $m$ has the form;
\[ m=A(\boldsymbol{t})\partial_\tau+\sum_{i\in I}B_i(\boldsymbol{t})\partial_{t_i} + \ad\tau - \tau^2\partial_\tau, \]
where $\tau,\boldsymbol{t}:=\tau,\{t_i\}_{i\in I}$ is the topological basis of $\Sigma V^*$ which is dual to the basis $\Sigma 1,\{\Sigma x_i\}_{i\in I}$ of $\Sigma V$. In this case we say that the \ai-algebra $V$ is unital or that it has a unit.
\item
Suppose that $V$ and $U$ are two unital \ai-algebras. We say that an \ai-morphism $\phi:\ctalg{U} \to \ctalg{V}$ is unital if $\phi$ has the form;
\begin{displaymath}
\begin{array}{ccc}
\phi(\tau') & = & \tau + A(\boldsymbol{t}), \\
\phi(t'_i) & = & B_i(\boldsymbol{t}); \\
\end{array}
\end{displaymath}
where $\tau,\boldsymbol{t}$ and $\tau',\boldsymbol{t}'$ are the topological bases of $\Sigma V^*$ and $\Sigma U^*$ which are dual to the bases $\Sigma 1_V,\{\Sigma x_i\}_{i\in I}$ and $\Sigma 1_U,\{\Sigma x'_j\}_{j\in J}$ of $\Sigma V$ and $\Sigma U$ respectively.
\end{enumerate}
\item[(ii)]
\begin{enumerate}
\item
We say that a \ci-structure $m:\clalg{V} \to \clalg{V}$ is unital if the corresponding \ai-structure (see Remark \ref{rem_calseq}) is unital.
\item
Suppose that $V$ and $U$ are two unital \ci-algebras. We say that a \ci-morphism $\phi:\clalg{U} \to \clalg{V}$ is unital if the corresponding \ai-morphism (see Remark \ref{rem_calseq}) is unital.
\end{enumerate}
\end{enumerate}
\end{defi}

\begin{rem} \label{rem_infdef}
There is an alternative definition of an $\infty$-structure on a free graded module $V$. This is the definition that was originally introduced by Stasheff in \cite{staha1} and \cite{staha2}. According to this definition an $\infty$-structure on $V$ is a system of maps
\[ \check{m}_i:V^{\otimes i} \to V, i \geq 1 \]
(where $|\check{m}_i|=2-i$) satisfying the higher homotopy axioms (and possibly some graded symmetry axioms as well). For \ai-algebras for instance these axioms imply that $\check{m}_2$ is associative up to homotopy (the homotopy being provided by $\check{m}_3$) whilst for \li-algebras they imply that $\check{m}_2$ satisfies the Jacobi identity up to homotopy. In particular, in all three cases they specify that the map $\check{m}_1$ is a differential and a graded derivation with respect to the multiplication $\check{m}_2$.

There is also an alternative definition of $\infty$-morphisms in this context. These are a system of maps
\[ \check{\phi}_i:V^{\otimes i} \to U, i \geq 1 \]
(where $|\check{\phi}_i|=1-i$) satisfying certain compatibility conditions with the $\check{m}_i$ and $\check{m}'_i$'s which specify the $\infty$-structures on $V$ and $U$ respectively. In particular the map $\check{\phi}_1$ must be a map of complexes; $\check{\phi}_1\circ\check{m}_1=\check{m}'_1\circ\check{\phi}_1$. We say that the $\infty$-morphism given by this system of maps is a weak equivalence if the map
\[ \check{\phi}_1: (V,\check{m}_1) \to (U,\check{m}'_1) \]
is a quasi-isomorphism. If $\check{m}_1=0$ then we say that the \ai-algebra is minimal. It is easy to see that a weak equivalence of two minimal $\infty$-algebras is in fact an isomorphism.
\end{rem}

Let us describe how this alternative style of definition is equivalent to that described in definitions \ref{def_infstr} and \ref{def_infmor} beginning with the {\ai} case. Firstly there is a one-to-one correspondence between systems of maps $\check{m}_i:V^{\otimes i} \to V, i\geq 1$ and systems of maps $m_i:\Sigma V^{\otimes i} \to \Sigma V, i\geq 1$ via the following commutative diagram:
\begin{equation} \label{fig_infdef}
\xymatrix{V^{\otimes i} \ar^{\check{m}_i}[r] & V \\ \Sigma V^{\otimes i} \ar^{(\Sigma^{-1})^{\otimes i}}[u] \ar_{m_i}[r] & \Sigma V \ar^{\Sigma^{-1}}[u]}
\end{equation}
Of course the $m_i$'s will inherit additional signs from the Koszul sign rule. It is well known that any system of maps $m_i:\Sigma V^{\otimes i} \to \Sigma V,i\geq 1$ can be uniquely extended to a coderivation $m$ on the tensor coalgebra $T\Sigma V$ which vanishes on $\gf\subset T\Sigma V$. Furthermore all coderivations vanishing on {\gf} are obtained in this way, hence there is a one-to-one correspondence
\[ \Hom_\gf(T\Sigma V/\gf,\Sigma V) \leftrightarrow \{ m \in \Coder(T\Sigma V) : m(\gf) = 0 \}. \]
The condition $m^2=0$ turns out to be equivalent to the higher homotopy associativity axioms for the $\check{m}_i$'s. Now simply observe that the dual of a coderivation on $T\Sigma V$ is a continuous derivation on $(T\Sigma V)^*=\ctalg{V}$. It follows from Proposition \ref{prop_antieq} that our two definitions of an \ai-structure are equivalent.

\ai-morphisms are dealt with in a similar manner. Again there is a one-to-one correspondence between systems of maps $\check{\phi}_i:V^{\otimes i} \to U, i\geq 1$ and systems of maps $\phi_i:\Sigma V^{\otimes i} \to \Sigma U,i\geq 1$. It is well known that systems of such maps are in one-to-one correspondence with coalgebra morphisms $\phi:T\Sigma V \to T\Sigma U$. The dual of a coalgebra morphism $\phi$ is a continuous algebra homomorphism $\phi^*:\ctalg{U} \to \ctalg{V}$. The condition in Definition \ref{def_infmor} which stipulates that $\phi^*$ commutes with the \ai-structures reflects the compatibility conditions placed on the $\check{\phi}_i$'s alluded to in Remark \ref{rem_infdef}. Further details on the formal passage to the dual framework can be found in Appendix \ref{app_todual}.

The {\ci} case is a restriction of the {\ai} case. Certain graded symmetry conditions are placed on the maps $\check{m}_i:V^{\otimes i} \to V$ in addition to the higher homotopy associativity conditions. These symmetry conditions are satisfied if and only if the corresponding \ai-structure $m:\ctalg{V} \to \ctalg{V}$ is a Hopf algebra derivation (i.e. a derivation and a coderivation), where \ctalg{V} is endowed with the cocommutative comultiplication (\emph{shuffle coproduct}) defined in Remark \ref{rem_ccmult}. Similarly certain graded symmetry conditions are placed on the maps $\check{\phi}_i:V^{\otimes i} \to U$ in addition to the conditions requiring that the $\check{\phi}_i$'s are compatible with the \ai-structures. These symmetry conditions are satisfied if and only if the corresponding \ai-morphism $\phi:\ctalg{U} \to \ctalg{V}$ is a Hopf algebra homomorphism. It now follows from Remark \ref{rem_calseq} that our two definitions of a \ci-structure are equivalent.

The {\li} case is slightly different. Systems of maps $m_i:S^i\Sigma V \to \Sigma V, i\geq 1$ are in one-to-one correspondence with coderivations $m:S\Sigma V \to S\Sigma V$ which vanish on $\gf\subset S\Sigma V$, where $S\Sigma V$ is equipped with the cocommutative comultiplication defined in Remark \ref{rem_ccmult}. Similarly systems of maps $\phi_i: S^i\Sigma V \to \Sigma U,i\geq 1$ are in one-to-one correspondence with coalgebra homomorphisms $\phi:S\Sigma V \to S\Sigma U$. The condition $m^2=0$ corresponds exactly to the higher homotopy Lie axioms. The dual of a coderivation on $S\Sigma V$ which vanishes on {\gf} is a continuous derivation $m^*:(S\Sigma V)^* \to (S\Sigma V)^*$ which vanishes at zero, where $(S\Sigma V)^*$ is endowed with the shuffle product. Similarly the dual of a coalgebra morphism $\phi:S\Sigma V \to S\Sigma U$ is a continuous algebra homomorphism $\phi^*:(S\Sigma U)^* \to (S\Sigma V)^*$.

Since coinvariants are dual to invariants we have the following identity:
\[ (S\Sigma V)^*=\gf \times (\Sigma V^*)^{S_1} \times ((\Sigma V^*)^{\cotimes 2})^{S_2} \times ((\Sigma V^*)^{\cotimes 3})^{S_3} \times \ldots. \]
Recall that in equation \eqref{eqn_symiso} we defined the map $i:\csalg{V} \to (S\Sigma V)^*$ which is the canonical map identifying coinvariants with invariants and is given by the formula
\begin{equation} \label{eqn_liduin}
i(x_1 \cotimes \ldots \cotimes x_n):= \sum_{\sigma \in S_n} \sigma \cdot x_1 \cotimes \ldots \cotimes x_n.
\end{equation}
This map has an inverse $\pi:(S\Sigma V)^* \to \csalg{V}$ given by the projection
\begin{equation} \label{eqn_lidupr}
\pi(x_1 \cotimes \ldots \cotimes x_n):= \frac{1}{n!} x_1 \cotimes \ldots \cotimes x_n.
\end{equation}
Using the maps $i$ and $\pi$ to identify the module $(S\Sigma V)^*$ with the module \csalg{V}, the shuffle product on $(S\Sigma V)^*$ is transformed into the canonical commutative multiplication on \csalg{V} which is inherited from the associative multiplication on \ctalg{V}. From this it follows that our two definitions of a \li-structure are equivalent.

Further details on $\infty$-algebras and their definitions can be found in \cite{keller}, \cite[\S 5]{getjon} and \cite{markl}.

\section{Minimal Infinity algebras}\label{minimal}

In this section we prove that any $C_\infty$-algebra is weakly equivalent to a minimal one. This theorem in the context of $A_\infty$-algebras was first proved by Kadeishvili in \cite{kadvil} but it had its precursors in the theory of minimal models in rational homotopy theory, cf. \cite{sulvan}. It has since been reproved by many authors, cf. \cite{keller} and references therein. The proof in the $L_\infty$ case was outlined in \cite{kontpm}. It seems that there is no published proof in the $C_\infty$ case in the literature. We will give here a short proof based on the notion of the Maurer-Cartan moduli space associated to a differential graded Lie algebra, cf. \cite{golmil}. With obvious modifications our proof works for $A_\infty$ and $L_\infty$ cases as well.

In this section we make an assumption that $\gf$ is a field, or, more generally, a graded field (i.e. a graded commutative ring whose homogeneous nonzero elements are invertible).

We first recall some standard facts from the Maurer-Cartan theory following \cite{golmil}. Let $\mathcal{G}$ be a differential graded Lie algebra which we assume to be nilpotent, or, more generally, an inverse limit of nilpotent differential graded Lie algebras. The set $\mathcal{MC}(\mathcal{G})\subset \mathcal{G}^1$ is by definition the set of solutions of the Maurer-Cartan equation
\begin{equation} \label{MC}
d\gamma+\frac{1}{2}[\gamma,\gamma]=0.
\end{equation}
Furthermore, the Lie algebra $\mathcal{G}^0$ acts on $\mathcal{MC}(\mathcal{G})$ by infinitesimal affine transformations: for $\alpha\in \mathcal{G}^0$ we have a vector field $\gamma \mapsto d\alpha+[\gamma,\alpha]$ on  $\mathcal{MC}(\mathcal{G})$. This action exponentiates to an action of the Lie group $\exp(\mathcal{G}^0)$ on $\mathcal{MC}(\mathcal{G})$. We will call the set of orbits with respect to this action \emph{the Maurer-Cartan moduli space} associated to the differential graded Lie algebra $\mathcal{G}$.

Now let $\mathcal{G}_1,\mathcal{G}_2$ be two nilpotent differential graded Lie algebras. We assume that they are endowed with finite filtrations $\{F_p(\mathcal{G}_1)\}$ and $\{F_p(\mathcal{G}_2)\}$, $p=1,2,\ldots, n$ and there is a map of filtered differential graded Lie algebras $f:\mathcal{G}_1\rightarrow \mathcal{G}_2$ which induces a quasi-isomorphism on the associated graded of $\mathcal{G}_1$ and $\mathcal{G}_2$. Under these assumptions we have the following result.

\begin{theorem} \label{getgm}
The map $f$ induces a bijection $\mathcal{MC}(\mathcal{G}_1)\rightarrow \mathcal{MC}(\mathcal{G}_2)$.
\end{theorem}

 This seems to be a well known result and is formulated in this form in \cite{getzlr}, Theorem 2.1. However Goldman-Millson's version \cite{golmil}, Theorem 2.4 does not readily carry over since we are dealing with differential graded Lie algebras not necessarily concentrated in nonnegative degrees. We therefore sketch a proof suitable for this more general situation. This proof is modelled on \cite{getzs}, Proposition 4.6. We start with the definition of the simplicial Maurer-Cartan set (also called the simplicial Deligne groupoid) associated to a differential graded Lie algebra $\mathcal{G}$. A more detailed discussion could be found in e.g. \cite{hinich}.  Let $\Omega_n$ be the algebra of polynomial differential forms on the standard $n$-simplex; the collection $\Omega_\bullet=\{\Omega_n\}_{n=0}^\infty$ forms a commutative simplicial differential graded algebra.  Note that $\mathcal{G}\otimes \Omega_n$ has the structure of a differential graded Lie algebra. 

\begin{defi}
For $n\geq 0$ set $\mathcal{MC}_n(\mathcal{G}):= \mathcal{MC}(\mathcal{G}\otimes \Omega_n)$. The simplicial structure on $\Omega_\bullet$ determines the structure of a simplicial set on  $\mathcal{MC}_\bullet(\mathcal{G}):=\{\mathcal{MC}_n(\mathcal{G})\}_{n=0}^\infty$ which will be called \emph{the Maurer-Cartan simplicial set} associated with $\mathcal{G}$.
\end{defi}

Then the main result of \cite{sstas} implies that there is a one-to-one correspondence between $\mathcal{MC}(\mathcal{G})$ and $\pi_0(\mathcal{MC}_\bullet(\mathcal{G}))$, the set of connected components of the simplicial set $\mathcal{MC}_\bullet(\mathcal{G})$. Furthermore, for a surjective map of differential graded Lie algebras $\mathcal{G}\rightarrow \mathcal{G}^\prime$ the induced map $\mathcal{MC}_\bullet(\mathcal{G})\rightarrow \mathcal{MC}_\bullet(\mathcal{G}^\prime)$ is a fibration of simplicial sets. This is proved (in a more general context for $L_\infty$-algebras) in \cite{getzs}.

Finally, using induction up the filtrations of $\mathcal{G}_1$ and $\mathcal{G}_2$ and comparing the associated towers of fibrations of simplicial sets we show that the simplicial sets $\mathcal{MC}(\mathcal{G}_1)$ and $\mathcal{MC}(\mathcal{G}_2)$ are weakly equivalent. In particular, their sets of connected components are in one-to-one correspondence. This finishes our sketch proof of Theorem \ref{getgm}.

Now let $V$ be a graded vector space and $m$ be a \ci-structure on $V$.  Then $m$ is a vector field on the Lie algebra $\widehat{L}(\Sigma V^*)$. We will denote $\Sigma V^*$ by $W$. Choose a topological basis $\{x_i\}_{i \in I}$ in $W$ and denote by $\mathcal{G}$ the Lie algebra of vector fields having the form $\sum_{i\in I} f_i\partial_{x_i}$, where the $f_i$'s are (possibly uncountably infinite) sums of Lie monomials in $x_i$ of order 2 or higher. Clearly, this definition of $\mathcal{G}$ does not depend on the choice of a basis and $\mathcal{G}$ is a formal Lie algebra. The Lie group $\exp(\mathcal{G}^0)$ is the group of formal Lie series whose linear term is the identity linear transformation.

We will view $\mathcal{G}$ as a differential graded algebra with respect to the operator $m_1$, the linear part of the vector field $m$. Denote by $Z^\bullet(W)$ and $H^\bullet(W)$ respectively, the cocycles and cohomology of $W$ with respect to $m_1$. 

 Note that $\mathcal{G}$ has a filtration given by the order of a vector field. This filtration is in fact a grading and $m_1$ preserves this grading since a commutator of a linear vector filed and a vector field of order $n$ is again a vector field of order $n$. Therefore the Lie algebra $H^\bullet(\mathcal{G})$ is bigraded.

Furthermore, direct inspection shows that $m-m_1$ is a Maurer-Cartan element in $\mathcal{G}$, i.e. it satisfies  equation \eqref{MC}. Moreover, two Maurer-Cartan elements are  equivalent with respect to the action of $\exp(\mathcal{G}^0)$ if and only if the corresponding \ci-structures are weakly equivalent through a formal diffeomorphism $f=(f_1,f_2,\ldots)$ of $\widehat{L}(W)$ whose linear part $f_1$ is the identity map.   Following \cite{polchu} we will call such equivalences \emph{strict $C_\infty$-isomorphisms}.

Let us now fix an integer $n>1$ and denote by $\mathcal{G}_n$ the quotient of $\mathcal{G}$ by the ideal of vector fields having order $>n$. Using the geometric langauge we say that $\mathcal{G}_n$ is the group of germs of formal diffeomorphisms of order $n$. Clearly $\mathcal{G}_n$ is a nilpotent differential graded algebra and $\mathcal{G}=\inlim{n} {\mathcal{G}_n}$.

\begin{prop} \label{form}
The differential graded Lie algebra $\mathcal{G}$  is quasi-isomorphic to its cohomology $H^\bullet(\mathcal{G})$ considered as a differential graded Lie algebra with trivial differential.
\end{prop}

\begin{rem}
The above proposition says that $\mathcal{G}$ is \emph{formal} in the sense of rational homotopy theory. We refrain from such a formulation, however, since the term `formal' has a different meaning for us.
\end{rem}

\begin{proof}
Let $Z^\bullet(W)$ be the space of $m_1$-cocycles in $W$ and choose a section of the projection map $Z^\bullet(W)\rightarrow H^\bullet(W)$. Then $H^\bullet(W)$ could be considered as a subspace in $Z^\bullet(W)$ and hence - also in $W$ itself.  Next, choose a section of the embedding  $H^\bullet(\mathcal{G})\hookrightarrow \mathcal{G}$. We thus have two maps $i:H^\bullet(\mathcal{G})\hookrightarrow \mathcal{G}$ and $j:\mathcal{G}\rightarrow H^\bullet(\mathcal{G})$ such that $j\circ i$ is the identity on $H^\bullet(\mathcal{G})$.

Denote by $\mathcal{H}$ the Lie algebra of vector fields on the Lie algebra $\widehat{L}(H^\bullet(W))$ which have order $\geq 2$. Define the maps $\tilde{i}:\mathcal{H} \rightarrow \mathcal{G}$ and $\tilde{j}:\mathcal{G}\rightarrow \mathcal{H}$ as follows: For $f\in \mathcal{H}, g\in \mathcal{G}, a\in \widehat{L}(W)$ and $b\in\widehat{L}(H^\bullet(W))$ set 
\[[\tilde{i}(f)](a):=i\circ g\circ j(a),\] \[[\tilde{j}(g)](b):=j\circ f\circ i(b).\]
It is immediate to check that $\tilde{j}\circ\tilde{i}=\id_{\mathcal{H}}$, i.e. that the Lie algebra $\mathcal{H}$ is a retract of $\mathcal{G}$. Note that $\mathcal{H}$ could be regarded as a differential graded Lie algebra with zero differential. We claim that the maps $\tilde{i}$ and $\tilde{j}$ induce mutually inverse homomorphisms between $H^\bullet(\mathcal{H})=\mathcal{H}$ and $H^\bullet(\mathcal{G})$. This will clearly imply the statement of the proposition.

Denote by $\widehat{L}^n(?), n=2,3,\ldots$ the space generated by the Lie monomials of length $n$ inside $\widehat{L}(?)$ and observe that $\tilde{i}$ and $\tilde{j}$ determine and are determined by the collection of maps
\[\tilde{i}_n:\Hom(H^\bullet(W),\hat{L}^n(H^\bullet(W))\leftrightarrows \Hom(V, \hat{L}^n(W)):\tilde{j}_n.\]

The functor $L$ associating to a graded vector space the free Lie algebra on it commutes with cohomology by \cite{quillen}, Appendix B, Proposition 2.1 (this is a consequence of the Poincar\'e-Birkhoff-Witt theorem). Since the inverse limit functor is exact on the category of finite dimensional vector spaces we conclude that $H^\bullet(\widehat{L}(W))\cong \widehat{L}(H^\bullet(W))$ and it follows that the maps $\tilde{i}_n,\tilde{j}_n$ are quasi-isomorphisms for all $n=2,3,\ldots$ so the claim follows.
\end{proof}

\begin{rem}
Note that our proof of Proposition \ref{form} actually gives slightly more than claimed. First, our proof shows the associated graded to $\mathcal{G}$ is quasi-isomorphic to its own cohomology. 
(This holds simply because the filtration on $\mathcal{G}$ is in fact a second grading.)

Second, denote by $\mathcal{H}_n$ the quotient of $\mathcal{H}$ by the ideal of vector fields of order $>n$. Then the restrictions of `formality maps' $\tilde{i}$ and $\tilde{j}$ determine  mutually (quasi-)inverse quasi-isomorphisms between $\mathcal{G}_n$ and $\mathcal{H}_n$. Invoking
again the analogy with rational homotopy theory one can express this by saying that the pronilpotent differential graded algebra $\mathcal{G}$ is \emph{continuously formal}. This is an important technical point which is necessary for the proof of the minimality theorem.
\end{rem}

\begin{cor} \label{kadeish}
(Minimality theorem) Let $(V,m_V)$ be a \ci-algebra. Then there exists a minimal \ci-algebra $(U, m_U)$ and a weak equivalence of \ci-algebras $(U, m_U)\rightarrow (V,m_V)$.
\end{cor}

\begin{proof}
Set $U:=H^\bullet(V)$ and $W:=\Sigma V^*$. Denote, as before, by $\mathcal{G}$ the Lie algebra of vector fields on $\widehat{L}(W)$ of order $\geq 2$ and similarly denote by $\mathcal{H}$ the Lie algebra of vector fields on $\widehat{L}(H^\bullet(W))$ whose order is $\geq 2$. Choosing a basis for representatives of cohomology classes, we will regard $H^\bullet(W)$ as a subspace in $W$. A choice of a complement will determine a map of differential graded Lie algebras $\tilde{i}: \mathcal{H}\rightarrow \mathcal{G}$ which is a quasi-isomorphism by Theorem \ref{form}. It suffices to show that $\tilde{i}$ induces a bijection on the Maurer-Cartan moduli spaces 
\[\mathcal{MC}(\mathcal{H})/\exp(\mathcal{H}^0) \overset{\cong}{\to} \mathcal{MC}(\mathcal{G})/\exp(\mathcal{G}^0).\]
Indeed, that would mean that any \ci-structure on $V$ could be reduced to a minimal one using a composition of a strict \ci-isomorphism and a linear projection $V\rightarrow U=H^\bullet(V)$.

To get the desired isomorphism note that the tower of differential graded Lie algebras $\mathcal{H}_2\leftarrow \mathcal{H}_3\leftarrow\ldots$ determines a tower of fibrations of simplicial sets $\mathcal{MC}_\bullet(\mathcal{H}_2)\leftarrow \mathcal{MC}_\bullet(\mathcal{H}_3)\ldots$. Similarly we have the tower of simplicial sets $\mathcal{MC}_\bullet(\mathcal{G}_2)\leftarrow \mathcal{MC}_\bullet(\mathcal{G}_3)\ldots$ associated with the tower of Lie algebras $\mathcal{G}_2\leftarrow \mathcal{G}_3\leftarrow\ldots$. Since $\mathcal{G}_n$ is quasi-isomorphic to $\mathcal{H}_n$ for all $n$ we conclude that these towers of fibrations are level-wise weakly equivalent. There are isomorphisms
\[ \inlim{n}{\mathcal{MC}_\bullet(\mathcal{G}_n)} \cong \mathcal{MC}_\bullet(\mathcal{G}), \]
\[ \inlim{n}{\mathcal{MC}_\bullet(\mathcal{H}_n)} \cong \mathcal{MC}_\bullet(\mathcal{H}) \]
and it follows that $\mathcal{MC}_\bullet(\mathcal{G})$ and $\mathcal{MC}_\bullet(\mathcal{H})$ are weakly equivalent simplicial sets. In particular, their sets of connected components are in one-to-one correspondence.
\end{proof}

\begin{rem}
Theorem \ref{form} and Corollary \ref{kadeish} extend in the context of $A_\infty$ and  \li-algebras. For \li-algebras the reference to the Poincar\'e-Birkhoff-Witt theorem is replaced by the fact that the functor of $S_n$-coinvariants is exact and therefore commutes with cohomology. For $A_\infty$-algebras the corresponding issue never arises and in fact the minimality theorem in the $A_\infty$ context is true in arbitrary characteristic.
\end{rem}

\section{The Cohomology of Infinity Algebras} \label{sec_iachom}

In this section we will define the various cohomology theories for \li, {\ai} and \ci-algebras that we will use throughout the rest of our paper. There will be quite a number of different cohomology theories to be defined. For \ai-algebras it will be useful to define additional quasi-isomorphic complexes computing the cohomology of the \ai-algebras. For example, this will allow us to describe a periodicity exact sequence for \ai-algebra cohomology. We will also prove some other basic facts about \ai-algebra cohomology which are the infinity-analogues of familiar results for strictly associative graded algebras.

\subsection{Hochschild, Harrison and Chevalley-Eilenberg theories for $\infty$-algebras}

We will begin by defining the cohomology theories which do not involve the action of the cyclic groups. These are the theories that control the deformations of $\infty$-algebras, cf. \cite{hamilt}, \cite{pensch} and \cite{penfia}.

\begin{defi}
Let $V$ be a free graded module:
\begin{enumerate}
\item[(a)]
Let $m:\csalg{V} \to \csalg{V}$ be an \li-structure. The Chevalley-Eilenberg complex of the \li-algebra $V$ with coefficients in $V$ is defined on the module consisting of all vector fields on \csalg{V}:
\[ \clac{V}{V}:=\Sigma^{-1}\Der(\csalg{V}). \]
The differential $d:\clac{V}{V} \to \clac{V}{V}$ is given by
\[ d(\xi):=[m,\xi], \quad \xi \in \Der(\csalg{V}). \]
The Chevalley-Eilenberg cohomology of $V$ with coefficients in $V$ is defined as the cohomology of the complex \clac{V}{V} and denoted by \hlac{V}{V}.
\item[(b)]
Let $m:\ctalg{V} \to \ctalg{V}$ be an \ai-structure. The Hochschild complex of the \ai-algebra $V$ with coefficients in $V$ is defined on the module consisting of all vector fields on \ctalg{V}:
\[ \choch{V}{V}:=\Sigma^{-1}\Der(\ctalg{V}). \]
The differential $d:\choch{V}{V} \to \choch{V}{V}$ is given by
\[ d(\xi):=[m,\xi], \quad \xi \in \Der(\ctalg{V}). \]
The Hochschild cohomology of $V$ with coefficients in $V$ is defined as the cohomology of the complex \choch{V}{V} and denoted by \hhoch{V}{V}.

\item[(c)]
Let $m:\clalg{V} \to \clalg{V}$ be a \ci-structure. The Harrison complex of the \ci-algebra $V$ with coefficients in $V$ is defined on the module consisting of all vector fields on \clalg{V}:
\[ \caq{V}{V}:=\Sigma^{-1}\Der(\clalg{V}). \]
The differential $d:\caq{V}{V} \to \caq{V}{V}$ is given by
\[ d(\xi):=[m,\xi], \quad \xi \in \Der(\clalg{V}). \]
The Harrison cohomology of $V$ with coefficients in $V$ is defined as the cohomology of the complex \caq{V}{V} and denoted by \haq{V}{V}.
\end{enumerate}
\end{defi}

\begin{rem}
The purpose of the desuspensions in the above definition is to make the grading consistent with the classical grading on these cohomology theories.
\end{rem}

\begin{rem}
The Chevalley-Eilenberg, Hochschild and Harrison complexes all have the structure of a differential graded Lie algebra under the commutator bracket $[-,-]$. It was Stasheff who realised in \cite{stasgb} that on the Hochschild cohomology of a strictly associative algebra this was the Gerstenhaber bracket introduced by Gerstenhaber in \cite{gerbra}. Stasheff considered the \ai-structure as a coderivation $m:T\Sigma V \to T\Sigma V$ and defined his bracket on $\Coder(T\Sigma V)$ as the commutator of coderivations. If we were to translate this structure directly to $\Der(\ctalg{V})$ via the contravariant functor of Proposition \ref{prop_antieq} then the induced Lie bracket $\{-,-\}$ would be the reverse of the commutator bracket
\[ \{\xi,\gamma\} = [\gamma,\xi]; \quad \xi,\gamma \in \Der(\ctalg{V}) \]
and the differential would be the map
\[ \xi \mapsto [\xi,m]; \quad \xi \in \Der(\ctalg{V}), \]
however, there is an isomorphism between this differential graded Lie structure and the differential graded Lie structure we defined on \choch{V}{V} above given by the map
\[ \xi \mapsto (-1)^{|\xi|(|\xi|+1)(2|\xi|+1)/6 + 1}\xi, \quad \xi \in \Der(\ctalg{V}). \]
\end{rem}

Now let us describe the cohomology theories which are dual to the corresponding cohomology theories defined above.

\begin{defi}
Let $V$ be a free graded module:
\begin{enumerate}
\item[(a)]
Let $m:\csalg{V} \to \csalg{V}$ be an \li-structure. The Chevalley-Eilenberg complex of the \li-algebra $V$ with coefficients in $V^*$ is defined on the module consisting of all 1-forms:
\[ \clac{V}{V^*}:=\Sigma\drof[Com]{\csalg{V}}. \]
The differential on this complex is the (suspension of the) Lie operator of the vector field $m$;
\[ L_m:\drof[Com]{\csalg{V}} \to \drof[Com]{\csalg{V}}. \]
The Chevalley-Eilenberg cohomology of $V$ with coefficients in $V^*$ is defined as the cohomology of the complex \clac{V}{V^*} and denoted by \hlac{V}{V^*}.
\item[(b)]
Let $m:\ctalg{V} \to \ctalg{V}$ be an \ai-structure. The Hochschild complex of the \ai-algebra $V$ with coefficients in $V^*$ is defined on the module consisting of all 1-forms:
\[ \choch{V}{V^*}:=\Sigma\drof[Ass]{\ctalg{V}}. \]
The differential on this complex is the (suspension of the) Lie operator of the vector field $m$;
\[ L_m:\drof[Ass]{\ctalg{V}} \to \drof[Ass]{\ctalg{V}}. \]
The Hochschild cohomology of $V$ with coefficients in $V^*$ is defined as the cohomology of the complex \choch{V}{V^*} and denoted by \hhoch{V}{V^*}.
\item[(c)]
Let $m:\clalg{V} \to \clalg{V}$ be a \ci-structure. The Harrison complex of the \ci-algebra $V$ with coefficients in $V^*$ is defined on the module consisting of all 1-forms:
\[ \caq{V}{V^*}:=\Sigma\drof[Lie]{\clalg{V}}. \]
The differential on this complex is the (suspension of the) Lie operator of the vector field $m$;
\[ L_m:\drof[Lie]{\clalg{V}} \to \drof[Lie]{\clalg{V}}. \]
The Harrison cohomology of $V$ with coefficients in $V^*$ is defined as the cohomology of the complex \caq{V}{V^*} and denoted by \haq{V}{V^*}.
\end{enumerate}
\end{defi}

\begin{rem}
Again the suspensions appear in order to keep the grading consistent with the classical grading on these cohomology theories.
\end{rem}

\begin{rem}
We did not prove that $L_m$ is indeed a differential, however this is obvious from Lemma \ref{lem_schf}:
\[ L_m^2=\frac{1}{2}[L_m,L_m]=L_{\frac{1}{2}[m,m]}=L_{m^2}=0. \]
\end{rem}

\begin{rem}
Our definition of the Harrison complex implies that in the case of an ungraded commutative algebra the
zeroth term of the complex is not present. This terminology differs slightly from what seems to be adopted in the modern literature (e.g. \cite{loday}) in that in the latter the algebra itself is considered as the zeroth term (note, however, that in his original paper \cite{harrison} Harrison only defined the 1st, 2nd and 3rd cohomology  groups). This distinction is unimportant since the zeroth term always splits off as a direct summand.
\end{rem}

Next we will define the bar cohomology of an \ai-algebra:

\begin{defi} \label{def_barhom}
Let $V$ be an \ai-algebra with \ai-structure $m:\ctalg{V} \to \ctalg{V}$. We define the map $b'$ as the restriction of (the suspension of) $m$ to $\cbr{V}:=\Sigma\left[\col{V}\right]$:
\begin{equation} \label{fig_barhom}
\xymatrix{\col{V} \ar^{b'}[d] \ar@<-0.5ex>@{^(->}[r] & \ctalg{V} \ar^{m}[d] \\ \col{V} \ar@<-0.5ex>@{^(->}[r] & \ctalg{V}}
\end{equation}
The map $b'$ is a differential on \cbr{V} and we define the bar cohomology of the \ai-algebra $V$ as the cohomology of the bar complex \cbr{V} and denote it by \hbr{V}.
\end{defi}

\begin{rem} \label{rem_chomid}
We can also define a map $b:\col{V} \to \col{V}$ by identifying \col{V} with the underlying module of the complex \choch{V}{V^*}. We do this via the module isomorphism \thiso[Ass] defined by Lemma \ref{lem_ofisom}:
\begin{equation} \label{fig_chomid}
\xymatrix{\col{V} \ar^{b}[d] \ar^{\thiso[Ass]}[r] & \drof[Ass]{\ctalg{V}} \ar^{L_m}[d] \\ \col{V} \ar^{\thiso[Ass]}[r] & \drof[Ass]{\ctalg{V}}}
\end{equation}
The complex whose underlying module is $\Sigma\left[\col{V}\right]$ and whose differential is (the suspension of) $b$ is by definition isomorphic to the Hochschild complex of $V$ with coefficients in $V^*$.
\end{rem}

Now we shall show that the bar cohomology of a unital \ai-algebra is trivial.

\begin{lemma} \label{lem_unitch}
Let $V$ be a unital \ai-algebra with unit $1 \in V$. There is a contracting homotopy $h:\cbr{V} \to \cbr{V}$ which is the dual of the map
\begin{displaymath}
\begin{array}{ccc}
\bigoplus_{i=1}^\infty \Sigma V^{\otimes i} & \to & \bigoplus_{i=1}^\infty \Sigma V^{\otimes i}, \\
x & \mapsto & 1 \otimes x. \\
\end{array}
\end{displaymath}
\end{lemma}

\begin{proof}
By Definition \ref{def_infstr} part (b) the \ai-structure $m:\fpsa \to \fpsa$ has the form
\[ m=A(\boldsymbol{t})\partial_\tau+\sum_{i\in I}B_i(\boldsymbol{t})\partial_{t_i} + \ad\tau - \tau^2\partial_\tau, \]
where $\tau$ is dual to the unit $1 \in V$. We have the following formula for the map $h$;
\begin{displaymath}
\begin{array}{cc}
h(\tau x) = x, & x \in \fpsa; \\
h(t_i x) = 0, & x \in \fpsa. \\
\end{array}
\end{displaymath}
We calculate
\begin{displaymath}
\begin{split}
b'h(t_i x)+hb'(t_i x) & = 0 + h[([\tau,t_i]+B_i(\boldsymbol{t}))x] \pm h(t_ib'(x)), \\
& = h(\tau t_i x)=t_i x. \\
b'h(\tau x)+hb'(\tau x) & = b'(x)+h(\tau^2 x) - h(\tau b'(x)), \\
& =b'(x)+ \tau x -b'(x) = \tau x. \\
\end{split}
\end{displaymath}
Similar calculations show that
\[ b'h(t_i)+hb'(t_i)= t_i \quad \text{and} \quad b'h(\tau)+hb'(\tau)= \tau. \]
therefore
\[ b'h+hb'=\id. \]
\end{proof}

Let $V$ be a unital \ai-algebra and let $\tau,\boldsymbol{t}$ be a topological basis of $\Sigma V^*$ where $\tau$ is dual to the unit $1 \in V$, then $\ctalg{V} = \fpsa$. We say a 1-form $\alpha \in \drof[Ass]{\ctalg{V}}$ is normalised if it is a linear combination of elements of the form $q\cdot dv$ where $v \in \Sigma V^*$ and $q \in \gf\langle\langle \boldsymbol{t} \rangle \rangle$;
\[ \alpha = A(\boldsymbol{t})d\tau + \sum_{i \in I} B_i(\boldsymbol{t})dt_i; \quad A(\boldsymbol{t}), B_i(\boldsymbol{t}) \in \fpsa. \]
We will denote the module of normalised 1-forms by \drnof[Ass]{\ctalg{V}}. It is clear to see that the map \thiso[Ass] defined by Lemma \ref{lem_ofisom} identifies the module of normalised 1-forms with the module $\Sigma V^* \cotimes \widehat{T}(\Sigma V/\gf)^*$.

\begin{prop} \label{prop_htpret}
Let $V$ be a unital \ai-algebra with \ai-structure $m:\ctalg{V} \to \ctalg{V}$. The normalised 1-forms $\Sigma\drnof[Ass]{\ctalg{V}}$ form a subcomplex of $\choch{V}{V^*}:= \left(\Sigma\drof[Ass]{\ctalg{V}}, L_m\right)$. Furthermore the subcomplex of normalised 1-forms is a chain deformation retract of \choch{V}{V^*}.
\end{prop}

\begin{proof}
By Definition \ref{def_infstr} the \ai-structure $m$ has the form
\[ m=A(\boldsymbol{t})\partial_\tau+\sum_{i\in I}B_i(\boldsymbol{t})\partial_{t_i} + \ad\tau - \tau^2\partial_\tau. \]
We calculate that for all $q \in \gf\langle\langle \boldsymbol{t} \rangle \rangle$ and $i \in I$;
\begin{displaymath}
\begin{array}{rcll}
L_m(q\cdot dt_i) & = & [\tau,q]\cdot dt_i +(-1)^{|q|+1}q\cdot d[\tau,t_i] & \mod \drnof[Ass]{\ctalg{V}}, \\
& = & [\tau,q]\cdot dt_i +(-1)^{|q|+1}q\cdot([d\tau,t_i]-[\tau,dt_i]) & \mod \drnof[Ass]{\ctalg{V}}, \\
& = & [\tau,q]\cdot dt_i +(-1)^{|q|}[q,t_i]\cdot d\tau +(-1)^{|q|}[q,\tau]\cdot dt_i & \mod \drnof[Ass]{\ctalg{V}}, \\
& = & (-1)^{|q|}[q,t_i]\cdot d\tau = 0 & \mod \drnof[Ass]{\ctalg{V}}. \\
\end{array}
\end{displaymath}
Similar calculations show that $L_m(q\cdot d\tau) \in \drnof[Ass]{\ctalg{V}}$ for all $q \in \gf\langle\langle \boldsymbol{t} \rangle \rangle$, hence \drnof[Ass]{\ctalg{V}} is a subcomplex of \choch{V}{V^*}.

Let $x_0,\ldots,x_n \in \Sigma V^*$, we say that the 1-form
\begin{equation} \label{eqn_htpretdummy}
x_1\ldots x_n\cdot dx_0
\end{equation}
is $i$-normalised if $x_1\ldots x_i \in \gf\langle\langle \boldsymbol{t} \rangle\rangle$. A generic $i$-normalised 1-form is a linear combination of $i$-normalised elements of the form \eqref{eqn_htpretdummy}. Obviously \eqref{eqn_htpretdummy} is $n$-normalised if and only if it is normalised.

Let $\gamma$ be the constant vector field $\gamma := 1\partial_\tau$. Define the map
\[ s_i:\drof[Ass]{\ctalg{V}} \to \drof[Ass]{\ctalg{V}} \]
by the formula
\[ s_i(x_1\ldots x_i\cdot x_{i+1} \cdot q\cdot dx_0) := (-1)^{|x_1|+\ldots+|x_i|}\gamma(x_{i+1})x_1\ldots x_i \cdot q\cdot dx_0, \]
where $x_0,\ldots,x_{i+1} \in \Sigma V^*$ and $q \in \fpsa$. We define the map
\[ h_i:\drof[Ass]{\ctalg{V}} \to \drof[Ass]{\ctalg{V}} \]
as $h_i:= \id + L_m s_i + s_i L_m$ which is of course chain homotopic to the identity.

We claim that $h_i$ takes $i$-normalised 1-forms to $i+1$-normalised 1-forms. Let $x,v \in \Sigma V^*,q \in \fpsa$ and let $p \in \gf\langle\langle \boldsymbol{t} \rangle\rangle$ be a power series of order $i$ so that $\alpha:=pxq\cdot dv$ is an $i$-normalised 1-form. Let $m_1$ be the linear part of the \ai-structure $m$. We will now calculate $h_i(\alpha)$:
\begin{equation} \label{eqn_htpretdummya}
\begin{split}
h_i(\alpha) = & \alpha + L_m((-1)^{|p|}\gamma(x)pq\cdot dv) + s_i(m(p)xq\cdot dv + (-1)^{|p|}pm(x)q\cdot dv \\
& + (-1)^{|p|+|x|}pxm(q)\cdot dv +(-1)^{|p|+|x|+|q|+1}pxq\cdot dm(v)), \\
= & \alpha + (-1)^{|p|}\gamma(x)[m(p)q\cdot dv + (-1)^{|p|}pm(q)\cdot dv + (-1)^{|p|+|q|+1}pq\cdot dm(v)] \\
& + s_i([\tau,p]q\cdot dv) + (-1)^{|p|+1}\gamma(x)m_1(p)q\cdot dv + (-1)^{|p|}s_i(pm(x)q\cdot dv) \\
& + (-1)^{|x|}\gamma(x)pm(q)\cdot dv + (-1)^{|p|+|x|+|q|+1}s_i(pxq\cdot dm(v)), \\
= & \alpha + (-1)^{|p|}\gamma(x)m(p)q\cdot dv + \gamma(x)pm(q)\cdot dv +(-1)^{|q|+1}\gamma(x)pq\cdot dm(v) - pxq\cdot dv \\
& + (-1)^{|p|+1}\gamma(x)m_1(p)q\cdot dv + (-1)^{|p|}s_i(pm(x)q\cdot dv) \\
& + (-1)^{|x|}\gamma(x)pm(q)\cdot dv + (-1)^{|p|+|x|+|q|+1}s_i(pxq\cdot dm(v)), \\
= & (-1)^{|p|}\gamma(x)m(p)q\cdot dv + (-1)^{|q|+1}\gamma(x)pq\cdot dm(v) + (-1)^{|p|+1}\gamma(x)m_1(p)q\cdot dv \\
& + (-1)^{|p|}s_i(pm(x)q\cdot dv) + (-1)^{|p|+|x|+|q|+1}s_i(pxq\cdot dm(v)). \\
\end{split}
\end{equation}

Suppose that $x,v \in \boldsymbol{t}:=\{t_i\}_{i \in I}$, then \eqref{eqn_htpretdummya} implies that
\begin{equation} \label{eqn_htpretdummyb}
\begin{split}
h_i(\alpha) & = (-1)^{|p|}s_i(p[\tau,x]q\cdot dv) + (-1)^{|p|+|x|+|q|+1}s_i(pxq\cdot d[\tau,v]), \\
& = \alpha + (-1)^{|p|+|x|+|q|}s_i([pxq,v]\cdot d\tau + [pxq,\tau]\cdot dv), \\
& = \alpha. \\
\end{split}
\end{equation}
A similar calculation shows that $h_i(\alpha)=\alpha$ when $x \in \boldsymbol{t}$ and $v = \tau$. This means that $h_i$ acts as the identity on $i+1$-normalised 1-forms.

Now suppose that $x = \tau$ and $v \in \boldsymbol{t}$, then
\begin{equation} \label{eqn_htpretdummyc}
\begin{split}
h_i(\alpha) = & (-1)^{|p|}[\tau,p]q\cdot dv + (-1)^{|p|}m_1(p)q\cdot dv + (-1)^{|q|+1}pq\cdot dm(v) + (-1)^{|p|+1}m_1(p)q\cdot dv \\
& + (-1)^{|p|}s_i(p\tau^2q\cdot dv) + (-1)^{|p|+|q|}s_i(p\tau q\cdot dm(v)) + (\text{$i+1$-normalised 1-forms}), \\
= & (-1)^{|p|}\tau pq\cdot dv - \alpha + (-1)^{|q|+1}pq\cdot dm(v) + \alpha \\
& + (-1)^{|p|+|q|}s_i(p\tau q\cdot dm(v)) + (\text{$i+1$-normalised 1-forms}), \\
= & (-1)^{|p|}\tau pq\cdot dv + (-1)^{|q|+1}pq\cdot d[\tau,v] \\
& + (-1)^{|p|+|q|}s_i(p\tau q\cdot d[\tau,v]) + (\text{$i+1$-normalised 1-forms}), \\
= & (-1)^{|p|}\tau pq\cdot dv + (-1)^{|q|}pqv\cdot d\tau + (-1)^{|q|}[pq,\tau]\cdot dv \\
& + (-1)^{|q|+1}pqv\cdot d\tau + (-1)^{|q|+1}pq\tau\cdot dv + (\text{$i+1$-normalised 1-forms}), \\
= & 0 + (\text{$i+1$-normalised 1-forms}). \\
\end{split}
\end{equation}
A similar calculation shows that
\[ h_i(\alpha) = 0 \mod (\text{$i+1$-normalised 1-forms}) \]
when $x = \tau$ and $v = \tau$. This means that $h_i$ takes $i$-normalised forms to $i+1$-normalised forms as claimed.

We define the map $H:\drof[Ass]{\ctalg{V}} \to \drnof[Ass]{\ctalg{V}}$ as
\[ H:= \ldots \circ h_n\circ\ldots\circ h_2\circ h_1. \]
By the definition of the $h_n$'s the map $H$ is homotopic to $\id$. Equations \eqref{eqn_htpretdummya}, \eqref{eqn_htpretdummyb} and \eqref{eqn_htpretdummyc} show that $H$ takes 1-forms to normalised 1-forms so $H$ is a well defined map. Lastly equation \eqref{eqn_htpretdummyb} shows that $H$ splits the inclusion
\[ \drnof[Ass]{\ctalg{V}} \hookrightarrow \drof[Ass]{\ctalg{V}}, \]
so $H$ is the chain homotopy retraction exhibiting \drnof[Ass]{\ctalg{V}} as a chain homotopy retract of \drof[Ass]{\ctalg{V}} that we sought.
\end{proof}

\subsection{Cyclic cohomology theories}

We will now define the cyclic cohomology theories for \li, {\ai} and \ci-algebras. In \cite{pensch} Penkava and Schwarz defined the cyclic cohomology of an \ai-algebra and showed that it controlled the deformations of the \ai-infinity algebra which preserve a fixed invariant inner product. Our approach to defining cyclic cohomology will be different however, in that we will use the framework of noncommutative geometry.

Let $V$ be a profinite graded module. We say $q(\boldsymbol{t}) \in \ctalg{V}$ vanishes at zero if it belongs to the ideal $\Sigma V^* \cdot \ctalg{V}$ of \ctalg{V}, i.e. the power series $q(\boldsymbol{t})$ has vanishing constant term.

\begin{defi}
Let $V$ be a free graded module:
\begin{enumerate}
\item[(a)]
Let $m:\csalg{V} \to \csalg{V}$ be an \li-structure. The cyclic Chevalley-Eilenberg complex of the \li-algebra $V$ is defined on the module of 0-forms vanishing at zero;
\[ \cclac{V}:=\Sigma\left[\{q \in \drzf[Com]{\csalg{V}}:\ q \text{ vanishes at zero} \}\right]. \]
The differential on this complex is the restriction of the (suspension of the) Lie operator of the vector field $m$;
\[ L_m:\drzf[Com]{\csalg{V}} \to \drzf[Com]{\csalg{V}} \]
to \cclac{V}. The cyclic Chevalley-Eilenberg cohomology of $V$ is defined as the cohomology of the complex \cclac{V} and is denoted by \hclac{V}.
\item[(b)]
Let $m:\ctalg{V} \to \ctalg{V}$ be an \ai-structure. The cyclic Hochschild complex of the \ai-algebra $V$ is defined on the module of 0-forms vanishing at zero:
\[ \cchoch{V}:=\Sigma\left[\{q \in \drzf[Ass]{\ctalg{V}}:\ q \text{ vanishes at zero} \}\right]. \]
The differential on this complex is the restriction of the (suspension of the) Lie operator of the vector field $m$;
\[ L_m:\drzf[Ass]{\ctalg{V}} \to \drzf[Ass]{\ctalg{V}} \]
to \cchoch{V}. The cyclic Hochschild cohomology of $V$ is defined as the cohomology of the complex \cchoch{V} and is denoted by \hchoch{V}.
\item[(c)]
Let $m:\clalg{V} \to \clalg{V}$ be a \ci-structure. The cyclic Harrison complex of the \ci-algebra $V$ is defined on the module consisting of all 0-forms:
\[ \ccaq{V}:=\Sigma\drzf[Lie]{\clalg{V}}. \]
The differential on this complex is the (suspension of the) Lie operator of the vector field $m$;
\[ L_m:\drzf[Lie]{\clalg{V}} \to \drzf[Lie]{\clalg{V}}. \]
The cyclic Harrison cohomology of $V$ is defined as the cohomology of the complex \ccaq{V} and is denoted by \hcaq{V}.
\end{enumerate}
\end{defi}

\begin{rem}
Note that a consequence of the above definition is that the cyclic Chevalley-Eilenberg complex is just the Chevalley-Eilenberg complex with \emph{trivial coefficients}.
\end{rem}

It will be useful for the purposes of the next section to describe equivalent formulations for the cyclic cohomology of an \ai-algebra. For this we will need the following lemmas:

\begin{lemma} \label{lem_difcom}
Let $V$ be an \ai-algebra and consider the maps $b'$ and $b$ defined by diagrams \eqref{fig_barhom} and \eqref{fig_chomid} respectively:
\begin{enumerate}
\item[(i)]
\[ bN=-Nb'. \]
\item[(ii)]
\[ b'(1-z)=-(1-z)b. \]
\end{enumerate}
\end{lemma}

\begin{proof}
\
\begin{enumerate}
\item[(i)]
This is a tautological consequence of Definition \ref{def_barhom}, Lemma \ref{lem_nrmdif} and Lemma \ref{lem_schf} part (v).
\item[(ii)]
Let $x \in \Sigma V^*$ and $y \in (\Sigma V^*)^{\cotimes n}$ for $n \geq 0$:
\begin{displaymath}
\begin{split}
b'(1-z)xy & = b'([x,y]) = [m(x),y] + (-1)^{|x|}[x,m(y)]. \\
(1-z)b(xy) & =(1-z)\thiso[Ass]^{-1}L_m(dx\cdot y), \\
& = -(1-z)\thiso[Ass]^{-1}(dm(x)\cdot y) - (-1)^{|x|}(1-z)\thiso[Ass]^{-1}(dx\cdot m(y)). \\
\end{split}
\end{displaymath}
By a calculation similar to that of equation \eqref{eqn_nrmdif} we could deduce that
\[ \thiso[Ass]^{-1}(d(x_1\cotimes\ldots\cotimes x_i)\cdot y) = (1+z+\ldots+z^{i-1}) \cdot x_1\cotimes\ldots\cotimes x_i \cdot y \]
and since $(1-z)(1+z+\ldots+z^{i-1})=1-z^i$ we use equation \eqref{eqn_cyccom} to conclude that
\[ (1-z)b(xy) = -\left([m(x),y] + (-1)^{|x|}[x,m(y)]\right), \]
therefore $b'(1-z)=-(1-z)b$ as claimed.
\end{enumerate}
\end{proof}

We can now define a new bicomplex computing the cyclic cohomology of an \ai-algebra:

\begin{defi} \label{def_tsygan}
Let $V$ be an \ai-algebra and consider the maps $b'$ and $b$ defined by diagrams \eqref{fig_barhom} and \eqref{fig_chomid} respectively. The Tsygan bicomplex \ctsygan{V} is the bicomplex;
\[ \xymatrix{\left(\Sigma\left[\col{V}\right],b\right) \ar^{1-z}[r] & \left(\Sigma\left[\col{V}\right],b'\right) \ar^{N}[r] & \left(\Sigma\left[\col{V}\right],b\right) \ar^-{1-z}[r] & \ldots}. \]
An element has homogeneous bidegree $(i,j)$ if it is in the $i$th column from the left (where the leftmost column has bidegree $(0,\bullet)$) and has degree $j$ in the graded profinite module $\Sigma\left[\col{V}\right]$. The Tsygan cohomology of the \ai-algebra $V$ is defined as the total cohomology of the bicomplex \ctsygan{V} \emph{formed by taking direct sums}, that is to say that it is the cohomology of the complex
\[ C:= \prod_{n \in \mathbb{Z}}\left(\bigoplus_{i+j=n}C^{ij}\right) \]
where $C^{ij}$ is the component of \ctsygan{V} of bidegree $(i,j)$.
\end{defi}

\begin{rem}
Note that the odd numbered columns are copies of the Bar complex \cbr{V} whilst the even numbered columns are isomorphic to the Hochschild complex \choch{V}{V^*}.
\end{rem}

\begin{rem}
Since the underlying module $V$ is $\mathbb{Z}$-graded, the direct product totalisation of the Tsygan complex might a priori lead to different cohomology. The direct product totalisation has the advantage that the Tsygan cohomology defined in this way is homotopy invariant in the sense that two weakly equivalent $A_\infty$-algebras will have isomorphic Tsygan cohomology. This follows easily from the standard spectral sequence arguments (which fail for the direct sum totalisation). On the other hand, our definition is clearly homotopy invariant for \emph{minimal} $A_\infty$-algebras since in this case a weak equivalence implies an isomorphism. We will see later on that, rather surprisingly, our definition turns out to be homotopy invariant after all for $C_\infty$-algebras. This is a consequence of the Hodge decomposition for the cyclic cohomology of a $C_\infty$-algebra.
\end{rem}

\begin{rem} \label{rem_connes}
It is also possible to define a Connes and a normalised Connes complex generalising the construction of the Connes and normalised Connes complex for strictly associative algebras described in 2.1.7 and 2.1.9 of \cite{loday}. This is done using Lemma \ref{lem_unitch} and Proposition \ref{prop_htpret} in exactly the same manner as it is performed in \cite{loday}.
\end{rem}

Let us now show that Tsygan cohomology computes the same cohomology as cyclic Hochschild cohomology which we defined earlier:

\begin{lemma} \label{lem_cyciso}
Let $V$ be an \ai-algebra:
\begin{enumerate}
\item[(i)]
The map $b':\col{V} \to \col{V}$ defined by diagram \eqref{fig_barhom} lifts uniquely to a map
\[ b' : \prod_{i=1}^\infty \left((\Sigma V^*)^{\cotimes i}\right)_{Z_i} \to \prod_{i=1}^\infty \left((\Sigma V^*)^{\cotimes i}\right)_{Z_i}. \]
This gives rise to the following identity:
\[ \cchoch{V} = \left(\Sigma\left[\prod_{i=1}^\infty \left((\Sigma V^*)^{\cotimes i}\right)_{Z_i}\right],b'\right). \]
\item[(ii)]
There is a quasi-isomorphism of complexes
\[ \left(\Sigma\left[\prod_{i=1}^\infty \left((\Sigma V^*)^{\cotimes i}\right)_{Z_i}\right],b'\right) \overset{\sim}{\to} \ctsygan{V} \]
which is given by mapping $\Sigma\left[\prod_{i=1}^\infty \left((\Sigma V^*)^{\cotimes i}\right)_{Z_i}\right]$ onto the leftmost column of \ctsygan{V} by
\[ \Sigma\left[\prod_{i=1}^\infty \left((\Sigma V^*)^{\cotimes i}\right)_{Z_i}\right] \overset{N}{\longrightarrow} \Sigma\left[\col{V}\right] = CC^{0,\bullet}_\mathrm{Tsygan}(V) \subset \ctsygan{V}. \]
\end{enumerate}
\end{lemma}

\begin{proof}
\
\begin{enumerate}
\item[(i)]
Recall that
\[ \drzf[Ass]{\ctalg{V}} = \ctalg{V}/[\ctalg{V},\ctalg{V}]. \]
The claim then follows tautologically from equation \eqref{eqn_cyccom} and the definitions.
\item[(ii)]
This map is a quasi-isomorphism because
\[ \xymatrix{ \left(\Sigma\left[\prod_{i=1}^\infty \left((\Sigma V^*)^{\cotimes i}\right)_{Z_i}\right],b'\right) \ar^-{N}[r] & \ctsygan{V}} \]
is a resolution of $\left(\Sigma\left[\prod_{i=1}^\infty \left((\Sigma V^*)^{\cotimes i}\right)_{Z_i}\right],b'\right)$.
\end{enumerate}
\end{proof}

Following Lemma \ref{lem_cyciso} we will denote the cohomology of the Tsygan complex by \hchoch{V}. We will now describe the periodicity long exact sequence linking Hochschild cohomology with cyclic cohomology:

\begin{prop} \label{prop_sbiseq}
Let $V$ be a unital \ai-algebra, then we have the following long exact sequence in cohomology:
\[ \xymatrix{ \ldots \ar^-{B}[r] & \hchoch[n-2]{V} \ar^{S}[r] & \hchoch[n]{V} \ar^{I}[r] & \hhoch[n]{V}{V^*} \ar^{B}[r] & \hchoch[n-1]{V} \ar^-{S}[r] & \ldots \\ }. \]
\end{prop}

\begin{proof}
There is a short exact sequence of complexes
\[ \xymatrix{ 0 \ar[r] & \Sigma^{-2}\ctsygan{V} \ar^{S}[r] & \ctsygan{V} \ar^{I}[r] & CC\{2\}^{\bullet\bullet}_{\mathrm{Tsygan}}(V) \ar[r] & 0 \\ } \]
where $CC\{2\}^{\bullet\bullet}_{\mathrm{Tsygan}}(V)$ is the total complex of the first two columns of \ctsygan{V};
\[ CC\{2\}^{\bullet\bullet}_{\mathrm{Tsygan}}(V) := \xymatrix{\left(\Sigma\left[\col{V}\right],b\right) \ar^{1-z}[r] & \left(\Sigma\left[\col{V}\right],b'\right)}. \]

By Lemma \ref{lem_unitch} the second column of $CC\{2\}^{\bullet\bullet}_{\mathrm{Tsygan}}(V)$ is acyclic, so $CC\{2\}^{\bullet\bullet}_{\mathrm{Tsygan}}(V)$ is quasi-isomorphic to its first column which is just the complex \choch{V}{V^*} by Remark \ref{rem_chomid}. The long exact sequence in cohomology of the proposition is derived from the above short exact sequence of complexes.
\end{proof}

\begin{rem} \label{rem_sbiseq}
We can describe the connecting map
\[ B: \hhoch[n]{V}{V^*} \to \hchoch[n-1]{V} \]
more explicitly as follows: First of all as mentioned in Proposition \ref{prop_sbiseq}, the projection
\[ \pi: CC\{2\}^{\bullet\bullet}_{\mathrm{Tsygan}}(V) \to \left(\Sigma\left[\col{V}\right],b\right) \]
of $CC\{2\}^{\bullet\bullet}_{\mathrm{Tsygan}}(V)$ onto the first column is a quasi-isomorphism. A simple check reveals that this map has a section (and hence a quasi-inverse),
\[ i:\left(\Sigma\left[\col{V}\right],b\right) \to CC\{2\}^{\bullet\bullet}_{\mathrm{Tsygan}}(V) \]
which is defined as follows;
\begin{displaymath}
\begin{array}{ccc}
\col{V} & \to & \col{V} \oplus \col{V}, \\
x & \mapsto & x \oplus -h(1-z)[x]; \\
\end{array}
\end{displaymath}
where $h:\col{V} \to \col{V}$ is the contracting homotopy of Lemma \ref{lem_unitch}. It follows from the definitions and a simple check that the connecting map $B$ is the map induced by the following map of complexes:
\begin{equation} \label{eqn_sbiseq}
\begin{array}{ccc}
\left(\Sigma\left[\col{V}\right],b\right) & \to & \left(\Sigma\left[\prod_{i=1}^\infty\left((\Sigma V^*)^{\cotimes i}\right)_{Z_i}\right],b'\right), \\
x & \mapsto & -h(1-z)[x]. \\
\end{array}
\end{equation}

A similar check also shows that the map $I:\hchoch[n]{V} \to \hhoch[n]{V}{V^*}$ is the map induced by the following map of complexes:
\begin{displaymath}
\begin{array}{ccc}
\left(\Sigma\left[\prod_{i=1}^\infty\left((\Sigma V^*)^{\cotimes i}\right)_{Z_i}\right],b'\right) & \to & \left(\Sigma\left[\col{V}\right],b\right), \\
x & \mapsto & N\cdot x. \\
\end{array}
\end{displaymath}
\end{rem}
\section{The Hodge Decomposition of Hochschild and Bar Cohomology} \label{sec_hdgdec}

In this section we will be concerned with constructing the Hodge decomposition of the Hochschild and bar cohomology of a \ci-algebra. Given a \ci-algebra $V$ we will construct both the Hodge decomposition of the Hochschild cohomology of $V$ with coefficients in $V^*$ as well as the Hodge decomposition of the bar cohomology of $V$.

The Hodge decomposition of the Hochschild (co)homology of a commutative algebra has been described by many authors such as \cite{gersch}, \cite[\S4.5,\S4.6]{loday} and \cite{natsch}. Our approach will be to determine a Hodge decomposition for \ci-algebras and this will naturally include a Hodge decomposition of the Hochschild cohomology of a commutative algebra. Our results however are more than just a mere generalisation of the results contained in \cite{loday} and \cite{natsch}. Considering the Hodge decomposition in the broader perspective of \ci-algebras leads us to use the framework of noncommutative geometry. This is a very natural setting in which to construct the Hodge decomposition and leads to a much more conceptual approach and results in a streamlining of the calculations. The pay off is that we are able to obtain new results, even for the Hodge decomposition of commutative algebras.

Let $W$ be a profinite graded module. For the rest of this section we will denote the canonical associative multiplication on $\widehat{T}W$ by $\mu$. The cocommutative comultiplication (\emph{shuffle product}) on $\widehat{T}W$ defined by Remark \ref{rem_ccmult} will be denoted by $\Delta$. This gives $\widehat{T}W$ the structure of a Hopf algebra. 

The Hodge decomposition will be constructed from a spectral decomposition of an operator which we will refer to as the \emph{modified shuffle operator} and which we will now define. The original `shuffle operator' was used in \cite{gersch} and \cite{natsch} and seems to have originally been introduced by Barr in \cite{barr}.

\begin{defi} \label{def_shufop}
Let $W$ be a profinite graded module. The \emph{modified shuffle operator} $s$ is defined as
\[ s:= \mu\Delta:\widehat{T}W \to \widehat{T}W. \]
We denote its components by
\[ s_n:W^{\cotimes n} \to W^{\cotimes n}. \]

We also define a second operator $\tilde{s}:\prod_{i=1}^\infty W^{\cotimes i} \to \prod_{i=1}^\infty W^{\cotimes i}$ by the commutative diagram;
\begin{displaymath}
\xymatrix{ W\cotimes\widehat{T}W \ar^{1\cotimes s}[r] \ar@{=}[d] & W\cotimes\widehat{T}W \ar@{=}[d] \\ \prod_{i=1}^\infty W^{\cotimes i} \ar^{\tilde{s}}[r] & \prod_{i=1}^\infty W^{\cotimes i}}
\end{displaymath}
We denote its components by
\[ \tilde{s}_n:W^{\cotimes n+1} \to W^{\cotimes n+1}. \]
\end{defi}

Next we will define a filtration of $\widehat{T}W$ which comes from the observation $\widehat{T}W = \widehat{\mathcal{U}}(\widehat{L}W)$. This will assist us in our calculations later:

\begin{defi} \label{def_pbwfil}
Let $W$ be a profinite graded module. The increasing filtration $\{F_p(\widehat{T}W)\}_{p=0}^\infty$ of $\widehat{T}W$ is defined as follows: $F_p(\widehat{T}W)$ is defined as the submodule of $\widehat{T}W$ which consists of all (possibly uncountably infinite) linear combinations of elements of the form
\[ g_1 \ldots g_i, \]
where $g_1,\ldots,g_i \in \widehat{L}W$ and $0 \leq i \leq p$. By convention $F_0(\widehat{T}W)= \gf$.
\end{defi}

\begin{rem}
Note that this filtration is not exhaustive but that $W^{\cotimes n} \subset F_n(\widehat{T}W)$. Also note that this filtration is preserved by the modified shuffle operator $s$. This follows from the fact that $\widehat{L}W$ consists of primitive elements of the Hopf algebra $(\widehat{T}W,\mu,\Delta)$.
\end{rem}

Now we will determine the spectral decomposition of the modified shuffle operator $s$. We will do this by exhibiting a certain polynomial that annihilates it. Let us define the polynomial $\nu_n(X) \in \mathbb{Z}(X)$  for $n \geq 0$ by the formula
\begin{equation} \label{eqn_poldef}
\nu_n(X):=\prod_{i=0}^n (X-\lambda_i), \quad \lambda_i:=2^i.
\end{equation}

\begin{lemma} \label{lem_shuann}
Let $W$ be a profinite graded module. For all $n \geq 0$;
\[ \nu_n(s_n)=\prod_{i=0}^n (s_n-\lambda_i\id)=0. \]
\end{lemma}

\begin{proof}
Let us prove the following equation:
\begin{equation} \label{eqn_shufev}
s_n(x)=\lambda_p x \mod F_{p-1}(\widehat{T}W), \quad \text{for all } x \in F_p(\widehat{T}W)\cap W^{\cotimes n}.
\end{equation}
Since the modified shuffle operator $s$ preserves the filtration, we may assume that $x$ is a linear combination of elements of the form $g_1\ldots g_p$ where $g_i \in \widehat{L}W$.

We calculate that for all $g_1,\ldots,g_p \in \widehat{L}W$;
\begin{displaymath}
\begin{split}
s_n(g_1\ldots g_p) & = \mu\Delta(g_1\ldots g_p), \\
& =\mu\left[(g_1\cotimes 1 + 1 \cotimes g_1)\ldots(g_p\cotimes 1 + 1 \cotimes g_p)\right], \\
& =2^p g_1\ldots g_p \mod F_{p-1}(\widehat{T}W). \\
\end{split}
\end{displaymath}
The last equality follows since by using commutators we can transform
\[ \mu\left[(g_1\cotimes 1 + 1 \cotimes g_1)\ldots(g_p\cotimes 1 + 1 \cotimes g_p)\right] \]
into $2^p g_1\ldots g_p$. By definition, these commutating elements are in $F_{p-1}(\widehat{T}W)$. This calculation implies equation \eqref{eqn_shufev} from which the lemma follows as a trivial consequence.
\end{proof}

\begin{rem}
Note $\nu_n(X)$ is not the minimal polynomial for $s_n$. The minimal polynomial of $s_n$ was constructed in \cite{gersch} but will not be necessary for our purposes.
\end{rem}

\begin{defi} \label{def_idemdf}
Let $W$ be a profinite graded module. We define a family of operators
\begin{displaymath}
\begin{array}{ll}
e(i):\widehat{T}W \to \widehat{T}W, & i \geq 0; \\
e_n(i):W^{\cotimes n} \to W^{\cotimes n}; \\
\end{array}
\end{displaymath}
as the Lagrange interpolation polynomials of the operator $s_n$:
\begin{displaymath}
e_n(i):=\left\{
\begin{array}{ccc}
\left[ \underset{\begin{subarray}{c} 0 \leq r \leq n \\ r \neq i \end{subarray}}{\prod} (\lambda_i - \lambda_r) \right]^{-1} \underset{\begin{subarray}{c} 0 \leq r \leq n \\ r \neq i \end{subarray}}{\prod} (s_n - \lambda_r\id) & , & 0 \leq i \leq n; \\
0 & , & i > n. \\
\end{array}
\right.
\end{displaymath}

We also define a family of operators
\[ \tilde{e}(i):\prod_{i=1}^\infty W^{\cotimes i} \to \prod_{i=1}^\infty W^{\cotimes i} \quad , \quad i \geq 0 \]
by the commutative diagram
\begin{displaymath}
\xymatrix{W\cotimes\widehat{T}W \ar^{1\cotimes e(i)}[r] \ar@{=}[d] & W\cotimes\widehat{T}W \ar@{=}[d] \\ \prod_{i=1}^\infty W^{\cotimes i} \ar^{\tilde{e}(i)}[r] & \prod_{i=1}^\infty W^{\cotimes i}}
\end{displaymath}
and denote their components by
\[ \tilde{e}_n(i):W^{\cotimes n+1} \to W^{\cotimes n+1}. \]
\end{defi}

\begin{lemma} \label{lem_shidem}
Let $W$ be a profinite graded module. We have the following identities:
\begin{displaymath}
\begin{array}{rl}
\textnormal{(a)} & s_n=\lambda_0 e_n(0) + \ldots + \lambda_n e_n(n). \\
\textnormal{(a*)} & s=\sum_{i=0}^\infty \lambda_i e(i). \\
\textnormal{(b)} & \id_n=e_n(0) + \ldots + e_n(n). \\
\textnormal{(b*)} & \id = \sum_{i=0}^\infty e(i). \\
\textnormal{(c)} & e_n(i) \circ e_n(j)=\left\{
\begin{array}{ccc}
e_n(i) & , & i=j; \\
0 & , & i \neq j. \\
\end{array}
\right.\\
\textnormal{(c*)} & e(i) \circ e(j)=\left\{
\begin{array}{ccc}
e(i) & , & i=j; \\
0 & , & i \neq j. \\
\end{array}
\right. \\
\end{array}
\end{displaymath}
\end{lemma}

\begin{proof}
(a), (b) and (c) are a formal consequence of Lemma \ref{lem_shuann}. Since we assume that our ground ring {\gf} contains the field $\mathbb{Q}$ and since the polynomial $\nu_n(X) \in \mathbb{Z}(X)$ defined by equation \eqref{eqn_poldef} annihilates $s_n$ and has no repeated roots, elementary linear algebra implies that $s_n$ is diagonalisable. The map $e(i)$ is the projection onto the eigenspace
\[ \{x \in \widehat{T}W : s(x)=\lambda_i x \}. \]
(a*), (b*) and (c*) are trivial consequences of (a), (b) and (c) respectively.
\end{proof}

\begin{rem}
Since $\tilde{e}(i):=1\cotimes e(i)$, the same identities hold when we replace $e(i)$ with $\tilde{e}(i)$, $s$ with $\tilde{s}$ and $\id_n$ with $\id_{n+1}$ in the above Lemma.
\end{rem}

We want to use the spectral decomposition of the modified shuffle operator $s$ to construct a decomposition of the relevant cohomology theories. For this we need the following lemma:

\begin{lemma} \label{lem_shucom}
Let $V$ be a \ci-algebra and consider the maps $b'$ and $b$ defined by diagrams \eqref{fig_barhom} and \eqref{fig_chomid} respectively. Also consider the maps $s,\tilde{s}:\col{V} \to \col{V}$ defined in Definition \ref{def_shufop}. We have the following identities:
\begin{enumerate}
\item[(i)]
\[ s \circ b'=b' \circ s. \]
\item[(ii)]
\[ \tilde{s} \circ b=b \circ \tilde{s}. \]
\end{enumerate}
\end{lemma}

\begin{proof}
Let $m:\clalg{V} \to \clalg{V}$ be the \ci-structure on $V$. Recall from Remark \ref{rem_calseq} that this could be considered to be an \ai-structure $m:\ctalg{V} \to \ctalg{V}$ which is also a Hopf algebra derivation of the Hopf algebra $(\ctalg{V},\mu,\Delta)$:
\begin{enumerate}
\item[(i)]
This is just a trivial consequence of the definition of $b'$ and the fact that $m$ is a derivation of the Hopf algebra $(\ctalg{V},\mu,\Delta)$:
\begin{equation} \label{eqn_shucom}
\mu\Delta m=\mu(m\cotimes 1+1\cotimes m)\Delta=m\mu\Delta.
\end{equation}
Since $b'$ is the restriction of $m$ to \col{V} we obtain $s \circ b'=b' \circ s$
\item[(ii)]
Let $\thiso[Ass]:\col{V} \to \drof[Ass]{\ctalg{V}}$ be the module isomorphism we defined in Lemma \ref{lem_ofisom}. Let $x \in \Sigma V^*$ and $y \in (\Sigma V^*)^{\cotimes n}$ for $n \geq 0$;
\begin{displaymath}
\begin{split}
b\tilde{s}(xy) & = b(xs(y)) = \thiso[Ass]^{-1}L_m(dx\cdot s(y)), \\
& = -\thiso[Ass]^{-1}(dm(x)\cdot s(y))-(-1)^{|x|}x\cdot sm(y). \\
\end{split}
\end{displaymath}
where the computation of the last term follows from equation \eqref{eqn_shucom}. Furthermore;
\begin{displaymath}
\begin{split}
\tilde{s}b(xy) & = \tilde{s}\thiso[Ass]^{-1}L_m(dx\cdot y), \\
& = -\tilde{s}\thiso[Ass]^{-1}(dm(x)\cdot y)-(-1)^{|x|}x\cdot sm(y). \\
\end{split}
\end{displaymath}

We would like to show that given any $u \in \clalg{V}$ and $w \in \ctalg{V}$ we have the following identity:
\begin{equation} \label{eqn_shucomdummy}
\thiso[Ass]^{-1}(du\cdot s(w))=\tilde{s}\thiso[Ass]^{-1}(du\cdot w).
\end{equation}
Since $m$ maps $\Sigma V^*$ to \clalg{V} it will follow from equation \eqref{eqn_shucomdummy} that $\tilde{s} \circ b=b \circ \tilde{s}$. Since \clalg{V} is generated by Lie monomials and \eqref{eqn_shucomdummy} is a tautology for a Lie monomial $u$ of order one (i.e. $u \in \Sigma V^*$), we can proceed by induction and assume there exists $u_1,u_2 \in \clalg{V}$ such that $u=[u_1,u_2]$ and such that \eqref{eqn_shucomdummy} holds for $u_1$ and $u_2$:
\begin{displaymath}
\begin{split}
du\cdot s(w) & = d[u_1,u_2]\cdot s(w) =[du_1,u_2]\cdot s(w) - (-1)^{|u_1||u_2|}[du_2,u_1] \cdot s(w), \\
& = du_1\cdot [u_2,s(w)] - (-1)^{|u_1||u_2|}du_2\cdot [u_1,s(w)]. \\
du \cdot w & = d[u_1,u_2]\cdot w = du_1\cdot [u_2,w] - (-1)^{|u_1||u_2|}du_2\cdot [u_1,w]. \\
\end{split}
\end{displaymath}
By the inductive hypothesis we obtain;
\begin{displaymath}
\begin{split}
\tilde{s}\thiso[Ass]^{-1}(du\cdot w) & = \tilde{s}\thiso[Ass]^{-1}(du_1\cdot [u_2,w]) - (-1)^{|u_1||u_2|}\tilde{s}\thiso[Ass]^{-1}(du_2\cdot [u_1,w]), \\
& = \thiso[Ass]^{-1}(du_1\cdot s([u_2,w])) - (-1)^{|u_1||u_2|}\thiso[Ass]^{-1}(du_2\cdot s([u_1,w])). \\
\end{split}
\end{displaymath}

In order to complete the proof we will need one final auxiliary calculation. Recall that \clalg{V} coincides with the Lie subalgebra of primitive elements of the Hopf algebra $(\ctalg{V},\mu,\Delta)$. Let $a \in \clalg{V}$ and $b \in \ctalg{V}$;
\begin{equation} \label{eqn_shucomdummya}
\begin{split}
s([a,b]) & = \mu[a\cotimes 1+1\cotimes a,\Delta b], \\
& = \mu(a\cotimes 1\cdot\Delta b -(-1)^{|a||b|}\Delta b \cdot 1\cotimes a) + \mu(1\cotimes a\cdot\Delta b - (-1)^{|a||b|}\Delta b \cdot a\cotimes 1), \\
& = a\cdot\mu\Delta b - (-1)^{|a||b|}\mu\Delta b \cdot a - 0, \\
& = [a,s(b)]. \\
\end{split}
\end{equation}
Finally we can establish equation \eqref{eqn_shucomdummy} and hence conclude the proof:
\begin{displaymath}
\begin{split}
\thiso[Ass]^{-1}(du\cdot s(w)) & = \thiso[Ass]^{-1}(du_1\cdot [u_2,s(w)]) - (-1)^{|u_1||u_2|}\thiso[Ass]^{-1}(du_2\cdot [u_1,s(w)]), \\
& = \thiso[Ass]^{-1}(du_1\cdot s([u_2,w])) - (-1)^{|u_1||u_2|}\thiso[Ass]^{-1}(du_2\cdot s([u_1,w])), \\
& = \tilde{s}\thiso[Ass]^{-1}(du\cdot w). \\
\end{split}
\end{displaymath}
\end{enumerate}
\end{proof}

\begin{cor} \label{cor_shucom}
Let $V$ be a \ci-algebra and consider the maps
\[ e(i),\tilde{e}(i):\col{V}\to\col{V} \]
defined by Definition \ref{def_idemdf}. We have the following identities:
\begin{enumerate}
\item[(i)]
For all $i \geq 1$,
\[ e(i) \circ b' = b' \circ e(i). \]
\item[(ii)]
For all $i \geq 0$,
\[ \tilde{e}(i) \circ b = b \circ \tilde{e}(i). \]
\end{enumerate}
\end{cor}

\begin{proof}
This is just a formal consequence of Lemma \ref{lem_shucom} and Lemma \ref{lem_shidem}. $e(i)$ is just the projection onto the eigenspace
\[ \{ x \in \col{V} : s(x)=\lambda_i x \} \]
whilst $\tilde{e}(i)$ is the projection onto the eigenspace
\[ \{ x \in \col{V} : \tilde{s}(x)=\lambda_i x \}. \]
Lemma \ref{lem_shucom} then tells us that $b'$ and $b$ preserve these eigenspaces respectively and hence commute with the projections $e(i)$ and $\tilde{e}(i)$ respectively.
\end{proof}

We are now in a position to state the main theorem of this section; the Hodge decomposition of the Hochschild and bar cohomology of a \ci-algebra:

\begin{theorem}
Let $V$ be a \ci-algebra:
\begin{enumerate}
\item[(i)]
The bar complex \cbr{V} of the \ci-algebra $V$ splits as the direct product of the subcomplexes $\Sigma\left(e(i)\left[\col{V}\right]\right)$:
\[ \cbr{V}=\prod_{i=1}^{\infty}\left(\Sigma\left(e(i)\left[\col{V}\right]\right),b'\right). \]
\item[(ii)]
The Hochschild complex of the \ci-algebra $V$ with coefficients in $V^*$ splits as the direct product of the subcomplexes $\Sigma\left(\tilde{e}(i)\left[\col{V}\right]\right)$:
\[ \choch{V}{V^*}=\prod_{i=0}^\infty\left(\Sigma\left(\tilde{e}(i)\left[\col{V}\right]\right),b\right). \]
\end{enumerate}
\end{theorem}

\begin{proof}
\
\begin{enumerate}
\item[(i)]
This is just a trivial consequence of Lemma \ref{lem_shidem} and Corollary \ref{cor_shucom}. Lemma \ref{lem_shidem} tells us that the module \col{V} splits as a product of submodules;
\[ \col{V}=\prod_{i=1}^\infty e(i)\left[\col{V}\right]. \]
Corollary \ref{cor_shucom} part (i) tells us that when we equip \col{V} with the differential $b'$, the modules $e(i)\left[\col{V}\right]$ are actually subcomplexes of $\left(\col{V},b'\right)$.
\item[(ii)]
Lemma \ref{lem_shidem} tells us that the module \col{V} splits as a product of submodules;
\[ \col{V}=\prod_{i=0}^\infty \tilde{e}(i)\left[\col{V}\right]. \]
Corollary \ref{cor_shucom} part (ii) tells us that when we equip \col{V} with the differential $b$ the modules $\tilde{e}(i)\left[\col{V}\right]$ are actually subcomplexes of $\left(\col{V},b\right)$. By Remark \ref{rem_chomid} the complex \choch{V}{V^*} is isomorphic to $\left(\Sigma\left[\col{V}\right],b\right)$.
\end{enumerate}
\end{proof}

\begin{defi}
Let $V$ be a \ci-algebra. We will define the complex \echoch{V}{V^*}{i} as the subcomplex of \choch{V}{V^*} given by the formula;
\[ \echoch{V}{V^*}{i}:= \left(\Sigma\left(\tilde{e}(i)\left[\col{V}\right]\right),b\right), \quad i \geq 0 \]
and denote its cohomology by \ehhoch{V}{V^*}{i}.
\end{defi}

Let $l: \drof[Lie]{\clalg{V}} \hookrightarrow \drof[Ass]{\ctalg{V}}$ be the map defined by equation \eqref{eqn_liasmp}. Suppose that $V$ is a \ci-algebra, then by Lemma \ref{lem_plcomm} we see that $l$ is a map of complexes;
\[ l:\caq{V}{V^*} \to \choch{V}{V^*}. \]
We shall now show that this map is in fact split:

\begin{prop} \label{prop_aqsplt}
Let $V$ be a \ci-algebra. There is an isomorphism
\[ \caq{V}{V^*} \cong \echoch{V}{V^*}{1}. \]
In particular the map $l:\caq{V}{V^*} \to \choch{V}{V^*}$ is split.
\end{prop}

\begin{proof}
We shall show that
\begin{equation} \label{eqn_idelie}
e(1)\left[\col{V}\right] = \clalg{V}.
\end{equation}
The proposition will then follow as a simple consequence of this. Since the Lie monomials generate the Lie subalgebra of primitive elements of the Hopf algebra $(\ctalg{V},\mu,\Delta)$ we have;
\[ \clalg{V} \subset \{ x \in \col{V}:s(x)=2x \} = e(1)\left[\col{V}\right]. \]

To prove the converse we consider the filtration $\{F_p(\ctalg{V})\}_{p=0}^\infty$ of \ctalg{V} defined by Definition \ref{def_pbwfil}. Let $x \in F_p(\ctalg{V})$ and suppose that $s(x)=2x$, then by equation \eqref{eqn_shufev};
\[ 2x = s(x) = 2^p x \mod F_{p-1}(\ctalg{V}) \]
It follows that if $p>1$ then $x \in F_{p-1}(\ctalg{V})$. By induction on $p$ we conclude that $x \in F_1(\ctalg{V})$. Equation \eqref{eqn_idelie} now follows from the fact that $(\Sigma V^*)^{\cotimes n} \subset F_n(\ctalg{V})$ and from the identity
\[ F_1(\ctalg{V}) = \gf \oplus \clalg{V}. \]

Equation \eqref{eqn_idelie} implies that
\[ \tilde{e}(1)\left[\col{V}\right] = \Sigma V^* \cotimes \clalg{V}. \]
We see from diagram \eqref{fig_clsdtf} that we have the following identity:
\begin{equation} \label{eqn_liespl}
l\left[\drof[Lie]{\clalg{V}}\right] = \thiso[Ass]\tilde{e}(1)\left[\col{V}\right].
\end{equation}
It follows that $\thiso[Ass]^{-1} \circ l$ is an isomorphism between \caq{V}{V^*} and \echoch{V}{V^*}{1} and that $l$ is split.
\end{proof}

\begin{rem}
Given a \ci-algebra $V$ with \ci-structure $m$ it is also possible to construct a Hodge decomposition of the Hochschild cohomology of $V$ with coefficients in $V$. We will describe it briefly here. Choose a topological basis $\{t_i\}_{i\in I}$ of $\Sigma V^*$. One introduces the operator $\bar{s}:\choch{V}{V} \to \choch{V}{V}$ as
\[ \bar{s}(\xi):= \sum_{i \in I} \mu\Delta\xi(t_i)\partial_{t_i}, \quad \xi \in \Der(\ctalg{V}); \]
which is obviously defined independently of the choice of basis. It is then possible to show, using equations \ref{eqn_shucom} and \ref{eqn_shucomdummya}, that this operator commutes with the differential $d:=\ad m$. It follows that the differential $d$ preserves the eigenspaces of the operator $\bar{s}$ and these eigenspaces can be identified using Lemma \ref{lem_shidem} in the usual way to give the Hodge decomposition of \choch{V}{V}. In particular it is possible to show that the Harrison complex \caq{V}{V} splits off of the Hochschild complex \choch{V}{V}.
\end{rem}
\section{The Hodge Decomposition of Cyclic Hochschild Cohomology}\label{cyclichodge}

In this section we will build on our results from the last section to construct the Hodge decomposition of the cyclic Hochschild cohomology of a \ci-algebra. Recall that in section \ref{sec_iachom} we defined two quasi-isomorphic complexes \cchoch{V} and \ctsygan{V} which compute the cyclic Hochschild cohomology of an \ai-algebra. We will describe a Hodge decomposition for both of these complexes based on the spectral decomposition of the modified shuffle operator $s$ and show that the quasi-isomorphism defined in Lemma \ref{lem_cyciso} respects this decomposition.

In order to apply the work we carried out in section \ref{sec_hdgdec} to the cyclic complexes, we will need to prove the following lemma:

\begin{lemma} \label{lem_nrmshu}
Let $W$ be a profinite graded module and consider the maps
\[ s,\tilde{s}:\prod_{i=1}^\infty W^{\cotimes i} \to \prod_{i=1}^\infty W^{\cotimes i} \]
defined in Definition \ref{def_shufop}. We have the following identities:
\begin{enumerate}
\item[(i)]
\[ 2\tilde{s} \circ N = N \circ s. \]
\item[(ii)]
\[ s\circ(1-z) = (1-z) \circ \tilde{s}. \]
\end{enumerate}
\end{lemma}

\begin{proof}
\
\begin{enumerate}
\item[(i)]
Let $\thiso[Ass]:\prod_{i=1}^\infty W^{\cotimes i} \to \drof[Ass]{\widehat{T}W}$ be the map defined by Lemma \ref{lem_ofisom} and let
\[ \mu:\de[Ass]{\widehat{T}W}\cotimes\de[Ass]{\widehat{T}W} \to \de[Ass]{\widehat{T}W} \]
be the multiplication map $\mu(x \cotimes y):= x \cdot y$. Let $\alpha_i := (-1)^{(|x_1|+\ldots+|x_{i-1}|)(|x_i|+\ldots+|x_n|)}$ and let $x:=x_1\cotimes\ldots\cotimes x_n \in W^{\cotimes n}$;
\begin{displaymath}
\begin{split}
\tilde{s} N(x) & = \tilde{s}\left(\sum_{i=1}^{n} \alpha_i x_i\cotimes\ldots\cotimes x_n\cotimes x_1 \cotimes \ldots \cotimes x_{i-1}\right), \\
& = \sum_{i=1}^{n} \alpha_i x_i\cotimes \mu\Delta(x_{i+1}\cotimes\ldots\cotimes x_n\cotimes x_1 \cotimes \ldots \cotimes x_{i-1}). \\
N s(x) & = \thiso[Ass]^{-1}d\mu\Delta(x), \\
& = \thiso[Ass]^{-1}d\mu(\Delta(x_1)\ldots\Delta(x_n)), \\
& = \thiso[Ass]^{-1}\mu(d\cotimes 1+1\cotimes d)[(x_1\cotimes 1+1\cotimes x_1)\ldots(x_n\cotimes 1+1\cotimes x_n)]. \\
\end{split}
\end{displaymath}
Now $d\cotimes 1$ and $1\cotimes d$ are derivations on the algebra $\de[Ass]{\widehat{T}W}\cotimes\de[Ass]{\widehat{T}W}$, therefore
\begin{equation} \label{eqn_nrmshudummy}
Ns(x)=  \thiso[Ass]^{-1}\mu\left(\sum_{i=1}^n (-1)^{|x_1|+\ldots+|x_{i-1}|}(x_1\cotimes 1+1\cotimes x_1)\ldots (dx_i\cotimes 1+1\cotimes dx_i)\ldots (x_n\cotimes 1+1\cotimes x_n)\right).
\end{equation}

We will need the following auxiliary calculation in order to complete the proof: Let $u,v,w \in \de[Ass]{\widehat{T}W}$;
\begin{displaymath}
\begin{split}
\mu[(u\cotimes 1+1\cotimes u)\cdot(v\cotimes w)] & = \mu(uv\cotimes w + (-1)^{|u||v|}v\cotimes uw), \\
& = uvw + (-1)^{|u||v|}vuw. \\
\mu[(v\cotimes w)\cdot(u\cotimes 1+1\cotimes u)] & = \mu((-1)^{|u||w|}vu\cotimes w + v\cotimes wu), \\
& = (-1)^{|u||w|}vuw + vwu. \\
\end{split}
\end{displaymath}
We see that
\[ \mu\left([u\cotimes 1 + 1\cotimes u , v\cotimes w]\right) = [u,\mu(v \cotimes w)] = 0 \mod [\de[Ass]{\widehat{T}W},\de[Ass]{\widehat{T}W}]. \]

Let $\beta_i:=(-1)^{(|x_i|+1)(|x_1|+\ldots+|x_{i-1}|+|x_{i+1}|+\ldots+|x_n|)}$. Applying the preceding calculation to equation \eqref{eqn_nrmshudummy} yields;
\begin{displaymath}
\begin{split}
Ns(x) = & \thiso[Ass]^{-1}\mu\left[\sum_{i=1}^n \alpha_i(dx_i\cotimes 1)\Delta(x_{i+1})\ldots \Delta(x_n)\Delta(x_1)\ldots\Delta(x_{i-1})\right] \\
& + \thiso[Ass]^{-1}\mu\left[ \sum_{i=1}^n \alpha_i\beta_i\Delta(x_{i+1})\ldots \Delta(x_n)\Delta(x_1)\ldots\Delta(x_{i-1})(1\cotimes dx_i)\right], \\
= & \thiso[Ass]^{-1}\left[\sum_{i=1}^n \alpha_i dx_i\cdot\mu\Delta(x_{i+1}\cotimes\ldots\cotimes x_n\cotimes x_1\cotimes\ldots\cotimes x_{i-1})\right] \\
& + \thiso[Ass]^{-1}\left[ \sum_{i=1}^n \alpha_i\beta_i\mu\Delta(x_{i+1}\cotimes\ldots\cotimes x_n \cotimes x_1 \cotimes \ldots \cotimes x_{i-1})\cdot dx_i \right], \\
= & 2\sum_{i=1}^n \alpha_i x_i\cotimes\mu\Delta(x_{i+1}\cotimes\ldots\cotimes x_n\cotimes x_1\cotimes\ldots\cotimes x_{i-1}), \\
= & 2\tilde{s}N(x). \\
\end{split}
\end{displaymath}
\item[(ii)]
Let $x \in W$ and $y \in W^{\cotimes n}$ for $n \geq 0$: From equation \eqref{eqn_shucomdummya} we see that
\[ s(1-z)(xy) = s([x,y]) = [x,s(y)] = (1-z)\tilde{s}(xy). \]
\end{enumerate}
\end{proof}

\begin{cor} \label{cor_nrmshu}
Let $W$ be a profinite graded module and consider the maps
\[ e(i),\tilde{e}(i):\prod_{i=1}^\infty W^{\cotimes i} \to \prod_{i=1}^\infty W^{\cotimes i} \]
defined by Definition \ref{def_idemdf}:
\begin{enumerate}
\item[(i)]
For all $i \geq 0$,
\[ \tilde{e}(i) \circ N = N \circ e(i+1). \]
\item[(ii)]
For all $i \geq 1$,
\[ e(i) \circ (1-z) = (1-z) \circ \tilde{e}(i). \]
\end{enumerate}
\end{cor}

\begin{proof}
This is just a formal consequence of Lemma \ref{lem_nrmshu} and Lemma \ref{lem_shidem}. $e(i)$ and $\tilde{e}(i)$ are just the projections onto the eigenspaces
\[ \{ x \in \prod_{i=1}^\infty W^{\cotimes i} : s(x)=\lambda_i x \} \quad \text{and} \quad \{ x \in \prod_{i=1}^\infty W^{\cotimes i} : \tilde{s}(x)=\lambda_i x \} \]
respectively (where $\lambda_i:=2^i$).

Part (ii) follows immediately from part (ii) of Lemma \ref{lem_nrmshu} and part (i) follows from part (i) of Lemma \ref{lem_nrmshu} and the observation $\lambda_{i+1}=2\lambda_i$.
\end{proof}

We are now in a position to construct the Hodge decomposition of cyclic cohomology:

\begin{theorem}
Let $V$ be a \ci-algebra:
\begin{enumerate}
\item[(i)]
The cyclic Hochschild complex \cchoch{V} splits as a direct product of the subcomplexes $\Sigma\left(e(i)\left[\prod_{i=1}^\infty \left((\Sigma V^*)^{\cotimes i}\right)_{Z_i}\right]\right)$:
\[ \cchoch{V}=\prod_{i=1}^\infty\left(\Sigma\left(e(i)\left[\prod_{i=1}^\infty \left((\Sigma V^*)^{\cotimes i}\right)_{Z_i}\right]\right),b'\right). \]
\item[(ii)]
Let $Q_i$ and $\widetilde{Q}_i$ denote the modules
\[ \Sigma\left(e(i)\left[\col{V}\right]\right) \quad \text{and} \quad \Sigma\left(\tilde{e}(i)\left[\col{V}\right]\right) \]
respectively. The Tsygan bicomplex \ctsygan{V} splits as a direct product of subcomplexes,
\[ \ctsygan{V} = \prod_{i=0}^\infty \Gamma_i; \]
where $\Gamma_i$ is the subcomplex,
\[ \Gamma_i:= \xymatrix{(\widetilde{Q}_i,b) \ar^{1-z}[r] & (Q_i,b') \ar^{N}[r] & (\widetilde{Q}_{i-1},b) \ar^{1-z}[r] & (Q_{i-1},b') \ar^{N}[r] & \ldots \ar^{1-z}[r] & (Q_1,b') \ar^{N}[r] & (\widetilde{Q}_0,b)}. \]
\end{enumerate}
\end{theorem}

\begin{proof}
\
\begin{enumerate}
\item[(i)]
Part (ii) of Corollary \ref{cor_nrmshu} tells us that the maps $e(i):\col{V} \to \col{V}$ could be lifted uniquely to maps on the module of coinvariants;
\[ e(i):\prod_{i=1}^\infty \left((\Sigma V^*)^{\cotimes i}\right)_{Z_i} \to \prod_{i=1}^\infty \left((\Sigma V^*)^{\cotimes i}\right)_{Z_i}, \quad i \geq 1. \]
By Lemma \ref{lem_cyciso} we know that
\[ \cchoch{V} = \left(\Sigma\left(\prod_{i=1}^\infty \left((\Sigma V^*)^{\cotimes i}\right)_{Z_i}\right),b'\right). \]
By Lemma \ref{lem_shidem} and Lemma \ref{lem_shucom} we see that this complex splits as a direct product of subcomplexes:
\[ \left(\Sigma\left(\prod_{i=1}^\infty \left((\Sigma V^*)^{\cotimes i}\right)_{Z_i}\right),b'\right)=\prod_{i=1}^\infty\left(\Sigma\left(e(i)\left[\prod_{i=1}^\infty \left((\Sigma V^*)^{\cotimes i}\right)_{Z_i}\right]\right),b'\right). \]
\item[(ii)]
Lemma \ref{lem_shidem} tells us that the bigraded module \ctsygan{V} splits as a direct product of bigraded submodules
\[ \ctsygan{V}=\prod_{i=0}^\infty \Gamma_i. \]
Corollary \ref{cor_nrmshu} tells us that the operators $(1-z)$ and $N$ restrict to maps
\begin{align*}
& N: e(i+1)\left[\col{V}\right] \to \tilde{e}(i)\left[\col{V}\right], \\
& (1-z): \tilde{e}(i)\left[\col{V}\right] \to e(i)\left[\col{V}\right] .
\end{align*}
Combining this with Lemma \ref{lem_shucom} we conclude that each $\Gamma_i$ is actually a subcomplex of \ctsygan{V} whence the result.
\end{enumerate}
\end{proof}

\begin{rem}
Note that since the bicomplex \ctsygan{V} splits as the direct product of the $\Gamma_i$'s, where the latter are
bicomplexes located within vertical strips of \emph{finite width}, it follows that both spectral sequences associated with \ctsygan{V} converge to its cohomology. In particular, it implies that \hchoch{V} is homotopy invariant.
\end{rem}

\begin{rem}
Recall that in Remark \ref{rem_connes} we asserted the existence of a normalised Connes complex computing the cyclic cohomology of a unital \ai-algebra $V$. It is possible to construct a Hodge decomposition of this complex using Corollary \ref{cor_nrmshu} to obtain an extension of Theorem 4.6.7 of \cite{loday} for unital \ai-algebras.
\end{rem}

\begin{rem}
Let $q:\cchoch{V} \overset{\sim}{\to} \ctsygan{V}$ be the quasi-isomorphism defined by Lemma \ref{lem_cyciso} part (ii). Corollary \ref{cor_nrmshu} part (i) implies that $q$ respects the Hodge decomposition of these complexes, that is to say that for all $i \geq 0$ it restricts to a map,
\[ q: \left(\Sigma\left(e(i+1)\left[\prod_{i=1}^\infty \left((\Sigma V^*)^{\cotimes i}\right)_{Z_i}\right]\right),b'\right) \overset{\sim}{\to} \Gamma_i. \]
In particular, the Hodge decomposition of \cchoch{V} and the Hodge decomposition of \ctsygan{V} agree on the level of cohomology.
\end{rem}

\begin{defi}
Let $V$ be a \ci-algebra. We will define the complex \ecchoch{V}{i} as the subcomplex of \cchoch{V} given by the formula;
\[ \ecchoch{V}{i}:=\left(\Sigma\left(e(i+1)\left[\prod_{i=1}^\infty \left((\Sigma V^*)^{\cotimes i}\right)_{Z_i}\right]\right),b'\right), \quad i\geq 0 \]
and denote its cohomology by \ehchoch{V}{i}.
\end{defi}

We will now describe how the Hodge decomposition splits the long exact sequence of Proposition \ref{prop_sbiseq}:

\begin{prop} \label{prop_sbidec}
Let $V$ be a unital \ci-algebra. The Hodge decomposition of the Hochschild and cyclic Hochschild cohomologies respects the long exact sequence of Proposition \ref{prop_sbiseq}, that is to say that for all $i\geq 0$ we have the following long exact sequence in cohomology:
\[  \ldots\xymatrix{ \ehchoch[n-2]{V}{i} \ar^{S}[r] & \ehchoch[n]{V}{i+1} \ar^{I}[r] & \ehhoch[n]{V}{V^*}{i+1} \ar^{B}[r] & \ehchoch[n-1]{V}{i} \\ }\ldots. \]
\end{prop}

\begin{proof}
It is a simple check using the definitions of the maps $I$ and $S$ (cf. Remark \ref{rem_sbiseq}) and the Hodge decompositions to see that $I$ and $S$ respect the Hodge decomposition as claimed. It only remains to prove that the map $B$ restricts to a map
\[ B: \ehhoch[n]{V}{V^*}{i+1} \to \ehchoch[n-1]{V}{i}. \]

By Proposition \ref{prop_htpret}, every cocycle in \choch{V}{V^*} is cohomologous to a normalised cocycle. Let $x$ be such a normalised cocycle, i.e. $x \in \Sigma V^* \cotimes \widehat{T}(\Sigma V/\gf)^*$ and $b(x)=0$: It follows from figure \eqref{eqn_sbiseq} that
\[ B(x) = -h(1-z)[x] = -h(x). \]
Clearly $\tilde{e}(i+1)[x]$ is still normalised so;
\begin{displaymath}
\begin{split}
B(\tilde{e}(i+1)[x]) & = -h(\tilde{e}(i+1)[x]), \\
& = -e(i+1)h(x); \\
\end{split}
\end{displaymath}
hence the image of $B$ lies in \ehchoch[n-1]{V}{i} as claimed.
\end{proof}

Let $V$ be a \ci-algebra and consider the map $l:\drzf[Lie]{\clalg{V}} \to \drzf[Ass]{\ctalg{V}}$ defined by equation \eqref{eqn_liasmp}. By Lemma \ref{lem_plcomm} this is a map of complexes,
\[ l: \ccaq{V} \to \cchoch{V}. \]
We shall now show that this map is in fact split:

\begin{prop} \label{prop_caqspt}
Let $V$ be a \ci-algebra. There is an isomorphism
\[ \ccaq{V} \cong \ecchoch{V}{1}. \]
In particular the map $l: \ccaq{V} \to \cchoch{V}$ is split.
\end{prop}

\begin{proof}
We shall establish the following identity:
\begin{equation} \label{eqn_caqsptdummy}
l\left[\drzf[Lie]{\clalg{V}}\right] = e(2)\left[\prod_{i=1}^\infty \left((\Sigma V^*)^{\cotimes i}\right)_{Z_i}\right].
\end{equation}
The proposition will then follow as a result.

Consider the following commutative diagram:
\begin{displaymath}
\xymatrix{
& 0 & 0 & 0 & \\
0 \ar[r] & \frac{\prod_{i=1}^\infty \left((\Sigma V^*)^{\cotimes i}\right)_{Z_i}}{e(2)\left[\prod_{i=1}^\infty \left((\Sigma V^*)^{\cotimes i}\right)_{Z_i}\right]} \ar[u] \ar^{N}[r] & \frac{\col{V}}{\tilde{e}(1)\left[\col{V}\right]} \ar[u] \ar[r] & \frac{\col{V}}{N \cdot \col{V} + \tilde{e}(1)\left[\col{V}\right]} \ar[u] \ar[r] & 0 \\
0 \ar[r] & \prod_{i=1}^\infty \left((\Sigma V^*)^{\cotimes i}\right)_{Z_i} \ar^{d}[r] \ar[u] & \drof[Ass]{\ctalg{V}} \ar^{\overline{\thiso[Ass]^{-1}}}[u] \ar[r] & \frac{\drof[Ass]{\ctalg{V}}}{d\left(\drzf[Ass]{\ctalg{V}}\right)} \ar^{\overline{\thiso[Ass]^{-1}}}[u] \ar[r] & 0 \\
0 \ar[r] & \drzf[Lie]{\clalg{V}} \ar^{l}[u] \ar^{d}[r] & \drof[Lie]{\clalg{V}} \ar^{l}[u] \ar[r] & \frac{\drof[Lie]{\clalg{V}}}{d\left(\drzf[Lie]{\clalg{V}}\right)} \ar^{l}[u] \ar[r] & 0 \\
& 0 \ar[u] & 0 \ar[u] & 0 \ar[u] & \\
}
\end{displaymath}
where $\overline{\thiso[Ass]^{-1}}$ denotes the map induced by $\thiso[Ass]^{-1}$ by composing it with the relevant projection. It follows from Lemma \ref{lem_shidem} and Corollary \ref{cor_nrmshu} that the top row is exact and it follows from Lemma \ref{lem_poinca} that the middle and bottom rows are also exact. Lemma \ref{lem_nrmdif} and Lemma \ref{lem_clsdtf} imply that the right column is exact whilst equation \eqref{eqn_liespl} implies that the middle column is also exact. It follows from the $3\times 3$-Lemma that the left column is exact which in turn implies \eqref{eqn_caqsptdummy} and this proposition.
\end{proof}

Let $V$ be a \ci-algebra. By Remark \ref{rem_infdef} it has the structure of a complex given by a differential
\[ \check{m}_1:V \to V \]
and of course the dual $V^*$ has the structure of a complex given by the dual of the map $\check{m}_1$. We will denote the cohomology of this complex by $H^\bullet(V^*)$. We will now use the results of this section to prove the following important result:

\begin{lemma} \label{lem_aqlesc}
Let $V$ be a \ci-algebra. We have the following long exact sequence in cohomology:
\[ \xymatrix{ \ldots \ar^-{B}[r] & H^{n-2}(V^*) \ar^{S}[r] & \hcaq[n]{V} \ar^-{I}[r] & \haq[n]{V}{V^*} \ar^-{B}[r] & H^{n-1}(V^*) \ar^-{S}[r] & \ldots \\ }. \]
\end{lemma}

\begin{proof}
This follows from setting $i=0$ in the long exact sequence of Proposition \ref{prop_sbidec}. Proposition \ref{prop_aqsplt} allows us to identify \ehhoch{V}{V^*}{1} with \haq{V}{V^*} whilst Proposition \ref{prop_caqspt} allows us to identify \ehchoch{V}{1} with \hcaq{V}.

It follows from Lemma \ref{lem_nrmdif} part (i) and equation  \eqref{eqn_idelie} that
\[ \ecchoch{V}{0} = \left(\Sigma\left(e(1)\left[\prod_{i=1}^\infty \left((\Sigma V^*)^{\cotimes i}\right)_{Z_i}\right]\right),b'\right) = \left( V^* , m_1 \right); \]
where $m_1$ is the linear part of the \ci-structure on $V$. The result now follows.
\end{proof}

\begin{rem} \label{rem_ordgrd}
Suppose now that $V$ is a strictly graded commutative algebra (in which case it is a \ci-algebra). In this case there is a bigrading on \caq{V}{V^*} and \ccaq{V}. We say a 0-form $\alpha \in \ccaq{V}$ has bidegree $(i,j)$ if it is a 0-form of order $i$ and has degree $j$ as an element in the profinite graded module $\Sigma\drzf[Lie]{\clalg{V}}$. Similarly we say a 1-form $\alpha \in \caq{V}{V^*}$ has bidegree $(i,j)$ if it is a 1-form of order $i$ and has degree $j$ as an element in the profinite graded module $\Sigma\drof[Lie]{\clalg{V}}$. The differentials on \ccaq{V} and \caq{V}{V^*} both have bidegree $(1,1)$. In this situation we can formulate and prove the following corollary:
\end{rem}

\begin{cor} \label{cor_caqiso}
Let $V$ be a unital strictly graded commutative algebra. The map
\[ I:\hcaq[i+1,j]{V} \to \haq[ij]{V}{V^*} \]
of Lemma \ref{lem_aqlesc} is;
\begin{enumerate}
\item[(i)]
a monomorphism if $i=1$,
\item[(ii)]
an epimorphism if $i=2$,
\item[(iii)]
an isomorphism if $i \geq 3$.
\end{enumerate}
\end{cor}

\begin{proof}
A straightforward check utilising the definitions shows that the long exact sequence of Proposition \ref{prop_sbiseq} respects the bigrading on the cohomology referred to above, that is to say we have a long exact sequence in cohomology
\[ \xymatrix{ \ldots \ar^-{B}[r] & \hchoch[i-1,j-2]{V} \ar^{S}[r] & \hchoch[i+1,j]{V} \ar^{I}[r] & \hhoch[ij]{V}{V^*} \ar^{B}[r] & \hchoch[i,j-1]{V} \ar^-{S}[r] & \ldots \\ }. \]

Since the long exact sequence of Lemma \ref{lem_aqlesc} is derived from this long exact sequence we obtain a bigraded long exact sequence
\[ \xymatrix{ \ldots \ar^-{B}[r] & H^{i-1,j-2}(V^*) \ar^{S}[r] & \hcaq[i+1,j]{V} \ar^-{I}[r] & \haq[ij]{V}{V^*} \ar^-{B}[r] & H^{i,j-1}(V^*) \ar^-{S}[r] & \ldots \\ }. \]
Since $H^{\bullet\bullet}(V^*)$ is obviously concentrated in bidegree $(1,\bullet)$ the map
\[ B: \haq[ij]{V}{V^*} \to H^{i,j-1}(V^*) \]
is zero for all $j \in \mathbb{Z}$ and $i \neq 1$ and the map
\[ S: H^{i-1,j-2}(V^*) \to \hcaq[i+1,j]{V} \]
is zero for all $j \in \mathbb{Z}$ and $i \neq 2$, whence the result.
\end{proof}

\section{Finite level structures and Obstruction theory} \label{sec_obstrc}

In this section we define $C_n$-algebras and discuss the problem of lifting a $C_n$-structure to a $C_{n+1}$-structure. The obstruction theory of $n$-algebras is in principle known to the experts, but it is hard to find an explicit reference in the literature, especially in the $C_n$-algebra case. The precursor for this type of obstruction theory is the seminal work of A. Robinson, \cite{robins}, although the main ideas go back to early work on algebraic deformation theory, \cite{gerste}, \cite{nirich}. In the context of unstable homotopy theory, somewhat similar constructions were employed by Stasheff in \cite{stalnm}. We shall only consider the case of \emph{minimal} $C_n$-structures as in this case it is possible to give a convenient interpretation of obstructions in cohomological terms.

Although it is possible to discuss the obstruction theory of $L_n$ and $A_n$-algebras we shall only describe the obstruction theory of $C_n$-algebras. This is in order to avoid repeating ourselves and since in our further applications of obstruction theory we will work exclusively with \ci-algebras. It should however be clear how to formulate and prove the $L_n$ and $A_n$ analogues of the results contained in this section with most of the definitions and theorems carrying over almost verbatim.

Let $V$ be a free graded module. Recall from equation \eqref{eqn_vecfrm} that any vector field $m\in \Der(\clalg{V})$ could be written in the form
\[ m = m_1 + m_2 + \ldots + m_n + \ldots \]
where $m_i$ is a vector field of order $i$. Since we often need to work modulo the endomorphisms of \clalg{V} of order $\geq n$ we will denote the module of such endomorphisms by $(n)$.

\subsection{Obstruction theory for $C_n$-algebra structures} \label{sec_obsalg}

In this subsection we will develop the obstruction theory for $C_n$-algebra structures. Let us begin with a definition:

\begin{defi} \label{def_minstr}
Let $V$ be a free graded module. A \emph{minimal} $C_n$-structure ($n \geq 3$) on $V$ is a vector field $m \in \Der(\clalg{V})$ of degree one which has the form
\[ m = m_2 + \ldots + m_{n-1} \quad , \quad m_i \text{ has order $i$} \]
and satisfies the condition $m^2 = 0 \mod (n+1)$.\footnote{Note that this definition entails a slight abuse of notation since the $m_i$'s are in fact \emph{dual} to the maps defined by diagram \eqref{fig_infdef}}

Let $m$ and $m'$ be two minimal $C_n$-structures on $V$. We say $m$ and $m'$ are equivalent if there is a diffeomorphism $\phi \in \Aut(\clalg{V})$ of the form
\begin{equation} \label{eqn_diffrm}
\phi = \id + \phi_2 + \phi_3 + \ldots + \phi_k +\ldots
\end{equation}
where $\phi_i$ is an endomorphism of order $i$, such that
\[ \phi \circ m \circ \phi^{-1} = m' \mod (n). \]
\end{defi}

\begin{rem}
Recall that we refer to diffeomorphisms of the form \eqref{eqn_diffrm} as \emph{pointed diffeomorphisms}.
\end{rem}

\begin{rem}
Obviously in the definition of minimal $L_n$ and $A_n$ structures one should replace \clalg{V} with \csalg{V} and \ctalg{V} respectively. The unital minimal $A_n$ and $C_n$ structures are defined precisely as in the $\infty$-case (cf. Definition \ref{def_infstr}).
\end{rem}

\begin{rem}
Note that this definition is slightly at odds with the definition of $n$-algebras (specifically $A_n$-algebras) given by Stasheff in \cite{staha1} and \cite{staha2}. Every \emph{minimal} $C_n$-algebra is a $C_n$-algebra under Stasheff's definition, however, given a \emph{minimal} $C_n$-structure $m=m_2+\ldots+m_{n-1}$ and an arbitrary vector field $m_n$ of order $n$, $m':=m+m_n$ is a $C_n$-algebra under Stasheff's definition which is obviously not a \emph{minimal} $C_n$-algebra. This distinction will be necessary in order to develop an obstruction theory.
\end{rem}

\begin{rem} \label{rem_ordgde}
Clearly if $m=m_2+\ldots+m_{n-1}$ is a minimal $C_n$-algebra then $m_2$ determines the structure of a strictly graded commutative algebra on the underlying free graded module $V$ (cf. Remark \ref{rem_infdef}) which we will call the underlying commutative algebra. Observe that two equivalent minimal $C_n$-structures have the same underlying commutative algebra.

Since the underlying algebra $A:=(V,m_2)$ is a strictly graded commutative algebra, \caq{A}{A} can be given a bigrading. We say a vector field $\xi \in \caq{A}{A}$ has bidegree $(i,j)$ if it is a vector field of order $i$ and has degree $j$ as an element in the profinite graded module $\Sigma^{-1}\Der(\clalg{V})$. The differential on \caq{A}{A} then has bidegree $(1,1)$.
\end{rem}

It will be useful to introduce the following definition:

\begin{defi}
Let $A:=(V,\mu_2)$ be a strictly graded commutative algebra, then the moduli space of minimal $C_n$-structures on $V$ fixing $\mu_2$ is denoted by \mstr{C}{n}{A} and defined as the quotient of the set
\[ \{ m:\clalg{V} \to \clalg{V} : m \text{ is a minimal } C_n\text{-structure and } m_2=\mu_2 \} \]
by the equivalence relation defined in Definition \ref{def_minstr}.
\end{defi}

We will now describe the appropriate terminology which is necessary in discussing extensions of $C_n$-structures to structures of higher order:

\begin{defi}
Let $V$ be a free graded module and let $m = m_2 + \ldots + m_{n-1}$ be a minimal $C_n$-structure on $V$. We say that $m$ is an extendable $C_n$-structure if there exists a vector field $m_n \in \Der(\clalg{V})$ of order $n$ and degree one such that $m + m_n$ is a $C_{n+1}$-structure on $V$ and we call $m_n$ an extension of $m$.

Let $m_n$ and $m'_n$ be two extensions of $m$. We say that $m_n$ and $m'_n$ are equivalent if there exists a diffeomorphism $\phi \in \Aut(\clalg{V})$ of the form
\[ \phi = \id + \phi_{n-1} + \phi_n + \ldots + \phi_{n+k} + \ldots \]
where $\phi_i$ is an endomorphism of order $i$, such that
\begin{equation} \label{eqn_exteqv}
\phi \circ (m + m_n) \circ \phi^{-1} = m + m'_n \mod (n+1).
\end{equation}
The quotient of the set of all extensions of $m$ by this equivalence relation will be denoted by \mext{n}{m}.
\end{defi}

\begin{rem} \label{rem_exteqv}
Let $m_n$ and $m'_n$ be two equivalent extensions of $m$, then there is a diffeomorphism $\phi$ satisfying \eqref{eqn_exteqv} and a vector field $\xi$ of order $n-1$ and degree zero such that $\phi$ has the form
\[ \phi = \id + \xi + \text{ endomorphisms of order}\geq n. \]
From this we conclude that $m'_n = m_n + [\xi,m_2]$.

Conversely if there exists a vector field $\xi$ of order $n-1$ and degree zero such that
\[ m'_n = m_n + [\xi,m_2] \]
then the diffeomorphism $\phi:=\exp(\xi)$ satisfies \eqref{eqn_exteqv}. This means that two extensions of $m$ are equivalent if and only if their difference is a Harrison coboundary.
\end{rem}

\begin{defi} \label{def_obstrc}
Let $V$ be a free graded module and let $m = m_2+\ldots + m_{n-1}$ be a $C_n$-structure on $V$. We define the vector field $\Obs(m)$ of order $n+1$ and degree 2 by
\begin{equation} \label{eqn_obstrc}
\Obs(m) := \frac{1}{2}\sum_{\begin{subarray}{c} i+j = n+2 \\ 3\leq i,j \leq n-1 \end{subarray}} [m_i, m_j].
\end{equation}
\end{defi}

We now formulate and prove a result analogous to theorems 7.14 and 7.15 of \cite{hamilt}:

\begin{theorem} \label{thm_obstrc}
Let $A:=(V,m_2)$ be a strictly graded commutative algebra. For all $n \geq 3$ \eqref{eqn_obstrc} induces a map
\begin{displaymath}
\begin{array}{ccc}
\mstr{C}{n}{A} & \to & \haq[n+1,3]{A}{A}, \\
m & \mapsto & \Obs(m); \\
\end{array}
\end{displaymath}
which we will denote by $\Obs_n$. The kernel of this map consists of precisely those $C_n$-structures which are extendable:
\[ \{ m \in \mstr{C}{n}{A} : m \textnormal{ is extendable} \} = \ker(\Obs_n). \]
\end{theorem}

\begin{proof}
First of all let us show that given any $C_n$-structure $m=m_2+\ldots+m_{n-1}$ on $V$, $\Obs(m)$ is a cocycle in \caq{A}{A}. Since $m$ is a $C_n$-structure we have the following equality for all $k\leq n-1$:
\[ [m_2,m_k] = -\frac{1}{2}\sum_{\begin{subarray}{c} i+j = k+2 \\ 3\leq i,j \leq k-1 \end{subarray}} [m_i,m_j]. \]
Using this we show that $\Obs(m)$ is a Harrison cocycle as follows:
\begin{displaymath}
\begin{split}
[m_2,\Obs(m)] & = \frac{1}{2}\sum_{\begin{subarray}{c} i+j = n+2 \\ 3\leq i,j \leq n-1 \end{subarray}} [[m_2,m_i],m_j] - \frac{1}{2}\sum_{\begin{subarray}{c} i+j = n+2 \\ 3\leq i,j \leq n-1 \end{subarray}} [m_i,[m_2,m_j]], \\
& = \sum_{\begin{subarray}{c} i+j = n+2 \\ 3\leq i,j \leq n-1 \end{subarray}} [[m_2,m_i],m_j], \\
& = \sum_{\begin{subarray}{c} i+j+k = n+4 \\ 3\leq i,j,k \leq n-2 \end{subarray}} [[m_i,m_j],m_k], \\
& = 0
\end{split}
\end{displaymath}
where the last equality follows from the Jacobi identity.

Next we need to show that if 
\[ m=m_2+m_3+\ldots+m_{n-1} \quad \text{and} \quad m'=m_2+m'_3+\ldots+m'_{n-1} \]
are two equivalent $C_n$-structures on $V$ then $\Obs(m)$ and $\Obs(m')$ are cohomologous Harrison cocycles.  Since $m$ and $m'$ are equivalent there exists a \emph{pointed} diffeomorphism $\phi \in \Aut(\clalg{V})$ such that
\[ m'=\phi\circ m \circ\phi^{-1} \mod(n). \]
This means there exists a vector field $\xi$ of order $n$ and degree 1 such that
\[ m' = \phi\circ m \circ\phi^{-1} + \xi \mod(n+1). \]

By the definition of $\Obs(m)$ (cf. equation \eqref{eqn_obstrc}) we have the equality
\[ m^2 = \Obs(m) \mod (n+2) \]
and likewise we have a similar equality for $m'$. Now from the calculation
\begin{displaymath}
\begin{split}
m'^2 & = \phi \circ m^2\circ \phi^{-1} + [m_2,\xi] \mod (n+2), \\
& = \Obs(m) + [m_2,\xi] \mod(n+2) \\
\end{split}
\end{displaymath}
we conclude that $\Obs(m') = \Obs(m) + [m_2,\xi]$.

So far we have proven that map $\Obs_n$ is well defined. In order to finish the proof we need to show that a $C_n$-structure $m=m_2+\ldots+m_{n-1}$ on $V$ is extendable if and only if $\Obs(m)$ is cohomologous to zero. The reason that this is true is because for any vector field $m_n$ of order $n$ and degree 1 the following identity holds:
\begin{equation} \label{eqn_obstrcdummy}
(m+m_n)^2 = \Obs(m) + [m_2,m_n] \mod(n+2).
\end{equation}
\end{proof}

\begin{rem}
Obviously when dealing with $L_n$ and $A_n$-algebras Harrison cohomology should be replaced with Chevalley-Eilenberg and Hochschild cohomology respectively.
\end{rem}

The next result, also reminiscent of deformation theory, analyses different extensions of a $C_n$-structure in terms of Harrison cohomology:

\begin{theorem} \label{thm_extend}
Let $A:=(V,m_2)$ be a strictly graded commutative algebra and let $m \in \mstr{C}{n}{A}$ be an extendable $C_n$-structure, then \haq[n,2]{A}{A} acts freely and transitively on \mext{n}{m}:
\begin{displaymath}
\begin{array}{ccc}
\haq[n,2]{A}{A} \times \mext{n}{m} & \to & \mext{n}{m}, \\
(\xi_n , m_n) & \mapsto & m_n + \xi_n. \\
\end{array}
\end{displaymath}
\end{theorem}

\begin{proof}
By equation \eqref{eqn_obstrcdummy} we see that $m_n$ is an extension of $m$ if and only if
\begin{equation} \label{eqn_extenddummy}
[m_2,m_n] = -\Obs(m).
\end{equation}
This means that if $m_n$ is an extension of $m$ and $\xi_n$ is a Harrison cocycle then $m_n + \xi_n$ is an extension of $m$. Furthermore if $m_n$ and $m'_n$ are equivalent extensions and $\xi_n$ and $\xi'_n$ are cohomologous cocycles then by Remark \ref{rem_exteqv}, $m_n + \xi_n$ and $m'_n + \xi'_n$ are equivalent extensions of $m$, therefore the above action is well defined.

Condition \eqref{eqn_extenddummy} shows us that if $m_n$ and $m'_n$ are two extensions of $m$ then $m_n - m'_n$ is a Harrison cocycle and hence the above action is transitive. Furthermore if $m_n$ is an extension of $m$ and $\xi_n$ is a Harrison cocycle such that $m_n$ and $m_n + \xi_n$ are equivalent extensions then by Remark \ref{rem_exteqv}, $\xi_n$ is a Harrison coboundary, thus the above action is free.
\end{proof}

\subsection{Obstruction theory for $C_n$-algebra morphisms} \label{sec_obsmor}

In this subsection we will develop the obstruction theory for morphisms between two $C_n$-algebras.

\begin{defi} \label{def_minmor}
Let $V$ be a free graded module and let $m$ and $m'$ be two \emph{minimal} $C_n$-structures on $V$. A \emph{minimal} $C_n$-morphism from $m$ to $m'$ is a diffeomorphism $\phi \in \Aut(\clalg{V})$ of degree zero such that
\[ \phi \circ m  = m' \circ \phi \mod (n). \]

Let $\phi$ and $\phi'$ be two such minimal $C_n$-morphisms. We say $\phi$ and $\phi'$ are homotopic if there exists a vector field $\eta$ of degree $-1$ such that
\[ \phi = \phi'\circ\exp([m,\eta]) \mod (n-1). \]
\end{defi}

\begin{rem}
There is a Lie power series $p(x,y):= x + y + \frac{1}{2}[x,y]+\ldots$ belonging to the pro-free Lie algebra on two generators $x$ and $y$ such that for two vector fields $\gamma$ and $\xi$ of degree zero
\[ \exp(p(x,y)) = \exp(x)\circ\exp(y). \]
It follows that homotopy is an equivalence relation.
\end{rem}

We will now introduce the moduli space of $C_n$-morphisms:

\begin{defi}
Let $V$ be a free graded module and let $m$ and $m'$ be two minimal $C_n$-structures on $V$ with the same underlying algebra, that is to say that $m_2 = m'_2$. The moduli space of minimal $C_n$-morphisms from $m$ to $m'$ is denoted by \mmor{n}{m}{m'} and defined as the quotient of the set
\[ \{ \phi:\clalg{V} \to \clalg{V} : \phi \text{ is a \emph{pointed} minimal $C_n$-morphism from $m$ to $m'$} \} \]
by the homotopy equivalence relation defined in Definition \ref{def_minmor}.
\end{defi}

Let $m = m_2 + \ldots + m_n$ be a $C_{n+1}$-structure on a free graded module $V$. We define the corresponding $C_n$-structure $\bar{m}$ on $V$ as
\[ \bar{m}:= m_2 + \ldots + m_{n-1}. \]
Now we introduce the terminology dealing with extensions of $C_n$-morphisms:

\begin{defi} \label{def_morext}
Let $V$ be a free graded module and let $m$ and $m'$ be two minimal $C_{n+1}$-structures on $V$ with the same underlying algebra. We say that a pointed $C_n$-morphism $\phi$ from $\bar{m}$ to $\bar{m}'$ is extendable if there exists a vector field $\gamma \in \Der(\clalg{V})$ of order $n-1$ and degree zero such that $\exp(\gamma)\circ\phi$ is a $C_{n+1}$-morphism from $m$ to $m'$ and we call $\gamma$ an extension of $\phi$.

Let $\gamma$ and $\gamma'$ be two extensions of $\phi$. We say that $\gamma$ and $\gamma'$ are equivalent if there exists a vector field $\eta$ of order $n-2$ and degree $-1$ such that 
\begin{equation} \label{eqn_morextcon}
\exp(\gamma)\circ\phi = \exp(\gamma')\circ\phi\circ\exp([m,\eta]) \mod(n).
\end{equation}
The quotient of the set of all extensions of $\phi$ by this equivalence relation will be denoted by \mmext{n}{\phi}{m}{m'}.
\end{defi}

\begin{rem} \label{rem_morext}
Let $\gamma$ and $\gamma'$ be two extensions of $\phi$ and let $\eta$ be a vector field of order $n-2$ and degree $-1$, then equation \eqref{eqn_morextcon} is satisfied if and only if
\[ \gamma = \gamma' + [m_2,\eta]. \]
It follows that two extensions are equivalent if and only if their difference is a Harrison coboundary in the Harrison complex of the underlying algebra $A:=(V,m_2)$.
\end{rem}

Next we will define the appropriate obstruction to the extension of $C_n$-morphisms.

\begin{defi} \label{def_morobs}
Let $V$ be a free graded module and let $m$ and $m'$ be two minimal $C_{n+1}$-structures on $V$ with the same underlying algebra. Let $\phi$ be a pointed $C_n$-morphism from $\bar{m}$ to $\bar{m}'$. We define the vector field $\obs(\phi)$ of order $n$ and degree $1$ by
\begin{equation} \label{eqn_morobs}
\obs(\phi):= \phi \circ m \circ \phi^{-1} - m' \mod (n+1).
\end{equation}
\end{defi}

We now have all the terminology in place to formulate the analogues of theorems \ref{thm_obstrc} and \ref{thm_extend}:

\begin{theorem} \label{thm_morobs}
Let $A:=(V,m_2)$ be a strictly graded commutative algebra and let $m$ and $m'$ be two minimal $C_{n+1}$-structures ($n \geq 3$) on $V$ whose underlying algebra is $A$, then \eqref{eqn_morobs} induces a map
\begin{displaymath}
\begin{array}{ccc}
\mmor{n}{\bar{m}}{\bar{m}'} & \to & \haq[n,2]{A}{A}, \\
\phi & \mapsto & \obs(\phi); \\
\end{array}
\end{displaymath}
which we will denote by $\obs_n$. The kernel of this map consists of precisely those $C_n$-morphisms which are extendable:
\[ \{ \phi \in \mmor{n}{\bar{m}}{\bar{m}'} : \phi \textnormal{ is extendable} \} = \ker(\obs_n). \]
\end{theorem}

\begin{proof}
Let $\phi$ be a pointed $C_n$-morphism from $\bar{m}$ to $\bar{m}'$, then by Theorem \ref{thm_extend} $\obs(\phi)$ is a cocycle in \caq{A}{A} and $\phi$ is extendable if and only if $\obs(\phi)$ is cohomologous to zero. We need only show that if $\phi'$ is another pointed $C_n$-morphism from $\bar{m}$ to $\bar{m}'$ which is homotopic to $\phi$ then $\obs(\phi)$ and $\obs(\phi')$ are cohomologous cocycles.

Since $\phi$ and $\phi'$ are homotopic there exists a vector field $\gamma$ of order $n-1$ and degree 0 and a vector field $\eta$ of degree $-1$ such that
\[ \phi = \exp(\gamma)\circ\phi'\circ\exp([m,\eta]) \mod (n). \]
Note that for any vector field $\eta$ of degree $-1$;
\[ \exp([m,\eta]) \circ m \circ \exp(-[m,\eta]) = \exp\left(\ad([m,\eta])\right)[m] = m \mod(n+2). \]
We use these facts to demonstrate that the relevant obstructions are cohomologous:
\begin{displaymath}
\begin{split}
\obs(\phi) & = \phi\circ m \circ\phi^{-1}-m' \mod(n+1), \\
& = \exp(\gamma)\circ\phi'\circ\exp([m,\eta]) \circ m \circ \exp(-[m,\eta])\circ\phi'^{-1}\circ\exp(-\gamma) - m' \mod(n+1), \\
& = \phi'\circ m \circ\phi'^{-1} + [\gamma,m_2] - m' \mod(n+1), \\
& = \obs(\phi') + [\gamma,m_2] \mod(n+1). \\
\end{split}
\end{displaymath}
\end{proof}

\begin{theorem} \label{thm_morext}
Let $A:=(V,m_2)$ be a strictly graded commutative algebra and let $m$ and $m'$ be two minimal $C_{n+1}$-structures on $V$ whose underlying algebra is $A$. Let $\phi \in \mmor{n}{\bar{m}}{\bar{m}'}$ be an extendable $C_n$-morphism, then \haq[n-1,1]{A}{A} acts freely and transitively on \mmext{n}{\phi}{m}{m'}:
\begin{displaymath}
\begin{array}{ccc}
\haq[n-1,1]{A}{A} \times \mmext{n}{\phi}{m}{m'} & \to & \mmext{n}{\phi}{m}{m'}, \\
(\xi,\gamma) & \mapsto & \gamma + \xi. \\
\end{array}
\end{displaymath}
\end{theorem}

\begin{proof}
A vector field $\gamma$ of order $n-1$ and degree zero is an extension of $\phi$ if and only if
\begin{equation} \label{eqn_morextdummya}
[m_2,\gamma] = \obs(\phi).
\end{equation}
It follows that if $\gamma$ is an extension of $\phi$ and $\xi$ is a Harrison cocycle then $\gamma + \xi$ is an extension of $\phi$. Moreover if $\gamma$ and $\gamma'$ are equivalent extensions and $\xi$ and $\xi'$ are cohomologous cocycles, then by Remark \ref{rem_morext} $\gamma+\xi$ and $\gamma'+\xi'$ are equivalent extensions, hence the above action is well defined.

If $\gamma$ and $\gamma'$ are two extensions of $\phi$ then by \eqref{eqn_morextdummya} their difference is a cocycle and therefore the above action is transitive. If $\gamma$ is an extension of $\phi$ and $\xi$ is a cocycle such that $\gamma$ and $\gamma+\xi$ are equivalent extensions, then by Remark \ref{rem_morext} $\xi$ is a coboundary and hence the above action is free.
\end{proof}

\section{Formal Noncommutative Symplectic Geometry and Symplectic $\infty$-Algebras}\label{symplectic}

In this section we will review the formal noncommutative \emph{symplectic} geometry introduced by Kontsevich in \cite{kontsg} and \cite{kontfd}. We will consider three types of this symplectic geometry corresponding to three types of algebra; commutative, associative and Lie algebras. We will define the notion of a (non)commutative symplectic form and formulate the Darboux theorem which says that there is a canonical form for a symplectic form on a `flat manifold'.

We will then recall the definition of an $\infty$-algebra with an invariant inner product. The main purpose of this section will be to apply Kontsevich's noncommutative symplectic geometry to obtain an equivalent formulation of these structures in geometrical terms: The main theorem of this section will be to show that such algebras could be equivalently defined as a \emph{symplectic} homological vector field on a formal flat \emph{symplectic} supermanifold. We will also formulate and prove some results for the cohomology of these $\infty$-algebras.

We will begin by recalling the basic terminology of symplectic geometry. Recall from Definition \ref{def_nondeg} that if $X$ is either a graded commutative, associative or Lie algebra then a homogeneous 2-form $\omega \in \drtf{X}$ is nondegenerate if and only if the map $\Phi:\Der(X) \to \drof{X}$ of degree $|\omega|$ defined by the formula
\begin{equation} \label{eqn_ofmvec}
\Phi(\xi):= i_\xi(\omega), \quad \xi \in \Der(X)
\end{equation}
is bijective. When $X$ is the commutative, associative or Lie algebra \csalg{V}, \ctalg{V} or \clalg{V} respectively then we will denote this map by \phiso[Com], \phiso[Ass] or \phiso[Lie] respectively.

\begin{defi} \label{def_symfrm}
Let $X$ be either a formal graded commutative, associative or Lie algebra and let $\omega \in \drtf{X}$ be a homogeneous 2-form. We say $\omega$ is a symplectic form if it is closed (i.e. $dw=0$) and nondegenerate.
\end{defi}

\begin{defi} \label{def_symvec}
Let both $X$ and $X'$ be either formal graded commutative, associative or Lie algebras and let $\omega \in \drtf{X}$ and $\omega' \in \drtf{X'}$ be homogeneous symplectic forms:
\begin{enumerate}
\item[(i)]
We say a vector field $\xi:X \to X$ is a \emph{symplectic vector field} if $L_\xi(\omega)=0$.
\item[(ii)]
We say a diffeomorphism $\phi:X \to X'$ is a \emph{symplectomorphism} if $\phi^*(\omega)=\omega'$.
\end{enumerate}
\end{defi}

We have the following simple proposition about symplectic vector fields which can be found in \cite{ginzsg}:

\begin{prop} \label{prop_symofm}
Let $X$ be either a formal graded commutative, associative or Lie algebra and let $\omega \in \drtf{X}$ be a symplectic form. Then the map $\Phi$ defined by equation \eqref{eqn_ofmvec} induces a one-to-one correspondence between \emph{symplectic} vector fields and \emph{closed} 1-forms:
\[ \Phi: \{ \xi \in \Der(X) : L_\xi(\omega) = 0 \} \to \{ \alpha \in \drof{X} : d\alpha = 0 \}. \]
\end{prop}

\begin{proof}
Let $\xi:X \to X$ be a vector field. By Lemma \ref{lem_schf} part (i) we see that $L_\xi(\omega)=0$ if and only if $di_\xi(\omega) = 0$.
\end{proof}

We will now formulate the Darboux theorem which says that on a `flat manifold' there is a canonical form for any symplectic form. A proof of this can be found in \cite[\S6]{ginzsg}:

\begin{theorem} \label{thm_drboux}
Let $W$ be a free graded module of finite rank and let $X$ be one of the three algebras $\widehat{S}W$, $\widehat{T}W$ or $\widehat{L}W$. Let $\omega \in \drtf{X}$ be a homogeneous 2-form and represent $\omega$ as
\[ \omega = \omega_0 + \omega_1 + \omega_2 + \ldots + \omega_n + \ldots, \]
where $\omega_i$ is a 2-form of order $i$:
\begin{enumerate}
\item[(i)]
$\omega$ is nondegenerate if and only if $\omega_0$ is nondegenerate.
\item[(ii)]
Suppose that $\omega$ is in fact a symplectic form, then there exists a diffeomorphism $\phi \in \Aut(X)$ such that $\phi^*(\omega) = \omega_0$.
\end{enumerate}
\end{theorem}
\noproof

We will now recall (cf. \cite{igusa} and \cite{kontfd}) the definition of an $\infty$-algebra with an invariant inner product. Recall from Remark \ref{rem_infdef} that $\infty$-structures on a free graded module $V$ can be described in terms of systems of maps
\[ \check{m}_i:V^{\otimes i} \to V, \quad i \geq 1 \]
satisfying certain higher homotopy axioms:

\begin{defi} \label{def_invinn}
Let $V$ be a free graded module of finite rank and let
\[ \innprod{V} \]
be an inner product on $V$. 
Let
\[ \check{m}_i:V^{\otimes i} \to V, \quad i \geq 1 \]
be a system of maps determining either an \li, {\ai} or \ci-structure on $V$. We say that this structure is invariant with respect to the inner product $\langle -,- \rangle$ if the following identity holds for all $x_0,\ldots,x_n \in V$;
\begin{equation} \label{eqn_invinn}
\langle \check{m}_n(x_1,\ldots,x_n),x_0 \rangle = (-1)^{n+|x_0|(|x_1|+\ldots+|x_n|)} \langle \check{m}_n(x_0,\ldots,x_{n-1}),x_n \rangle
\end{equation}
and we say that $V$ is an \li, {\ai} or \ci-algebra with an invariant inner product respectively.
\end{defi}

Given a free graded module $V$ we can define an (inhomogeneous) linear bijection
\[ \kappa:(T\Sigma V)^* \to (TV)^* \]
as follows: Given a linear map $\alpha:\Sigma V^{\otimes n} \to \gf$ of degree $d$ the map $\kappa(\alpha):V^{\otimes n} \to \gf$ of degree $d-n$ is defined by the following commutative diagram:
\begin{equation} \label{fig_tfmaps}
\xymatrix{V^{\otimes n} \ar^{\kappa(\alpha)}[r] & \gf \\ \Sigma V^{\otimes n} \ar^{(\Sigma^{-1})^{\otimes n}}[u] \ar_{\alpha}[ru]}
\end{equation}

It is well known that the \emph{unsigned} action of $S_n$ on $(\Sigma V^{\otimes n})^*$ corresponds to the \emph{signed} action of $S_n$ on $(V^{\otimes n})^*$ under the map $\kappa$. In particular $\kappa$ gives us the following bijection of degree $-2$:
\[ \kappa:(\Lambda^2 \Sigma V)^* \to (S^2 V)^* \]

We will now prove the main result of this section which links \emph{invariant} $\infty$-structures to \emph{symplectic} vector fields:

\begin{theorem} \label{thm_syminv}
Let $V$ be a free graded module of finite rank. Let
\begin{equation} \label{eqn_syminvdummy}
\check{m}_i:V^{\otimes i} \to V, \quad i\geq 1
\end{equation}
be either an \li, {\ai} or \ci-structure on $V$ and let
\[ \innprod{V} \]
be an inner product on $V$.

Consider the vector field $m$ (belonging to $\Der(\csalg{V})$, $\Der(\ctalg{V})$ or $\Der(\clalg{V})$)
corresponding to the \li, {\ai} or \ci-structure respectively and consider the symplectic form $\omega$ (belonging to \drtf[Com]{\csalg{V}}, \drtf[Ass]{\ctalg{V}} or \drtf[Lie]{\clalg{V}} respectively) corresponding to the inner product $\langle -,- \rangle$; $\langle -,- \rangle = \kappa\zeta(\omega)$. The $\infty$-structure \eqref{eqn_syminvdummy} is invariant with respect to the inner product $\langle -,- \rangle$ if and only if $L_m \omega = 0$.
\end{theorem}

\begin{proof}
Start by choosing a basis $x_1,\ldots,x_n$ of the free graded module $V$ and let $y_i:=\Sigma x_i$, then there exist coefficients $a_{ij} \in \gf$ such that
\[ \omega = \sum_{1\leq i,j \leq n} a_{ij}dy_i^* dy_j^* \]
regardless of whether $\omega$ belongs to \drtf[Com]{\csalg{V}}, \drtf[Ass]{\ctalg{V}} or \drtf[Lie]{\clalg{V}}. The corresponding inner product $\langle -,- \rangle:= \kappa\zeta(\omega)$ is given by the formula,
\begin{equation} \label{eqn_syminvdummya}
\langle a,b \rangle = -\sum_{1 \leq i,j \leq n} a_{ij}\left[x_i^*(a)x_j^*(b) +(-1)^{|x_i||x_j|}x_j^*(a)x_i^*(b)\right].
\end{equation}

Let us treat the {\ai} case first since it is the simplest. Let
\[ m_i:\Sigma V^{\otimes i} \to \Sigma V, \quad i\geq 1 \]
be the system of maps defined by figure \eqref{fig_infdef} which corresponds to the system of maps in \eqref{eqn_syminvdummy} which define the \ai-structure: The vector field $m \in \Der(\ctalg{V})$ then has the form
\[ m = m_1^* + m_2^* + \ldots + m_n^* + \ldots \]
where the map $m_k^*: \Sigma V^* \to (\Sigma V^*)^{\cotimes k}$ is extended uniquely to a vector field on \ctalg{V}, cf. Proposition \ref{prop_dergen}. We then calculate $L_m(\omega)$ as
\[ L_m(\omega) = -d i_m(\omega) = -d\left( \sum_{k=1}^\infty \sum_{1 \leq i,j \leq n} a_{ij}\left[ (-1)^{|x_i||x_j|}dy_j^*\cdot (y_i^* \circ m_k) + dy_i^* \cdot (y_j^* \circ m_k) \right] \right). \]

Recall that the \emph{unsigned} action of $S_n$ on $(\Sigma V^*)^{\cotimes n}$ corresponds to the \emph{signed} action of $S_n$ on $(V^*)^{\cotimes n}$ under the map $\kappa$ defined by diagram \eqref{fig_tfmaps}. Let \zetiso[Ass] be the isomorphism defined by Theorem \ref{thm_tfmaps}: Since $L_m\omega$ is a closed 2-form, this means that $L_m(\omega) = 0$ if and only if $\kappa\zetiso[Ass]L_m(\omega) = 0$. It follows from the definition of \zetiso[Ass] and $\kappa$ that $L_m \omega = 0$ if and only if for all $k \geq 1$,
\[ (1-z_k)\cdot\sum_{1 \leq i,j \leq n}a_{ij}\left[(-1)^{k|x_j|}(-1)^{|x_i||x_j|}x_j^*\cdot x_i^*\circ \check{m}_k + (-1)^{k|x_i|}x_i^*\cdot x_j^*\circ \check{m}_k\right] = 0. \]
Since $\langle -,- \rangle$ is graded symmetric it follows from equation \eqref{eqn_syminvdummya} and a simple calculation that the above holds if and only if the following map is invariant under the (signed) action of the cyclic group;
\begin{displaymath}
\begin{array}{ccc}
V^{\otimes k+1} & \to & \gf, \\
u_0 \otimes \ldots \otimes u_k & \mapsto & \langle m_k(u_0,\ldots,u_{k-1}),u_k \rangle. \\
\end{array}
\end{displaymath}
The invariance of the above map under the action of the cyclic group is condition \eqref{eqn_invinn} of Definition \ref{def_invinn} and hence we have established the theorem for \ai-algebras.

Let us treat the {\ci} case next. This will follow as a corollary of the theorem for \ai-algebras. Let $l:\dr[Lie]{\clalg{V}} \to \dr[Ass]{\ctalg{V}}$ be the map defined by equation \eqref{eqn_liasmp} and let $m \in \Der(\clalg{V})$ be the \ci-structure on $V$. By lemmas \ref{lem_clsdtf} and \ref{lem_plcomm} we see that
\[ L_m(\omega) =0 \Leftrightarrow lL_m(\omega) = 0 \Leftrightarrow L_m l(\omega) = 0. \]
By the preceding argument $L_m l(\omega) = 0$ if and only if the inner product
\[ \langle -,- \rangle := \kappa\zetiso[Lie](\omega) = \kappa\zetiso[Ass]l(\omega) \]
is invariant with respect to the system of maps in \eqref{eqn_syminvdummy} which describe the \ci-structure. Again this is condition \eqref{eqn_invinn} of Definition \ref{def_invinn} and hence we have proven the theorem for \ci-algebras.

Lastly we will treat the {\li} case. Let $j:\drof[Com]{\csalg{V}} \to \drof[Ass]{\ctalg{V}}$ be the map defined by diagram \eqref{fig_caofmp} and let \zetiso[Com] be the isomorphism defined by Theorem \ref{thm_tfmaps}. Let 
\[ m_i: S^i\Sigma V \to \Sigma V, \quad i \geq 1 \]
be the system of maps defined by figure \eqref{fig_infdef} which corresponds to the system of maps in \eqref{eqn_syminvdummy} which define the \li-structure: The vector field $m \in \Der(\csalg{V})$ then satisfies the identity,
\[ i \circ m \circ \pi = \sum_{k=1}^\infty m_k^*; \]
where the map $m_k^*: \Sigma V^* \to (S^k\Sigma V)^*$ is extended uniquely to a vector field on $(S\Sigma V)^*$ (which is equipped with the shuffle product) and the maps $i$ and $\pi$ are defined by equations \eqref{eqn_liduin} and \eqref{eqn_lidupr} respectively.

Now $L_m(\omega)=0$ if and only if $\kappa\zetiso[Com]L_m(\omega) = 0$. We calculate $\kappa\zetiso[Com]L_m(\omega)$ as;
\begin{displaymath}
\begin{split}
\kappa\zetiso[Com]L_m(\omega) & = -\kappa\zetiso[Com]di_m(\omega) = -\kappa\zetiso[Ass]dji_m(\omega), \\
& = -\kappa\zetiso[Ass]dj\left(\sum_{1 \leq i,j \leq n} a_{ij}\left[ (-1)^{|x_i||x_j|}dy_j^*\cdot m(y_i^*) + dy_i^* \cdot m(y_j^*) \right]\right), \\
& = -\kappa\zetiso[Ass]d\left(\sum_{1 \leq i,j \leq n} a_{ij}\left[ (-1)^{|x_i||x_j|}dy_j^*\cdot im(y_i^*) + dy_i^* \cdot im(y_j^*) \right]\right), \\
& = -\kappa\zetiso[Ass]d\left(\sum_{k=1}^\infty\sum_{1 \leq i,j \leq n} a_{ij}\left[ (-1)^{|x_i||x_j|}dy_j^*\cdot (y_i^* \circ m_k) + dy_i^* \cdot (y_j^* \circ m_k) \right]\right). \\
\end{split}
\end{displaymath}
Again it follows from the definition of $\zetiso[Ass]$ and the same simple calculations that $L_m(\omega)=0$ if and only if \eqref{eqn_invinn} holds for all $k \geq 1$. This proves the theorem for \li-algebras.
\end{proof}

This result motivates the following definition of a \emph{symplectic} $\infty$-algebra:

\begin{defi}
Let $V$ be a free graded module of finite rank:
\begin{enumerate}
\item[(i)]
A \emph{symplectic} \li-structure on $V$ is a symplectic form $\omega \in \drtf[Com]{\csalg{V}}$ together with a \emph{symplectic} vector field
\[ m:\csalg{V} \to \csalg{V} \]
of degree one and vanishing at zero; such that $m^2=0$.
\item[(ii)]
A \emph{symplectic} \ai-structure on $V$ is a symplectic form $\omega \in \drtf[Ass]{\ctalg{V}}$ together with a \emph{symplectic} vector field
\[ m:\ctalg{V} \to \ctalg{V} \]
of degree one and vanishing at zero; such that $m^2=0$.
\item[(iii)]
A \emph{symplectic} \ci-structure on $V$ is a symplectic form $\omega \in \drtf[Lie]{\clalg{V}}$ together with a \emph{symplectic} vector field
\[ m:\clalg{V} \to \clalg{V} \]
of degree one; such that $m^2=0$.
\end{enumerate}
\end{defi}

\begin{rem}
Of course, in order to avoid homotopy non invariant constructions we should also insist that the underlying complex of our symplectic $\infty$-algebra be cellular, cf. Remark \ref{rem_infcel}.
\end{rem}

\begin{defi}
Let $V$ and $U$ be free graded modules of finite ranks:
\begin{enumerate}
\item[(i)]
Let $(m,\omega)$ and $(m',\omega')$ be \emph{symplectic} \li-structures on $V$ and $U$ respectively: A \emph{symplectic} \li-morphism from $V$ to $U$ is a continuous algebra homomorphism
\[ \phi:\csalg{U} \to \csalg{V} \]
of degree zero such that $\phi \circ m' = m \circ \phi$ and $\phi^*(\omega') = \omega$.
\item[(ii)]
Let $(m,\omega)$ and $(m',\omega')$ be \emph{symplectic} \ai-structures on $V$ and $U$ respectively: A \emph{symplectic} \ai-morphism from $V$ to $U$ is a continuous algebra homomorphism
\[ \phi:\ctalg{U} \to \ctalg{V} \]
of degree zero such that $\phi \circ m' = m \circ \phi$ and $\phi^*(\omega') = \omega$.
\item[(iii)]
Let $(m,\omega)$ and $(m',\omega')$ be \emph{symplectic} \ci-structures on $V$ and $U$ respectively: A \emph{symplectic} \ci-morphism from $V$ to $U$ is a continuous algebra homomorphism
\[ \phi:\clalg{U} \to \clalg{V} \]
of degree zero such that $\phi \circ m' = m \circ \phi$ and $\phi^*(\omega') = \omega$.
\end{enumerate}
\end{defi}

\begin{rem}
We have the following diagram of functors,
\[ \xymatrix{ \ci\mathrm{-}algebras \ar[r] & \ai\mathrm{-}algebras \ar[r] & \li\mathrm{-}algebras \\ \emph{Symplectic }\ci\mathrm{-}algebras \ar[r] \ar[u] & \emph{Symplectic }\ai\mathrm{-}algebras \ar[r] \ar[u] & \emph{Symplectic }\li\mathrm{-}algebras \ar[u] \\ } \]
where the top row is just diagram \eqref{fig_calseq} and the vertical arrows are the forgetful functors. In order to define the arrows in the bottom row it only remains to describe what happens to the symplectic forms in the symplectic $\infty$-algebras. This mapping is provided by the maps $l$ and $p$ defined in section \ref{relations}, cf. \eqref{eqn_liasmp} and \eqref{eqn_ascomp}:
\[ \xymatrix{ \drtf[Lie]{\clalg{V}} \ar^{l}[r] & \drtf[Ass]{\ctalg{V}} \ar^{p}[r] & \drtf[Com]{\csalg{V}} \\ }. \]
\end{rem}

We have the following useful lemmas regarding the cohomology of a symplectic $\infty$-algebra:

\begin{lemma} \label{lem_duaiso}
Let $(V,m,\omega)$ be
\begin{enumerate}
\item[(a)] a \emph{symplectic} \li-algebra,
\item[(b)] a \emph{symplectic} \ai-algebra or
\item[(c)] a \emph{symplectic} \ci-algebra.
\end{enumerate}
The map $\Phi$ defined by equation \eqref{eqn_ofmvec} is an isomorphism of complexes:
\begin{enumerate}
\item[(a)] $\phiso[Com]:\clac{V}{V} \to \Sigma^{|\omega|-2}\clac{V}{V^*}$.
\item[(b)] $\phiso[Ass]:\choch{V}{V} \to \Sigma^{|\omega|-2}\choch{V}{V^*}$.
\item[(c)] $\phiso[Lie]:\caq{V}{V} \to \Sigma^{|\omega|-2}\caq{V}{V^*}$.
\end{enumerate}
\end{lemma}

\begin{proof}
We shall treat all three cases simultaneously. By Definition \ref{def_nondeg} the map $\Phi$ must be bijective. It only remains to prove that this map commutes with the differentials:
          
Let $X$ be the algebra \csalg{V}, \ctalg{V} or \clalg{V} depending on whether $V$ is a \emph{symplectic} \li, {\ai} or \ci-algebra respectively. The complex computing the cohomology of $V$ with coefficients in $V$ is
\[ C^\bullet(V,V):= \left(\Sigma^{-1}\Der(X),\ad m\right). \]
The complex computing the cohomology of $V$ with coefficients in $V^*$ is
\[ C^\bullet(V,V^*):= \left(\Sigma\drof{X},L_m\right). \]

Since $m$ is a symplectic vector field, Lemma \ref{lem_schf} implies that $\Phi$ respects the differentials:
\begin{displaymath}
\begin{split}
\Phi(\ad m(\xi)) & = i_{[m,\xi]}(\omega) = [L_m,i_\xi](\omega), \\
& = L_m i_\xi(\omega) = L_m\Phi(\xi). \\
\end{split}
\end{displaymath}
\end{proof}

We need an analogue of Lemma \ref{lem_duaiso} in the context of the cyclic theories:

\begin{lemma} \label{lem_duciso}
Let $(V,m,\omega)$ be
\begin{enumerate}
\item[(a)] a \emph{symplectic} \li-algebra,
\item[(b)] a \emph{symplectic} \ai-algebra or
\item[(c)] a \emph{symplectic} \ci-algebra;
\end{enumerate}
then the Lie subalgebra of
\begin{enumerate}
\item[(a)] $\clac{V}{V}:=\left(\Der(\csalg{V}),\ad m\right)$,
\item[(b)] $\choch{V}{V}:=\left(\Der(\ctalg{V}),\ad m\right)$ or
\item[(c)] $\caq{V}{V}:=\left(\Der(\clalg{V}),\ad m\right)$
\end{enumerate}
consisting of all \emph{symplectic} vector fields forms a subcomplex denoted by
\begin{enumerate}
\item[(a)] $S\clac{V}{V}$,
\item[(b)] $S\choch{V}{V}$ or
\item[(c)] $S\caq{V}{V}$
\end{enumerate}
respectively. Furthermore, there is an isomorphism
\begin{enumerate}
\item[(a)] $\upiso[Com]:\Sigma^{|\omega|-2}\cclac{V} \to S\clac{V}{V}$,
\item[(b)] $\upiso[Ass]:\Sigma^{|\omega|-2}\cchoch{V} \to S\choch{V}{V}$ or
\item[(c)] $\upiso[Lie]:\Sigma^{|\omega|-2}\ccaq{V} \to S\caq{V}{V}$
\end{enumerate}
defined by the formula,
\begin{equation} \label{eqn_duciso}
\Upsilon(\alpha):=(-1)^{|\alpha|}\Phi^{-1}d(\alpha);
\end{equation}
where $\Phi$ is the map defined by equation \eqref{eqn_ofmvec}.
\end{lemma}

\begin{proof}
Again the proof is the same in all three cases. Since the $\infty$-structure $m$ is a symplectic vector field, the Lie subalgebra of symplectic vector fields does indeed form a subcomplex as claimed.

Lemma \ref{lem_schf} part (v) and Lemma \ref{lem_duaiso} tell us that $\Upsilon$ is a map of complexes, i.e. it respects the differentials. It remains to prove that $\Upsilon$ is bijective: This follows from Lemma \ref{lem_poinca} and Proposition \ref{prop_symofm}.
\end{proof}

\section{Examples} \label{examples}

In this section we give a few explicit examples of \emph{symplectic} $\infty$-algebras (or in fact $\infty$-algebras with an invariant inner product). A large class of examples of symplectic $C_\infty$-algebras is provided by taking minimal models of the cochain algebras of Poincar\'e duality spaces. This follows from Theorem \ref{thm_main} and we will have more to say about these algebras in section \ref{strings}.

Let us now describe a series of symplectic \ai-algebras which were introduced (without symplectic structure) in the second author's paper \cite{lazmod} under the name `Moore algebras'. It should be mentioned that the Moore algebras are not $A_\infty$-algebras under Definition \ref{def_infstr} but rather $\mathbb{Z}/2$-graded (or super) $A_\infty$-algebras. This means that their underlying modules are $\mathbb{Z}/2$-graded rather than $\mathbb{Z}$-graded and the $A_\infty$-structure will simply be an odd vector field having square $0$. Clearly any $\mathbb{Z}$-graded $\infty$-algebra determines a super $\infty$-algebra. Conversely, a super $\infty$-algebra gives rise to a $\mathbb{Z}$-graded one, albeit over a larger ring $\gf[v^{\pm 1}]$ where $|v|=2$. All our results carry over in an obvious fashion to the supercase. 

`Even' Moore algebras give examples of symplectic \ai-algebras with odd symplectic forms. Consider the free module
\[ V:= \Sigma^{-1}\gf \oplus \gf \]
with generators $y$ and $1$ in degrees 1 and 0 respectively. We have $\ctalg{V}=\gf\langle\langle\tau,t\rangle\rangle$ where $t$ is dual to the generator $\Sigma y$ and $\tau$ is dual to the generator $\Sigma 1$. In \cite{lazmod} the second author proved that any vector field $m$ of the form
\[ m = \sum_{i=1}^\infty u_i t^i\partial_\tau + \ad\tau - \tau^2\partial_\tau \]
is an \ai-structure on $V$ (i.e. $m^2=0$). We endow this module with an odd symplectic form
\[ \omega:= d\tau dt. \]
The corresponding inner product $\langle -,- \rangle:=\kappa\zetiso[Ass](\omega)$ is given by the formula;
\begin{displaymath}
\begin{array}{ccccc}
\langle y,1 \rangle & = & \langle 1,y \rangle & = & -1, \\
\langle y,y \rangle & = & \langle 1,1 \rangle & = & 0. \\
\end{array}
\end{displaymath}
Note that if we disregard higher multiplications, this even Moore algebra is essentially the cohomology algebra of an odd-dimensional sphere together with its Poincar\'e duality form. 

Examples of symplectic \ai-algebras with even symplectic forms are provided by `odd' Moore algebras. Consider the free module
\[ V:= \gf \oplus \gf \]
with generators $y$ and $1$ in degree 0. Again we have $\ctalg{V}=\gf\langle\langle\tau,t\rangle\rangle$ where $t$ is dual to the generator $\Sigma y$ and $\tau$ is dual to the generator $\Sigma 1$. In \cite{lazmod} it is proved that any vector field $m$ of the form
\[ m = \sum_{i=1}^\infty v_i t^{2i}\partial_t + \sum_{i=1}^\infty w_i t^{2i}\partial_\tau + \ad\tau - \tau^2\partial_\tau \]
is an \ai-structure on $V$. Here there is a choice for an even symplectic form. One option is
\[ \omega := d\tau d\tau + dtdt. \]
The corresponding inner product $\langle -,- \rangle:=\kappa\zetiso[Ass](\omega)$ is given by the formula;
\begin{displaymath}
\begin{array}{ccccc}
\langle 1,1 \rangle & = & \langle y,y \rangle & = & -1, \\
\langle y,1 \rangle & = & \langle 1,y \rangle & = & 0. \\
\end{array}
\end{displaymath}
The nonunital version of this $A_\infty$-algebra was introduced in \cite{kontfd} and was shown to give rise to the Morita-Miller-Mumford classes in the cohomology of moduli spaces of complex algebraic curves.

The other option is to take
\[ \omega:= d\tau dt.\]
The corresponding inner product $\langle -,- \rangle:=\kappa\zetiso[Ass](\omega)$ will then be given by the formula;
\begin{displaymath}
\begin{array}{ccccc}
\langle 1,1 \rangle & = & \langle y,y \rangle & = & 0, \\
\langle y,1 \rangle & = & \langle 1,y \rangle & = & 1. \\
\end{array}
\end{displaymath}
Note that (disregarding higher multiplications) the odd Moore algebra is essentially the cohomology algebra of an even-dimensional sphere; moreover the second choice of the symplectic form corresponds to the Poincar\'e duality form.

\begin{prop}
\
\begin{enumerate}
\item[(i)]
Consider the free graded module
\[ V:= \Sigma^{-1}\gf \oplus \gf \]
endowed with the symplectic form
\[ \omega:= d\tau dt. \]
The \ai-structure
\[ m := \sum_{i=1}^\infty u_i t^i\partial_\tau + \ad\tau - \tau^2\partial_\tau \]
is a symplectic vector field and hence gives $V$ the structure of an \ai-algebra with an (odd) invariant inner product.
\item[(ii)]
Consider the free graded module
\[ V:= \gf \oplus \gf. \]
\begin{enumerate}
\item
Endow $V$ with the symplectic form
\[ \omega:= d\tau d\tau + dtdt. \]
The \ai-structure
\[ m := \sum_{i=1}^\infty v_i t^{2i}\partial_t + t^2\partial_\tau + \ad\tau - \tau^2\partial_\tau \]
is a symplectic vector field and hence gives $V$ the structure of an \ai-algebra with an (even) invariant inner product.
\item
Endow $V$ with the symplectic form
\[ \omega:= d\tau dt. \]
The \ai-structure
\[ m:= \sum_{i=1}^\infty w_i t^{2i}\partial_\tau + \ad\tau-\tau^2\partial_\tau \]
is a symplectic vector field and hence gives $V$ the structure of an \ai-algebra with an (even) invariant inner product.
\end{enumerate}
\end{enumerate}
\end{prop}

\begin{proof}
\
\begin{enumerate}
\item[(i)]
We calculate $L_m \omega$ as
\begin{displaymath}
\begin{split}
L_m \omega & = -\left(d(u(t)+\tau^2)\cdot dt + d\tau\cdot d[\tau,t]\right), \\
& = -\left(d(u(t))\cdot dt + d\tau\cdot\tau\cdot dt - \tau\cdot d\tau\cdot dt + d\tau\cdot[d\tau,t] - d\tau\cdot[\tau,dt]\right), \\
& = -\left(d(u(t))\cdot dt + d\tau\cdot\tau\cdot dt - \tau\cdot d\tau\cdot dt + [d\tau,d\tau]\cdot t - d\tau\cdot\tau\cdot dt - d\tau\cdot dt\cdot \tau \right), \\
& = -d(u(t))\cdot dt. \\
\end{split}
\end{displaymath}
From Theorem \ref{thm_tfmaps} it follows that $d(u(t))\cdot dt=0$ and therefore $m$ is a symplectic vector field as claimed.
\item[(ii)]
\
\begin{enumerate}
\item
We calculate $L_m \omega$ as
\begin{displaymath}
\begin{split}
L_m \omega = & -2\left(dt\cdot d(v(t)+[\tau,t]) + d\tau\cdot d(\tau^2 + t^2)\right), \\
= & -2\left(dt\cdot d(v(t)) + dt\cdot[d\tau,t] - dt\cdot[\tau,dt]\right) \\
& -2\left(d\tau\cdot d\tau\cdot\tau - d\tau\cdot\tau\cdot d\tau + d\tau\cdot dt\cdot t - d\tau\cdot t\cdot dt\right), \\
= & -2\left(dt\cdot d(v(t)) + dt\cdot d\tau\cdot t - dt\cdot t\cdot d\tau + [dt,dt]\cdot\tau + d\tau\cdot dt\cdot t - d\tau\cdot t\cdot dt\right), \\
= & -2dt\cdot d(v(t)). \\
\end{split}
\end{displaymath}
Again, it follows from Theorem \ref{thm_tfmaps} that $dt\cdot d(v(t)) = 0$ and hence $m$ is a symplectic vector field.
\item
Again we calculate $L_m \omega$ as
\begin{displaymath}
\begin{split}
L_m \omega & = -d(\tau^2)dt-d(w(t))dt - d\tau d[\tau,t], \\
& = d\tau \cdot\tau\cdot dt - \tau d\tau dt + d\tau[d\tau,t] - d\tau[\tau,dt] - d(w(t))dt, \\
& = - d(w(t))dt. \\
\end{split}
\end{displaymath}
As was shown above, $d(w(t))dt = 0$ and therefore $m$ is a symplectic vector field.
\end{enumerate}
\end{enumerate}
\end{proof}
We can construct further examples of symplectic \li-algebras from \emph{even} Moore algebras (it is easy to see that the odd ones lead to trivial examples of such). Namely, consider the commutative superalgebra $\gf[[t]]\otimes \Lambda_\gf(\tau)$ together with the vector field $m:=\sum_{i=1}^\infty u_it^i\partial_{\tau}$. Then $m^2=0$ and $m$ preserves the symplectic form $\omega=dtd\tau$. In other words we have a symplectic \li-algebra on the supermodule $\gf^{1|1}$. 

The last symplectic $L_\infty$-algebra admits a curious topological interpretation. Namely, set $u_i=0$ for $i\neq k>1$ and $u_k=1$. In other words, we have a \li-structure on the module $\gf\langle x\rangle\oplus\gf\langle y\rangle$ where $x$ is odd whereas $y$ is even. This \li-structure is determined by $m_k(x,\ldots,x)=y$ and the other $m_i$'s vanish. Denote this $L_\infty$-algebra by $L$. It is easy to see that $L$ is weakly equivalent to the Lie model of the rational homotopy type of $\mathbb{C}P^{k-1}$, the complex projective space of dimension $k-1$. Here $x$ corresponds to the generator of $\pi_2(\mathbb{C}P^{k-1})$ and $y$ corresponds to the generator of $\pi_{2k-1}(\mathbb{C}P^{k-1})$. It would be interesting to find out what topological structure on $\mathbb{C}P^{k-1}$ corresponds to the invariant inner product on $L$.

\section{Obstruction theory for Symplectic \ci-Algebras} \label{sec_symobs}

In this section we develop the obstruction theory in the symplectic context which is parallel to our treatment in section \ref{sec_obstrc}. Once again it is possible to develop an obstruction theory for symplectic $L_n$ and $A_n$-algebras but we shall only discuss the theory for symplectic $C_n$-algebras. This is in order to avoid repeating ourselves and since the goal of this section is to lay the groundwork for section \ref{correspondence} in which we prove our main result and this result will only hold for \ci-algebras.

In this section we shall only consider \emph{symplectic} $C_n$-algebras \emph{whose symplectic form $\omega$ has order zero}. This is not really a significant restriction however since the Darboux theorem (Theorem \ref{thm_drboux}) tells us that any symplectic $C_n$-algebra is isomorphic to one of this form. It follows from Theorem \ref{thm_syminv} that a strictly graded commutative \emph{symplectic} \ci-algebra $A:=(V,m_2,\omega)$ is the same thing as a commutative Frobenius algebra and we hereafter refer to them as such.

Recall that given a \emph{symplectic} \ci-algebra $A:=(V,m,\omega)$, Lemma \ref{lem_duciso} gives us an isomorphism \upiso[Lie] between the (shifted) cyclic Harrison complex $\Sigma^{|\omega|-2}\ccaq{A}$ and the subcomplex of \caq{A}{A} consisting of all symplectic vector fields which we denoted by $S\caq{A}{A}$. Furthermore if $A$ is in fact a strictly graded commutative Frobenius algebra then \upiso[Lie] respects the natural bigrading induced on the cohomology defined in remarks \ref{rem_ordgde} and \ref{rem_ordgrd}:
\[ \upiso[Lie]: \ccaq[n+1,j+|\omega|-2]{A} \to S\caq[nj]{A}{A}. \]

Recall that given a free graded module $V$, the module consisting of all continuous endomorphisms
\[ f:\clalg{V}\to\clalg{V} \]
of order $\geq n$ is denoted by $(n)$.

\subsection{Obstruction theory for symplectic $C_n$-algebra structures}

In this section we will formulate the \emph{symplectic} analogues of the theorems and definitions given in section \ref{sec_obsalg}.

\begin{defi} \label{def_minsym}
Let $V$ be a free graded module of finite rank. A \emph{minimal symplectic} $C_n$-structure ($n \geq 3$) on $V$ consists of a \emph{symplectic} form $\omega \in \drtf[Lie]{\clalg{V}}$ \emph{of order zero} and a \emph{symplectic} vector field $m \in \Der(\clalg{V})$ of degree one which has the form
\[ m = m_2 + \ldots + m_{n-1}, \quad m_i \text{ has order $i$} \]
and satisfies the condition $m^2 = 0 \mod (n+1)$.

Let $m$ and $m'$ be two minimal symplectic $C_n$-structures on $V$. We say $m$ and $m'$ are equivalent if there is a \emph{symplectomorphism} of degree zero $\phi \in \Aut(\clalg{V})$ of the form
\begin{equation} \label{eqn_sdifrm}
\phi = \id + \phi_2 + \phi_3 + \ldots + \phi_k + \ldots
\end{equation}
where $\phi_i$ is an endomorphism of order $i$ and such that
\[ \phi \circ m \circ \phi^{-1} = m' \mod (n). \]
\end{defi}

\begin{rem}
Following the convention made in section \ref{sec_notcon} we shall call a symplectomorphism of the form \eqref{eqn_sdifrm} a \emph{pointed symplectomorphism}. Note that under the definition of equivalence, it would be impossible to have two equivalent minimal symplectic $C_n$-structures on $V$ which had different symplectic forms, hence the reason for the omission of the symplectic forms in the definition.
\end{rem}

\begin{rem}
Obviously if $m=m_2+\ldots+m_{n-1}$ is a minimal symplectic $C_n$-algebra with symplectic form $\omega$ then $A:=(V,m_2,\omega)$ is a strictly graded commutative Frobenius algebra which we will call the underlying Frobenius algebra. Observe that two equivalent minimal symplectic $C_n$-structures have the same underlying Frobenius algebra.
\end{rem}

It will be useful to introduce the following definition:

\begin{defi}
Let $A:=(V,\mu_2,\omega)$ be a strictly graded commutative Frobenius algebra, then the moduli space of minimal \emph{symplectic} $C_n$-structures on $V$ fixing $\mu_2$ and $\omega$ is denoted by \mstr[S]{C}{n}{A} and defined as the quotient of the set
\[ \{ m:\clalg{V} \to \clalg{V} : m \text{ is a minimal \emph{symplectic} $C_n$-structure with respect to $\omega$ and $m_2=\mu_2$} \} \]
by the equivalence relation defined in Definition \ref{def_minsym}.
\end{defi}

We will now describe the appropriate terminology necessary in discussing extensions of symplectic $C_n$-structures to structures of higher order:

\begin{defi}
Let $V$ be a free graded module of finite rank and let $\omega \in \drtf[Lie]{\clalg{V}}$ be a symplectic form of order zero. We say that a minimal symplectic $C_n$-structure $m = m_2 + \ldots + m_{n-1}$ on $V$ is extendable if there exists a \emph{symplectic} vector field $m_n \in \Der(\clalg{V})$ of order $n$ such that $m + m_n$ is a symplectic $C_{n+1}$-structure on $V$ and we call $m_n$ an extension of $m$.

Let $m_n$ and $m'_n$ be two extensions of $m$. We say that $m_n$ and $m'_n$ are equivalent if there exists a \emph{symplectomorphism} $\phi \in \Aut(\clalg{V})$ of degree zero of the form
\[ \phi = \id + \phi_{n-1} + \phi_n + \ldots + \phi_{n+k} + \ldots \]
where $\phi_i$ is an endomorphism of order $i$ and such that
\[ \phi \circ (m + m_n) \circ \phi^{-1} = m + m'_n \mod (n+1). \]
The quotient of the set of all extensions of $m$ by this equivalence relation will be denoted by \mext[S]{n}{m}.
\end{defi}

We now identify the appropriate obstruction to extending \emph{symplectic} $C_n$-structures. We will use the same notation as in Definition \ref{def_obstrc}. We will later justify this abuse of notation by Lemma \ref{lem_obscom}:

\begin{defi}
Let $V$ be a free graded module of finite rank and let $\omega \in \drtf[Lie]{\clalg{V}}$ be a symplectic form of order zero. Let $m = m_2+\ldots + m_{n-1}$ be a minimal symplectic $C_n$-structure on $V$. We define the 0-form $\Obs(m)$ of order $n+2$ and degree $2+|\omega|$ by
\begin{equation} \label{eqn_symobs}
\Obs(m) := \frac{1}{2}\upiso[Lie]^{-1}\left(\sum_{\begin{subarray}{c} i+j = n+2 \\ 3\leq i,j \leq n-1 \end{subarray}} [m_i,m_j] \right).
\end{equation}
\end{defi}

We can now formulate the symplectic analogues of theorems \ref{thm_obstrc} and \ref{thm_extend}. We shall omit the proofs since they are the same as the proofs of section \ref{sec_obstrc} verbatim except that we must use the fact from Lemma \ref{lem_duciso} that $\upiso[Lie]$ is an isomorphism.

\begin{theorem} \label{thm_symobs}
Let $A:=(V,m_2,\omega)$ be a strictly graded commutative Frobenius algebra. For all $n \geq 3$ \eqref{eqn_symobs} induces a map
\begin{displaymath}
\begin{array}{ccc}
\mstr[S]{C}{n}{A} & \to & \hcaq[n+2,1+|\omega|]{A}, \\
m & \mapsto & \Obs(m); \\
\end{array}
\end{displaymath}
which we will denote by $\Obs_n$. The kernel of this map consists of precisely those symplectic $C_n$-structures which are extendable:
\[ \{ m \in \mstr[S]{C}{n}{A} : m \textnormal{ is extendable} \} = \ker(\Obs_n). \]
\end{theorem}
\noproof

\begin{theorem} \label{thm_symext}
Let $A:=(V,m_2,\omega)$ be a strictly graded commutative Frobenius algebra and let $m \in \mstr[S]{C}{n}{A}$ be an extendable symplectic $C_n$-structure, then \hcaq[n+1,|\omega|]{A} acts freely and transitively on \mext[S]{n}{m}:
\begin{displaymath}
\begin{array}{ccc}
\hcaq[n+1,|\omega|]{A} \times \mext[S]{n}{m} & \to & \mext[S]{n}{m}, \\
(\alpha , m_n) & \mapsto & m_n + \upiso[Lie](\alpha). \\
\end{array}
\end{displaymath}
\end{theorem}
\noproof

\begin{rem}
The analogues of theorems \ref{thm_symobs} and \ref{thm_symext} apply to symplectic $L_n$ and $A_n$-algebras as well. Cyclic Harrison cohomology is replaced with cyclic Chevalley-Eilenberg and cyclic Hochschild cohomology respectively. Recall that cyclic Chevalley-Eilenberg cohomology is simply Chevalley-Eilenberg cohomology with \emph{trivial coefficients}.
\end{rem}

\subsection{Obstruction theory for symplectic $C_n$-algebra morphisms}

In this section we will develop the obstruction theory for morphisms between two \emph{symplectic} $C_n$-algebras. Again, the formulation is entirely analogous to that of section \ref{sec_obsmor}.

\begin{defi} \label{def_symmor}
Let $V$ be a free graded module of finite rank and let $\omega \in \drtf[Lie]{\clalg{V}}$ be a symplectic form of order zero. Let $m$ and $m'$ be two minimal \emph{symplectic} $C_n$-structures on $V$. A minimal \emph{symplectic} $C_n$-morphism from $m$ to $m'$ is a \emph{symplectomorphism} $\phi \in \Aut(\clalg{V})$ of degree zero such that
\[ \phi \circ m = m'\circ\phi \mod (n). \]

Let $\phi$ and $\phi'$ be two such minimal symplectic $C_n$-morphisms. We say $\phi$ and $\phi'$ are homotopic if there exists a \emph{symplectic} vector field $\eta$ of degree $-1$ such that
\[ \phi = \phi'\circ\exp([m,\eta]) \mod(n-1). \]
\end{defi}

We will now introduce the moduli space of \emph{symplectic} $C_n$-morphisms:

\begin{defi}
Let $V$ be a free graded module of finite rank and let $\omega \in \drtf[Lie]{\clalg{V}}$ be a symplectic form of order zero. Let $m$ and $m'$ be two minimal symplectic $C_n$-structures on $V$ with the same underlying Frobenius algebra (i.e. $m_2 = m'_2$). The moduli space of minimal \emph{symplectic} $C_n$-morphisms from $m$ to $m'$ is denoted by \mmor[S]{n}{m}{m'} and defined as the quotient of the set
\[ \{ \phi:\clalg{V} \to \clalg{V} : \phi \text{ is a minimal \emph{pointed symplectic} $C_n$-morphism from $m$ to $m'$} \} \]
by the homotopy equivalence relation defined in Definition \ref{def_symmor}.
\end{defi}

Let $m = m_2 + \ldots + m_n$ be a symplectic $C_{n+1}$-structure on a free graded module $V$. Recall that we define the corresponding symplectic $C_n$-structure $\bar{m}$ on $V$ as
\[ \bar{m}:= m_2 + \ldots + m_{n-1}. \]
Now we introduce the terminology dealing with extensions of symplectic $C_n$-morphisms:

\begin{defi}
Let $V$ be a free graded module of finite rank and let $\omega \in \drtf[Lie]{\clalg{V}}$ be a symplectic form of order zero. Let $m$ and $m'$ be two minimal symplectic $C_{n+1}$-structures on $V$ with the same underlying Frobenius algebra. We say that a pointed symplectic $C_n$-morphism $\phi$ from $\bar{m}$ to $\bar{m}'$ is extendable if there exists a \emph{symplectic} vector field $\gamma$ of order $n-1$ and degree zero such that $\exp(\gamma)\circ\phi$ is a symplectic $C_{n+1}$-morphism from $m$ to $m'$ and we call $\gamma$ an extension of $\phi$.

Let $\gamma$ and $\gamma'$ be two extensions of $\phi$. We say that $\gamma$ and $\gamma'$ are equivalent if there exists a \emph{symplectic} vector field $\eta$ of order $n-2$ and degree $-1$ such that
\[ \exp(\gamma)\circ\phi = \exp(\gamma')\circ\phi\circ\exp([m,\eta]) \mod(n). \]
The quotient of the set of all extensions of $\phi$ by this equivalence relation will be denoted by \mmext[S]{n}{\phi}{m}{m'}.
\end{defi}

Next we will define the appropriate obstruction to the extension of \emph{symplectic} $C_n$-morphisms. We shall use the same notation as in Definition \ref{def_morobs}. This abuse of notation will later be justified by Lemma \ref{lem_mobscm}.

\begin{defi}
Let $V$ be a free graded module of finite rank and let $\omega \in \drtf[Lie]{\clalg{V}}$ be a symplectic form of order zero. Let $m$ and $m'$ be two minimal symplectic $C_{n+1}$-structures on $V$ with the same underlying Frobenius algebra. Let $\phi$ be a pointed symplectic $C_n$-morphism from $\bar{m}$ to $\bar{m}'$. We define the 0-form $\obs(\phi)$ of order $n+1$ and degree $1+|\omega|$ by
\begin{equation} \label{eqn_smorob}
\obs(\phi):= \upiso[Lie]^{-1}\left(\phi \circ m \circ \phi^{-1} - m'\right) \mod (n+2).
\end{equation}
\end{defi}

We will now prove the symplectic analogues of theorems \ref{thm_morobs} and \ref{thm_morext}. We shall omit the proofs since they are the same as those of section \ref{sec_obstrc} verbatim except that we must use the fact from Lemma \ref{lem_duciso} that $\upiso[Lie]$ is an isomorphism.

\begin{theorem} \label{thm_smorob}
Let $A:=(V,m_2,\omega)$ be a strictly graded commutative Frobenius algebra and let $m$ and $m'$ be two minimal symplectic $C_{n+1}$-structures ($n \geq 3)$ on $V$ whose underlying Frobenius algebra is $A$, then \eqref{eqn_smorob} induces a map
\begin{displaymath}
\begin{array}{ccc}
\mmor[S]{n}{\bar{m}}{\bar{m}'} & \to & \hcaq[n+1,|\omega|]{A}, \\
\phi & \mapsto & \obs(\phi); \\
\end{array}
\end{displaymath}
which we will denote by $\obs_n$. The kernel of this map consists of precisely those symplectic $C_n$-morphisms which are extendable:
\[ \{ \phi \in \mmor[S]{n}{\bar{m}}{\bar{m}'} : \phi \textnormal{ is extendable} \} = \ker(\obs_n). \]
\end{theorem}
\noproof

\begin{theorem} \label{thm_smorex}
Let $A:=(V,m_2,\omega)$ be a strictly graded commutative Frobenius algebra and let $m$ and $m'$ be two minimal symplectic $C_{n+1}$-structures on $V$ whose underlying Frobenius algebra is $A$. Let $\phi \in \mmor[S]{n}{\bar{m}}{\bar{m}'}$ be an extendable symplectic $C_n$-morphism, then \hcaq[n,|\omega|-1]{A} acts freely and transitively on \mmext[S]{n}{\phi}{m}{m'}:
\begin{displaymath}
\begin{array}{ccc}
\hcaq[n,|\omega|-1]{A} \times \mmext[S]{n}{\phi}{m}{m'} & \to & \mmext[S]{n}{\phi}{m}{m'}, \\
(\alpha,\gamma) & \mapsto & \gamma + \upiso[Lie](\alpha). \\
\end{array}
\end{displaymath}
\end{theorem}
\noproof

\section{The Main Theorem: A Correspondence between $\ci$ and Symplectic \ci-structures} \label{correspondence}

In this section we will prove our main result; that a symplectic \ci-algebra $(V,m,\omega)$ is uniquely determined by its underlying \ci-algebra $(V,m)$ together with the structure of a commutative Frobenius algebra on $(V,m_2)$. The point is that the obstruction theories for symplectic and nonsymplectic $C_\infty$-algebras turn out to be `isomorphic' thanks to Corollary \ref{cor_caqiso}.

Once again, in this section we shall only consider symplectic \ci-algebras \emph{whose symplectic form $\omega$ has order zero}. It is important to note that whilst the results of sections \ref{sec_obstrc} and \ref{sec_symobs} apply equally well to {\li} and \ai-algebras, the results of this section will apply \emph{exclusively to \ci-algebras}.

It will be necessary to introduce the following definitions in order to formulate our main theorem later in this section:

\begin{defi}
\
\begin{enumerate}
\item[(i)]
Let $A:=(V,\mu_2)$ be a strictly graded commutative algebra, then the moduli space of minimal \ci-structures on $V$ fixing $\mu_2$ is denoted by \mstr{C}{\infty}{A} and defined as the quotient of the set
\[ \{ m:\clalg{V} \to \clalg{V} : m \text{ is a minimal \ci-structure and } m_2=\mu_2 \} \]
by the action under conjugation of the group $G$ consisting of all diffeomorphisms $\phi \in \Aut(\clalg{V})$ of the form
\[ \phi = \id + \phi_2 + \phi_3 + \ldots + \phi_k + \ldots, \]
where $\phi_i$ is an endomorphism of order $i$.
\item[(ii)]
Let $A:=(V,\mu_2,\omega)$ be a strictly graded commutative Frobenius algebra, then the moduli space of minimal \emph{symplectic} \ci-structures on $V$ fixing $\mu_2$ and $\omega$ is denoted by \mstr[S]{C}{\infty}{A} and defined as the quotient of the set
\[ \{ m:\clalg{V} \to \clalg{V} : m \text{ is a minimal \emph{symplectic} \ci-structure and } m_2=\mu_2 \} \]
by the action under conjugation of the group $G$ consisting of all \emph{symplectomorphisms} $\phi \in \Aut(\clalg{V})$ of the form
\[ \phi = \id + \phi_2 + \phi_3 + \ldots + \phi_k + \ldots, \]
where $\phi_i$ is an endomorphism of order $i$.
\end{enumerate}
\end{defi}

Recall that we call a diffeomorphism (symplectomorphism) of the form
\[ \phi = \id + \phi_2 + \phi_3 + \ldots + \phi_k + \ldots \]
a \emph{pointed} diffeomorphism (symplectomorphism).

\begin{defi}
\
\begin{enumerate}
\item[(i)]
Let $V$ be a free graded module and let $m$ and $m'$ be two minimal \ci-structures on $V$ with the same underlying algebra. We say that two \ci-morphisms $\phi$ and $\phi'$ from $m$ to $m'$ are homotopic if there exists a vector field $\eta$ of degree $-1$ such that
\[ \phi = \phi'\circ\exp([m,\eta]). \]

We denote the moduli space of \ci-morphisms from $m$ to $m'$ by \mmor{\infty}{m}{m'} and define it as the quotient of the set
\[ \{ \phi:\clalg{V} \to \clalg{V} : \phi \text{ is a \emph{pointed} \ci-morphism from $m$ to $m'$} \} \]
by the homotopy equivalence relation defined above.
\item[(ii)]
Let $V$ be a free graded module of finite rank and let $\omega \in \drtf[Lie]{\clalg{V}}$ be a symplectic form of order zero. Let $m$ and $m'$ be two minimal \emph{symplectic} \ci-structures on $V$ with the same underlying Frobenius algebra. We say that two \emph{symplectic} \ci-morphisms $\phi$ and $\phi'$ from $m$ to $m'$ are homotopic if there exists a \emph{symplectic} vector field $\eta$ of degree $-1$ such that
\[ \phi = \phi'\circ\exp([m,\eta]). \]

We denote the moduli space of \emph{symplectic} \ci-morphisms from $m$ to $m'$ by \mmor[S]{\infty}{m}{m'} and define it as the quotient of the set
\[ \{ \phi:\clalg{V} \to \clalg{V} : \phi \text{ is a \emph{pointed symplectic} \ci-morphism from $m$ to $m'$} \} \]
by the homotopy equivalence relation defined above.
\end{enumerate}
\end{defi}

Recall that in sections \ref{sec_obstrc} and \ref{sec_symobs} we defined moduli spaces \mstr{C}{n}{A} and \mstr[S]{C}{n}{A}. Given a strictly graded commutative Frobenius algebra $A$ we can define a map
\[ \iota : \mstr[S]{C}{n}{A} \to \mstr{C}{n}{A} \]
which is induced by the canonical inclusion of the Lie subalgebra of symplectic vector fields into the Lie algebra $\Der(\clalg{V})$ of all vector fields. This map is well defined because the group of symplectomorphisms is a subgroup of the group $\Aut(\clalg{V})$ of all diffeomorphisms. This map is also defined for $n=\infty$.

Similarly if $(m,\omega)$ is a symplectic $C_n$-structure on a free graded module $V$ of finite rank then we can define a map
\[ \iota : \mext[S]{n}{m} \to \mext{n}{m} \]
which is again induced by the canonical inclusion of the Lie subalgebra of symplectic vector fields into $\Der(\clalg{V})$.

The same applies to the moduli spaces \mmor{n}{m}{m'} and \mmor[S]{n}{m}{m'} which we also introduced in sections \ref{sec_obstrc} and \ref{sec_symobs}. Given a symplectic form $\omega \in \drtf[Lie]{\clalg{V}}$ of order zero and two symplectic $C_n$-structures $m$ and $m'$ with the same underlying Frobenius algebra we can define a map
\[ \iota : \mmor[S]{n}{m}{m'} \to \mmor{n}{m}{m'} \]
which is induced by the canonical inclusion of the subgroup of symplectomorphisms into the group $\Aut(\clalg{V})$ of all diffeomorphisms. The fact that this map is well defined modulo the homotopy equivalence relations follows tautologically from the fact that the symplectic vector fields are a Lie subalgebra of the Lie algebra of all vector fields. Likewise this map is also defined for $n=\infty$.

Given a symplectic form $\omega \in \drtf[Lie]{\clalg{V}}$ of order zero, two symplectic $C_n$-structures $m$ and $m'$ and a symplectic $C_n$-morphism $\phi$ from $m$ to $m'$ we can also define a map
\[ \iota : \mmext[S]{n}{\phi}{m}{m'} \to \mmext{n}{\phi}{m}{m'} \]
which is induced by the canonical inclusion of the Lie subalgebra of symplectic vector fields into $\Der(\clalg{V})$. Again it is tautological to check that this map is well defined modulo the homotopy equivalence relations.

If $A:=(V,m_2,\omega)$ is a strictly graded commutative Frobenius algebra then recall that there is a natural bigrading on cohomology which was defined in remarks \ref{rem_ordgde} and \ref{rem_ordgrd}. The isomorphisms of lemmas \ref{lem_duaiso} and \ref{lem_duciso} respect this bigrading. This gives us the following commutative diagram for all $j \in \mathbb{Z}$ and all $n\geq 1$:
\begin{equation} \label{fig_isocom}
\xymatrix{ \hcaq[n+1,j+|\omega|-2]{A} \ar^{I}[r] \ar^{\Psi}[rd] \ar_{\upiso[Lie]}[d] & \haq[n,j+|\omega|-2]{A}{A^*} \\ S\haq[nj]{A}{A} \ar@{^{(}->}[r] & \haq[nj]{A}{A} \ar^{\phiso[Lie]}[u] \\ }
\end{equation}
where \phiso[Lie] is the isomorphism which is defined by equation \eqref{eqn_ofmvec}, \upiso[Lie] is the isomorphism which is defined by equation \eqref{eqn_duciso} and $I$ is the map of Corollary \ref{cor_caqiso} which comes from the periodicity exact sequence. The map $\Psi$ is defined as $\Psi:=\phiso[Lie]^{-1}\circ I$. Any action of the group \haq[nj]{A}{A} could be pulled back along $\Psi$ to an action of the group \hcaq[n+1,j+|\omega|-2]{A}. We will now need the following auxiliary lemmas in order to prove our main result:

\begin{lemma} \label{lem_obscom}
Let $A:=(V,m_2,\omega)$ be a strictly graded commutative Frobenius algebra:
\begin{enumerate}
\item[(i)]
For all $n\geq 3$ the following diagram commutes:
\[ \xymatrix{ \mstr{C}{n}{A} \ar^-{\Obs_n}[r] & \haq[n+1,3]{A}{A} \\ \mstr[S]{C}{n}{A} \ar^{\iota}[u] \ar^-{\Obs_n}[r] & \hcaq[n+2,1+|\omega|]{A} \ar^{\Psi}[u] \\ } \]
\item[(ii)]
Let $m \in \mstr[S]{C}{n}{A}$ be an extendable symplectic $C_n$-structure. The map
\[ \iota : \mext[S]{n}{m} \to \mext{n}{m} \]
is \hcaq[n+1,|\omega|]{A}-equivariant.
\end{enumerate}
\end{lemma}

\begin{proof}
This is a tautological consequence of diagram \eqref{fig_isocom}.
\end{proof}

\begin{lemma} \label{lem_mobscm}
Let $A:=(V,m_2,\omega)$ be a strictly graded commutative Frobenius algebra and let $m$ and $m'$ be two minimal symplectic $C_{n+1}$-structures on $V$ whose underlying Frobenius algebra is $A$:
\item[(i)]
The following diagram commutes:
\[ \xymatrix{ \mmor{n}{\bar{m}}{\bar{m}'} \ar^{\obs_n}[r] & \haq[n,2]{A}{A} \\ \mmor[S]{n}{\bar{m}}{\bar{m}'} \ar^{\iota}[u] \ar^{\obs_n}[r] & \hcaq[n+1,|\omega|]{A} \ar^{\Psi}[u] \\ } \]
\item[(ii)]
Let $\phi \in \mmor[S]{n}{\bar{m}}{\bar{m}'}$ be an extendable symplectic $C_n$-morphism. The map
\[ \iota : \mmext[S]{n}{\phi}{m}{m'} \to \mmext{n}{\phi}{m}{m'} \]
is \hcaq[n,|\omega|-1]{A}-equivariant.
\end{lemma}

\begin{proof}
Again this follows tautologically from diagram \eqref{fig_isocom}.
\end{proof}

We are now ready to formulate our main result:

\begin{theorem} \label{thm_main}
Let $A:=(V,m_2,\omega)$ be a strictly graded \emph{unital} commutative Frobenius algebra:
\begin{enumerate}
\item[(i)]
The map
\[ \iota : \mstr[S]{C}{\infty}{A} \to \mstr{C}{\infty}{A} \]
is a bijection.
\item[(ii)]
Let $m$ and $m'$ be two minimal symplectic \ci-structures on $V$ whose underlying Frobenius algebra is $A$. The map
\[ \iota : \mmor[S]{\infty}{m}{m'} \to \mmor{\infty}{m}{m'} \]
is a surjection.
\end{enumerate}
\end{theorem}

\begin{proof}
Let us begin by proving that the map $\iota : \mstr[S]{C}{\infty}{A} \to \mstr{C}{\infty}{A}$ is surjective. Let \[ m = m_2 + m_3 + \ldots + m_n + \ldots \]
be a minimal \ci-structure on $V$. We will inductively construct a sequence of \emph{symplectic} vector fields $m'_i, \ 3\leq i < \infty$ and a sequence of vector fields $\gamma_i , \ 2\leq i < \infty$, where $m'_i$ has order $i$ and degree one and $\gamma_i$ has order $i$ and degree zero, such that
\begin{enumerate}
\item[(i)]
\[ m':= m_2 + m'_3 + \ldots + m'_n + \ldots \]
is a minimal \emph{symplectic} \ci-structure.
\item[(ii)]
\[ \phi:=\ldots\circ\exp(\gamma_n)\circ\ldots\circ\exp(\gamma_3)\circ\exp(\gamma_2) \]
is a \ci-morphism from $m$ to $m'$.
\end{enumerate}

Let us assume that we have constructed a sequence of \emph{symplectic} vector fields $m'_3,\ldots,m'_{n-1}$ of degree one and a sequence of vector fields $\gamma_2,\ldots,\gamma_{n-2}$ of degree zero, where $m'_i$ and $\gamma_i$ have order $i$, satisfying
\begin{enumerate}
\item[(i)]
\[ m':= m_2 + m'_3 + \ldots + m'_{n-1} \]
is a minimal \emph{symplectic} $C_n$-structure.
\item[(ii)]
\[ \phi:= \exp(\gamma_{n-2})\circ\ldots\circ\exp(\gamma_2) \]
is a minimal $C_n$-morphism from $\bar{m}:= m_2 + \ldots + m_{n-1}$ to $m'$.
\end{enumerate}
Note that the base case $n=3$ is trivial. The $C_n$-structures $\bar{m}$ and $m'$ represent the same class in \mstr{C}{n}{A}, therefore by Lemma \ref{lem_obscom} and Theorem \ref{thm_obstrc} we see that
\[ \Psi(\Obs_n(m')) = \Obs_n(\iota(m')) = \Obs_n(\bar{m}) = 0 \]
because the $C_n$-structure $\bar{m}$ is extendable. Corollary \ref{cor_caqiso} now implies that $\Obs_n(m') = 0$ and therefore by Theorem \ref{thm_symobs}, $m'$ is an extendable symplectic $C_n$-structure (recall from diagram \eqref{fig_isocom} that $\Psi$ was defined in terms of the map $I$ of Corollary \ref{cor_caqiso}).

Consider the \ci-structure $\phi\circ m\circ\phi^{-1}$; there exists a vector field $m_n$ (not necessarily symplectic) of order $n$ and degree one such that
\[ \phi\circ m\circ\phi^{-1} = m' + m_n \mod(n+1) \]
and therefore $m'+m_n$ is a $C_{n+1}$-structure on $V$. We also know from the above argument that there exists a \emph{symplectic} vector field $m'_n$ of order $n$ and degree one such that $m'+m'_n$ is a \emph{symplectic} $C_{n+1}$-structure. Let us denote the corresponding class of $m'_n$ in \mext[S]{n}{m'} by $\tilde{m}'_n$ and the corresponding class of $m_n$ in \mext{n}{m'} by $\tilde{m}_n$. By Theorem \ref{thm_extend} there exists a cohomology class $\xi_n \in \haq[n,2]{A}{A}$ such that
\[ \iota(\tilde{m}'_n) + \xi_n = \tilde{m}_n. \]
By Corollary \ref{cor_caqiso} and Lemma \ref{lem_obscom} there exists a cohomology class $\alpha \in \hcaq[n+1,|\omega|]{A}$ such that
\[ \iota(\tilde{m}'_n + \Upsilon(\alpha)) = \tilde{m}_n. \]
We see that by modifying our choice of symplectic vector field $m'_n$ appropriately we can assume that there exists a vector field $\gamma_{n-1}$ of order $n-1$ and degree zero such that
\[ \exp(\gamma_{n-1})\circ(m'+m_n)\circ\exp(-\gamma_{n-1}) = m' + m'_n \mod(n+1). \]
This completes the inductive step and proves that the map $\iota : \mstr[S]{C}{\infty}{A} \to \mstr{C}{\infty}{A}$ is surjective.

The proof that the map $\iota : \mmor[S]{\infty}{m}{m'} \to \mmor{\infty}{m}{m'}$ is surjective proceeds in a similar fashion to the above proof. We shall outline the main ideas of the proof and leave the reader to provide the details:
\begin{enumerate}
\item[(i)]
Given a \ci-morphism $\phi$ from $m$ to $m'$ the idea is to construct by induction a sequence of \emph{symplectic} vector fields $\gamma'_i, \ 2 \leq i < \infty$ of degree zero and a sequence of vector fields $\eta_i, \ 1 \leq i < \infty$ of degree $-1$, where $\gamma'_i$ and $\eta_i$ have order $i$, such that
\begin{enumerate}
\item
\[ \phi':=\ldots\circ\exp(\gamma'_n)\circ\ldots\circ\exp(\gamma'_3)\circ\exp(\gamma'_2) \]
is a \emph{symplectic} \ci-morphism from $m$ to $m'$.
\item
\[ \phi=\phi'\circ\left[\ldots\circ\exp([m,\eta_n])\circ\ldots\circ\exp([m,\eta_2])\circ\exp([m,\eta_1])\right]. \]
\end{enumerate}
\item[(ii)]
At the $n$th stage use Lemma \ref{lem_mobscm}, Corollary \ref{cor_caqiso} and theorems \ref{thm_morobs} and \ref{thm_smorob} to show that the obstruction to extending the given symplectic $C_n$-morphism to the next level vanishes.
\item[(iii)]
Use Lemma \ref{lem_mobscm}, Theorem \ref{thm_morext} and \emph{both} parts (ii) \emph{and} (iii) of Corollary \ref{cor_caqiso} to modify this extended symplectic $C_{n+1}$-morphism to one which is homotopy equivalent to $\phi$ modulo $(n)$.
\end{enumerate}

Finally we show that the map $\iota : \mstr[S]{C}{\infty}{A} \to \mstr{C}{\infty}{A}$ is injective. Let $m,m'\in\mstr[S]{C}{\infty}{A}$ be two symplectic \ci-structures and suppose that there exists a pointed \ci-morphism $\phi$ from $m$ to $m'$, then by part (ii) of this theorem there exists a \emph{symplectic} pointed \ci-morphism $\phi'$ from $m$ to $m'$ which is homotopy equivalent to $\phi$.
\end{proof}

\begin{rem}
A natural question is whether a result similar to Theorem \ref{thm_main} holds in the $A_\infty$ and $L_\infty$ contexts. The answer is no. The crucial point on which the proof of Theorem \ref{thm_main} turns is that the cyclic Harrison theory essentially coincides with the (noncyclic) Harrison theory, cf. Corollary \ref{cor_caqiso}. This fails badly for both the associative and Lie cases. For example, let $g$ be a semisimple Lie algebra which could be considered as a symplectic \li-algebra together with its Killing form. Then $\hlac[i]{g}{g}=0$ unless $i=0$; however the cyclic theory $H^\bullet(g,\gf)$ is not zero in dimensions $\geq 3$.

In the associative case a relevant counterexample is provided by the group ring of a finite nonabelian group $G$. Indeed, in this case \hchoch{\gf[G]} is isomorphic to the direct sum of copies of $\hchoch{\gf}\cong \gf[u]$ where the summation is over all conjugacy classes of $G$. On the other hand $\hhoch[i]{\gf[G]}{\gf[G]}=0$ for all $i>0$.
\end{rem}

\begin{rem}
There is a construction, cf. \cite{kontfd} which associates to any minimal symplectic $C_\infty$-algebra a cycle in an appropriate version of the graph complex. Moreover, two weakly equivalent symplectic $C_\infty$-algebras give rise to homologous cycles. Therefore, our result shows that the corresponding graph homology class only depends on the underlying $C_\infty$-algebra together with the structure of a graded commutative Frobenius algebra on its underlying module. In particular, a graph homology class could be associated to any Poincar\'e duality space. It would be interesting to express this construction in classical homotopy theoretic terms, i.e. Massey products.
\end{rem}

\section{Symplectic \ai-algebras and string topology} \label{strings}

In this section we apply our  results to establish the existence of certain structures on the ordinary and equivariant homology of the free loop space on a manifold or, more generally a formal Poincar\'e duality space. These structures; namely the loop product, the loop bracket and the string bracket, were introduced and studied in \cite{chasul} under the general heading `string topology'. Since the original work of Sullivan and Chas a number of other papers appeared which investigated string topology structures from different perspectives. In particular, the approach of \cite{cohenj} and  \cite{jklein} was purely homotopy theoretic while that of \cite{fethvi} was based on the theory of minimal models.

In order to be able to directly apply Theorem \ref{thm_main} and to keep this paper within reasonable limits we will assume that our ground ring $\gf$ is the field of rational numbers $\mathbb{Q}$. However it seems likely that a more careful analysis should yield similar results in positive characteristic as well. The `string topology' operations that we define resemble formally those of Sullivan and Chas and we expect that they in fact agree although this issue is not considered here.

For a topological space $X$ (always assumed to be a $CW$-complex) we denote by $C^\bullet(X)$ the singular cochain algebra of $X$ with coefficients in $\mathbb{Q}$. The homology and cohomology of $X$ with coefficients in $\mathbb{Q}$ will be denoted by $H_\bullet(X)$ and $H^\bullet(X)$ respectively.   The space of all maps $S^1\rightarrow X$ from the circle $S^1$ to $X$ will be denoted by $LX$. The space $LX$ has an action of $S^1$ and the corresponding equivariant homology will be denoted by $H^{S^1}_\bullet(LM)$:
\[  H^{S^1}_\bullet(LM):= H^\bullet(ES^1\times _{S^1} LM).\]
\begin{theorem} \label{loopproduct}
Let $M$ be a simply-connected  Poincar\'e duality space of formal dimension $d$. Then the graded vector space $\mathbb{H}_\bullet(M):=H_{\bullet+d}(LM)$ has the structure of a graded Gerstenhaber algebra. Two homotopy equivalent Poincar\'e duality spaces give rise to isomorphic Gerstenhaber algebras.
\end{theorem}

\begin{rem}
The graded Gerstenhaber algebra structure on $H_{\bullet+d}(LM)$ consists of two operations: a graded commutative and associative product
\[a\otimes b\mapsto a\bullet b:H_p(LM)\otimes H_q(LM)\rightarrow H_{p+q-d}(LM)\]
and the graded Lie bracket
\[a\otimes b\mapsto \{a, b\}:H_p(LM)\otimes H_q(LM)\rightarrow H_{p+q-d+1}(LM).\]
Moreover, $\{a,?\}$ is a graded derivation of $\bullet$ for any $a\in H_\bullet(LM)$. The operations $\bullet$ and $\{,\}$ are called the loop product and the loop bracket respectively.
\end{rem}

\begin{proof}
There exists a minimal $C_\infty$-algebra $(H^\bullet(M),m)$ weakly equivalent to the cochain algebra $C^\bullet(M)$. Moreover, by Theorem \ref{thm_main} we could assume that the Poincar\'e duality form extends to an invariant inner product of order $d$ on $(H^\bullet(M),m)$. To ease notation we will denote $H^\bullet(M)$ by $V$. Then \hhoch{V}{V} will denote the Hochschild cohomology of $(H^\bullet(M),m)$ with coefficients in itself and \hhoch{V}{V^*} the Hochschild cohomology of $(H^\bullet(M),m)$ with coefficients in $V^*$. 

Then by Lemma \ref{lem_duaiso} we have an isomorphism $\hhoch{V}{V^*} \cong \hhoch[\bullet+d]{V}{V}$. It is well known that the Hochschild cohomology of a differential graded algebra with coefficients in itself has the structure of a graded Gerstenhaber algebra; moreover two weakly equivalent DGAs give rise to isomorphic Gerstenhaber algebras, cf. \cite{femeth}. Since there exists a DGA weakly equivalent to the \ai-algebra $V$ we conclude that $\hhoch{V}{V} \cong \hhoch[\bullet-d]{V}{V^*}$ supports the structure of a graded Gerstenhaber algebra.

Since $C^\bullet(M)$ is weakly equivalent to $(V,m)=(H^\bullet(M),m)$ we have an isomorphism of graded vector spaces $\hhoch{C^\bullet(M)}{[C^\bullet(M)]^*} \cong \hhoch{V}{V^*}$. Therefore  \hhoch[\bullet-d]{C^\bullet(M)}{[C^\bullet(M)]^*} has the structure of a graded Gerstenhaber algebra. 

Furthermore, the graded space \hhoch{C^\bullet(M)}{[C^\bullet(M)]^*} is isomorphic to the $\gf$-vector dual to \Hhoch{C^\bullet(M)}{C^\bullet(M)} (Hochschild homology of $C^\bullet(M)$) and the latter is isomorphic to $H^\bullet(LM)$ by \cite{burfi}, \cite{goodwi}. 
Since the space $LM$ is of finite type we have  $[H^\bullet(LM)]^*\cong H_\bullet(LM)$ and we conclude that  \[\mathbb{H}_\bullet(M)=H_{\bullet+d}(LM)\cong \hhoch[\bullet-d]{C^\bullet(M)}{[C^\bullet(M)]^*}\]
has the structure of a graded Gerstenhaber algebra as claimed.

Homotopy invariance is likewise clear since for two homotopy equivalent spaces $M$ and $N$ the cochain algebras $C^\bullet(M)$ and $C^\bullet(N)$ are weakly equivalent.
\end{proof}

\begin{rem}
Although two homotopy equivalent Poincar\'e duality spaces $M, M^\prime$ give rise to isomorphic Gerstenhaber algebras on their string homology this structure is not natural in the sense that a map $M\rightarrow M^\prime$ does not lead to a map between the corresponding Gerstenhaber algebras. An analogous situation arises when one consider the monoid of self-maps $\Map(X,X)$ of a topological space $X$; the association $X\mapsto \Map(X,X)$ is not a functor even though homotopy equivalent spaces give rise to homotopy equivalent monoids of self-maps. The same remark applies to the string bracket considered below.
\end{rem}

The other part of the string topology operations is called the \emph{string bracket} and is defined on the equivariant homology of $LM$: $H_\bullet^{S^1}(LM):=LM\times_{S^1}ES^{1}$. To put the string bracket in the proper context we need to introduce \emph{negative} cyclic cohomology of $A_\infty$-algebras.

Let $(V,M)$ be an $A_\infty$-algebra, not necessarily unital. 

\begin{defi}
The negative cyclic complex $CC^\bullet_-(V)$ of $(V,m)$ is the total complex of the following bicomplex formed by taking direct sums:
\begin{equation}\label{negativecyclic}
\xymatrix{\ldots \ar^-{1-z}[r] & \cbr{V} \ar^-N[r] & \choch{V}{V^*} \ar^-{1-z}[r] & \cbr{V} }.
\end{equation}
The cohomology of $CC^\bullet_-(V)$ will be denoted by $HC^\bullet_-(V)$.
\end{defi}

\begin{rem}
Standard spectral sequence arguments show that two weakly equivalent \ai-algebras have isomorphic negative cyclic cohomology.
\end{rem}

\begin{lemma}\label{comparison}
Let $(V,m)$ be a unital $A_\infty$-algebra for which there exists an integer $N$ such that $\hhoch[k]{V}{V^*}=0$ for $k>N$. Then for any integer $n$ we have
\[HC^n_-(V)\cong \hchoch[n+1]{V}.\]
\end{lemma}

\begin{proof}
Note that the even-numbered columns isomorphic to \cbr{V} are acyclic since the $A_\infty$-algebra $(V,m)$ is unital. This together with
the assumption that $\hhoch[k]{V}{V^*}=0$ for sufficiently large $k$ ensures that $CC^\bullet_-(V)$ is quasi-isomorphic to the subcomplex $\overline{CC}^\bullet_-(V)$ whose columns are appropriate truncations of the columns of $CC^\bullet_-(V)$  and which has no nonzero terms above a certain horizontal line.  

Consider the auxiliary complex $\widetilde{CC}^\bullet_-(V)$ formed by the \emph{direct product} totalisation of \eqref{negativecyclic}
Then the complex  $\overline{CC}^\bullet_-(V)$ can also be considered as a subcomplex of $\widetilde{CC}^\bullet_-(V)$. Note that for the complex $\overline{CC}^\bullet_-(V)$ there is no difference between the direct product or direct sum totalisation and therefore both its spectral sequences converge. 

Comparing the appropriate spectral sequences for $\overline{CC}^\bullet_-(V)$ and $\widetilde{CC}^\bullet_-(V)$ we see, using the exactness of  the rows, that these complexes are quasi-isomorphic and therefore they are both quasi-isomorphic to $CC^\bullet_-(V)$.
That shows that under our assumptions the direct product and direct sum totalisations of \eqref{negativecyclic} have the same cohomology.

Next consider the following bicomplex $\widehat{CC}^\bullet_-(V)$:
\[\xymatrix{ \ldots \ar^-{N}[r] & \choch{V}{V^*} \ar^-{1-z}[r] & \cbr{V} \ar^-{N}[r] & \choch{V}{V^*} }.\]
Its cohomology will be denoted by $\widehat{HC}^\bullet_-(V)$.
Similarly to the above, it makes no difference which totalisation (direct sum or direct product) to take. 
Now zig-zag arguments (or the appropriate spectral sequence)  show, using the exactness of the rows,  that  $\widehat{HC}^\bullet_-(V)$ is isomorphic to the complex formed by the image of the horizontal differential in the rightmost column. The latter complex is clearly isomorphic to \cchoch{V}. Therefore   $\widehat{HC}^\bullet_-(V)\cong HC^\bullet_-(V)$.                                                               

Finally, using the acyclicity of \cbr{V} we see that the projection map
\[CC^\bullet_-(V)\longrightarrow \Sigma\widehat{CC}^\bullet_-(V)\]
is a quasi-isomorphism.  This finishes the proof.
\end{proof}

Next we shall make a link between negative cyclic cohomology of an $A_\infty$-algebra and a more customary notion of the negative cyclic \emph{homology} of a differential graded algebra. This departure from our convention to work with cohomology rather than homology is necessary because the equivariant cohomology of a loop space is expressed in \cite{jjones} in terms of cyclic homology rather than cyclic cohomology.

Recall that for a DGA $V$ its cyclic Tsygan complex is the following bicomplex $CC_\bullet^-(V)$ lying in the right half-plane and formed by taking direct products:
\[\xymatrix{ \ldots & \Cbr{V} \ar_-{1-z}[l] & \Choch{V}{V} \ar_-N[l] & \Cbr{V} \ar_-{1-z}[l] }.\]
Here \Choch{V}{V} is the usual homological Hochschild complex of $V$ and \Cbr{V} is the bar-complex of $V$ (acyclic in the unital case). The operators $1-z$ and $N$ are formed as in the cohomological cyclic complex.

Then we have the following almost obvious result.

\begin{lemma}\label{finiteness}
Let $V$ be a (unital) differential graded algebra whose Hochschild homology \Hhoch{V}{V} is finite dimensional in each degree. Then $HC^-_\bullet(V)$ is isomorphic to the graded dual of the graded $\gf$-vector space $HC^\bullet_-(V)$.
\end{lemma}

\begin{proof}
Note that the finite dimensionality assumption ensures that the complex $[\choch{V}{V^*}]^*$ is quasi-isomorphic to \Choch{V}{V}. Since the functor $?^*:=\Hom_\gf(?,\gf)$ takes direct sums to direct products we conclude that $[CC^\bullet_-(V)]^*$ is quasi-isomorphic to $CC^-_\bullet(V)$. 
\end{proof}

\begin{rem}\label{simpleconnectivity}
Let $V$ be a differential graded algebra such that $H^n(V)=0$ for $n<0$, $H^0(V)=\gf$, $H^1(V)=0$ and $H^n(V)$ is finite dimensional for all $n$. For example, the cochain algebra $C^\bullet(X)$ of a simply-connected space $X$ of finite type satisfies these conditions. Then it is easy to see from the spectral sequence associated with the normalised Hochschild complex $\overline{C}^\bullet_\mathrm{Hoch}(V,V^*)$ that \hhoch{V}{V^*} is finite dimensional in each degree.
\end{rem}

With these preparations we are ready to formulate our last result.

\begin{theorem}
Let $M$ be a rational Poincar\'e duality space of formal dimension $d$. Then $H_\bullet^{S^1}(LM)$ has the structure of a graded Lie algebra of degree $2-d$. If $d\neq 0 \mod 4$ then two homotopy equivalent spaces give rise to isomorphic graded Lie algebras. If $d=0\mod 4$ then the choice of the fundamental cycle in $H_d(M)$ leads to two possibly nonisomorphic  graded Lie algebra structures on  $H_\bullet^{S^1}(LM)$.
\end{theorem} 

\begin{proof}
We know from \cite{jjones} that $H^n_{S^1}(LM)\cong HC_{-n}^-(C^\bullet(M))$ and using Lemma \ref{finiteness} together with Remark \ref{simpleconnectivity} we conclude that $H_n^{S^1}(LM)\cong HC^{-n}_-(C^\bullet(M))$.

Furthermore, the differential graded algebra $C^\bullet(M)$ clearly satisfies the conditions of Lemma \ref{comparison} with $N=0$ since $\hhoch{C^\bullet(M)}{[C^\bullet(M)]^*} \cong H_{-\bullet}(LM)$ and $H_{\bullet}(LM)$ is concentrated in nonnegative degrees. Therefore we have an isomorphism
\begin{equation} \label{shift}
H_n^{S^1}(M)\cong HC^{-n}_-(C^\bullet(M))\cong HC^{-n+1}(C^\bullet(M)).
\end{equation}
Similar to the proof of Theorem \ref{loopproduct} let $(V,m)$ be a minimal symplectic $A_\infty$-algebra weakly equivalent to $C^\bullet(M)$. Since by Lemma \ref{lem_duciso} the complex \cchoch{V} is isomorphic (with an appropriate shift) to the complex $S\choch{V}{V}$ consisting of symplectic vector fields its homology supports the structure of a graded Lie algebra (the Gerstenhaber bracket). Since $\hchoch[n]{C^\bullet(M)} \cong \hchoch[n]{V}$ we have the Lie bracket on \hchoch[n]{C^\bullet(M)}:
\[\hchoch[n]{C^\bullet(M)} \otimes \hchoch[k]{C^\bullet(M)} \rightarrow \hchoch[n+k-d+1]{C^\bullet(M)}.\]
The latter is translated into the string bracket via isomorphism \eqref{shift}:
\[H^{S^1}_n(LM)\otimes H^{S^1}_k(LM)\rightarrow H^{S^1}_{n+k-d+2}(LM).\]

Now let us consider the problem of homotopy invariance.  For a given rational Poincar\'e duality space $N$ homotopy equivalent to $M$, the cohomology algebras $H^\bullet(M)$ and $H^\bullet(N)$ are isomorphic but the matrices of inner products on them differ by a nonzero rational number coming from the choice of the fundamental class. If $d\neq 0\mod 4$ then we see, using appropriate rescalings that $H^\bullet(M)$ and $H^\bullet(N)$ are isomorphic as Frobenius algebras, after all. It follows that in this case the graded Lie algebra structures on  $H^{S^1}_\bullet(LM)$ and $H^{S^1}_\bullet(LN)$ are isomorphic. If  $d=0\mod 4$ then the signatures of $M$ and $N$ might differ by a sign depending on whether the given homotopy equivalence $M\rightarrow N$ preserves orientation or changes it.  So we have precisely two nonisomorphic structures  of a symplectic $A_\infty$-algebra on $H^\bullet(M)$ depending on the choice of orientation. This leads to two possibly nonisomorphic graded Lie algebras on $H^{S^1}_\bullet(LM)$.
\end{proof}

\begin{rem}
Clearly, the loop bracket and the loop product as defined above are nontrivial as a quick calculation of the Hochschild cohomology of $H^\bullet(S^n)$ makes clear. Since the cyclic Hochschild cohomology of a truncated polynomial algebra is concentrated in even dimensions, we conclude that for spaces such as $S^n$ and $\mathbb{C}P^n$ the string bracket is zero. This is in agreement with the calculations made in \cite{fethvi}. To obtain an example of a nontrivial string bracket consider the space $X:=S^3\times S^3$ and let $A:=H^\bullet(X)$. Then the Lie algebra consisting of symplectic vector fields of degree $0$ is identified with the Lie algebra of the group of automorphisms of $H^3(X)$ preserving the Poincar\'e duality form. In other words this is the algebra of traceless $2\times 2$ matrices and it clearly has nontrivial commutators, therefore the string bracket is nontrivial in $H_\bullet^{S^1}(X)$.
\end{rem}

\appendix
\section{Formal $\gf$-algebras} \label{app_todual}

We felt it necessary to include a section dealing with our formal passage to the dual language of formal \gf-algebras that we use in our paper since it is difficult to find any references for the material in the literature. 

We will denote the category of free graded $\gf$-modules by $Mod_\gf$.  In what follows we shall omit the adjective `graded' when talking about graded $\gf$-algebras and graded $\gf$-modules since the ungraded ones are not considered.

\begin{defi} \label{def_promod}
A profinite $\gf$-module $V$ is a $\gf$-module which is an inverse limit of a diagram of free $\gf$-modules of finite rank; $V=\inlim{\alpha}{V_\alpha}$. A fundamental system of neighbourhoods of zero of $V$ is generated by the kernels of the projections $V\rightarrow V_\alpha$. The induced topology is known as the inverse limit topology. Profinite $\gf$-modules form a category {\pvect} in which the morphisms are \emph{continuous} $\gf$-linear maps between $\gf$-modules.
\end{defi}

\begin{rem}
A profinite $\gf$-module can also be defined as a \emph{topologically free} $\gf$-module, i.e. a module $M$ for which there exist a collection  of elements $\{t_i\}_{i\in I}\subset M$ (topological basis) such that any element could be uniquely represented as a (possibly uncountably infinite) linear combination of the $t_i$'s.
\end{rem}

\begin{rem}
When we talk about a submodule of a profinite module generated by a set $X$ we mean the module generated by all convergent infinite linear combinations of elements in the set $X$. Note also that a submodule as well as a quotient of a profinite $\gf$-module need not be a profinite $\gf$-module in general. However in all cases that we encounter the $\gf$-modules and their submodules are obtained from profinite $\mathbb{Q}$-vector spaces by extending the scalars to $\gf$ and so this complication never arises.
\end{rem} 

We would like to define the appropriate notion of tensor product in the category \pvect. This is given by the following definition:

\begin{defi}
Let $V=\inlim{\alpha}{V_\alpha}$ and $U=\inlim{\beta}{U_\beta}$ be profinite $\gf$-modules, then their completed tensor product $V \cotimes U$ is given by the formula,
\[ V \cotimes U:= \inlim{\alpha ,\beta}{V_\alpha \otimes V_\beta}. \]
\end{defi}

\begin{rem}
The completed tensor product could also be introduced as the universal object solving the problem of factorising continuous bilinear forms. Given two continuous linear maps between profinite modules $\phi:V \to W$ and $\psi:U \to X$ we could form the continuous linear map $\phi \cotimes \psi:V\cotimes U \to W \cotimes X$ in an obvious way. With this definition the category {\pvect} becomes a \emph{symmetric monoidal} category. If we endow {\pvect} with the direct product $\Pi$ then it also becomes an additive category.

Note that the construction of $V \cotimes U$ is not canonical as it depends on a preferred system $V_\alpha$, $U_\beta$ or rather a preferred choice of topological basis. Obviously different choices will produce continuously linearly isomorphic profinite modules.
\end{rem}

\begin{prop} \label{prop_antieq}
The functor $F:Mod_\gf \to \pvect$ given by sending $V$ to its linear dual $F(V):=\Hom_\gf(V,\gf)$ establishes an anti-equivalence of additive symmetric monoidal categories whose inverse functor $G$ is given by sending $U$ to the continuous linear dual $G(U):=\Hom_\mathrm{cont}(U,\gf)$.
\end{prop}
\noproof

\begin{rem}
Since a free $\gf$-module is a direct limit (union) of its finite rank free submodules, it is easy to see that the anti-equivalence of proposition \ref{prop_antieq} identifies the tensor product in $Mod_\gf$ with the completed tensor product in \pvect.
\end{rem}

We shall now discuss the notions of \emph{formal} associative, commutative and Lie algebras. To treat them uniformly we consider \emph{nonunital} associative and commutative $\gf$-algebras. When we say `$\gf$-algebra', that will simply mean that the corresponding statement could be applied to either commutative, associative or Lie algebras.

\begin{defi}Let $A$ be a $\gf$-algebra whose underlying $\gf$-module is free of finite rank. For technical convenience we assume that $A$ is obtained from $\mathbb{Q}$ (or any subfield of $\gf$) by the extension of scalars. We say that $A$ is \emph{nilpotent} if there exists a positive integer $N$ for which the product (or bracket in the Lie case) of any $N$ elements in $A$ is zero.
An inverse limit of nilpotent $\gf$-algebras is called a formal $\gf$-algebra.    
A morphism between two formal $\gf$-algebras is simply a continuous homomorphism of corresponding structures. The categories of formal associative, commutative and Lie algebras $\gf$-algebras will be denoted by $\mathcal{F}_{Ass}Alg, \mathcal{F}_{Com}Alg, \mathcal{F}_{Lie}Alg$ respectively. The notation $\mathcal{F}Alg$ will be used to denote either of the three categories.  A formal $\gf$-algebra supplied with a \emph{continuous} differential will be called a formal differential graded $\gf$-algebra.
\end{defi}

\begin{rem}
Every nonunital $\gf$-algebra (commutative or associative) gives rise to a unital one obtained by the well known procedure of adjoining a unit. We will call a formal associative or commutative algebra with an adjoined unit a \emph{formal augmented $\gf$-algebra}. In the main text by a `formal commutative or associative algebra' we will always mean a `formal commutative or associative augmented $\gf$-algebra'. 
\end{rem}

An important example of a formal $\gf$-algebra is the so-called \emph{pro-free} $\gf$-algebra. To start, we define the \emph{free} associative $\gf$-algebra on a free $\gf$-module $V$ as $TV_+:=\bigoplus_{i=1}^\infty V^{\otimes i}$. Similarly the free commutative algebra on $V$ is defined as $SV_+:=\bigoplus_{i=1}^\infty (V^{\otimes i})_{S_i}$. Here we used the subscript $+$ to avoid confusion with the unital free associative algebra $TV:=\bigoplus_{i=0}^\infty V^{\otimes i}$ and similarly in the commutative case. The free Lie algebra on $V$ could be defined as the submodule in $T_+(V)$ spanned by all Lie monomials in $V$.

\begin{defi}
Let $V$ be a free $\gf$-module of finite rank.  Then the pro-free (associative, commutative or Lie) $\gf$-algebra on $V$ is the $\gf$-algebra formed by formal power series (associative, commutative or Lie) in elements of $V$. It will be denoted  by $\widehat{T}_+(V), \widehat{S}_+(V), \widehat{L}(V)$ in the associative. commutative and Lie cases respectively.

If $V=\inlim{\alpha}{V_\alpha}$ is a profinite $\gf$-module then we define $\widehat{T}_+(V)$ as $\inlim{\alpha} {\widehat{T}_+(V)}$ and similarly in the commutative and Lie cases. We will denote by $\widehat{T}(V)$ and $ \widehat{S}(V)$ the unital versions of $\widehat{T}_+(V),$ and $\widehat{S}_+(V)$ respectively.
 \end{defi}

Clearly pro-free $\gf$-algebras are formal $\gf$-algebras. Furthermore we have the following result whose proof is a simple check.

\begin{prop}
 The functor $V\mapsto F(V)$ from  $\pvect$ to $\mathcal{F}Alg$ is left adjoint to the forgetful functor.
\end{prop}\noproof

\begin{rem}
The category of formal $\gf$-algebras is equivalent to the category of \emph{cocomplete} $\gf$-coalgebras,
i.e. $\gf$-coalgebras $C$ for which the kernels of iterated coproducts $\Delta^n:C\rightarrow C^{\otimes n}$ form an exhaustive filtration.  The functor from cocomplete $\gf$-coalgebras to formal $\gf$-algebras is simply taking the $\gf$-linear dual, and the 
inverse functor is taking the \emph{continuous} dual. Because of this equivalence the theory of $\infty$-algebras is often formulated in terms of coalgebras and coderivations.
\end{rem}

We will often consider vector fields (=continuous derivations) and diffeomorphisms (=continuous invertible homomorphisms) of formal $\gf$-algebras. The following proposition is straightforward:

\begin{prop} \label{prop_dergen}
Let $\widehat{T}(V)$ ($\widehat{S}(V)$, $\widehat{L}(V)$) be a pro-free associative (commutative or Lie) algebra on a profinite $\gf$-module $V$.
Then any vector field or diffeomorphism of $\widehat{T}(V)$ ($\widehat{S}(V)$, $\widehat{L}(V)$) is uniquely determined by its restriction on $V$. In particular, any vector field $\xi$ has the form
\[\xi=\sum f_i(\bf t)\partial_{t_i},\]
where $\boldsymbol{t}:=\{t_i\}_{i \in I}$ is a topological basis of $V$ and $f_i(\boldsymbol{t})$ is a formal power series (associative, commutative or Lie) in the $t_i$'s.  
\end{prop}\noproof

We will now define the formal universal enveloping algebra  of a formal Lie algebra $\mathcal{G}$.  

\begin{defi}
\
\begin{enumerate}
\item
First let $\mathcal{G}$ be a nilpotent Lie $\gf$-algebra whose underlying $\gf$-module is free of finite rank. Denote by ${\mathcal{U}}\mathcal{G}$ its usual universal enveloping algebra and by $I(\mathcal{G})$ the augmentation ideal in ${\mathcal{U}}\mathcal{G}$. Then  $\overline{\mathcal{U}}\mathcal{G}:=\inlim{n}I(\mathcal{G})/I^n(\mathcal{G})$.
\item 
Now let ${\mathcal{G}}=\inlim{} {\mathcal{G}}_i$ be a formal Lie $\gf$-algebra. Here $\{{\mathcal{G}}_i\}_{i\in i}$ is an inverse system of nilpotent Lie $\gf$-algebras. Then $\overline{\mathcal{U}}\mathcal{G} := \inlim{n}\overline{\mathcal{U}}\mathcal{G}_i$.
\end{enumerate}
\end{defi}

\begin{rem}
Note that  $\overline{\mathcal{U}}\mathcal{G}$ is a formal associative algebra, in particular it is nonunital. We will denote by $\widehat{\mathcal{U}}\mathcal{G}$ the unital $\gf$-algebra obtained from $\overline{\mathcal{U}}\mathcal{G}$ by adjoining a unit.
\end{rem}

Just like the usual universal enveloping algebra the formal one is characterised by a certain universal property. Namely, associated with any formal associative algebra $A$ is a formal Lie algebra $l(A)$ whose underlying $\gf$-module coincides with that of $A$ and the Lie bracket is defined as
\[ [a,b]:=ab-(-1)^{|a||b|}ba; \quad a,b\in A. \]
Thus, we have a functor $l$ from formal associative $\gf$-algebras to formal Lie $\gf$-algebras. Then we have the following result.

\begin{prop}\label{adjoint}
The functor $\overline{\mathcal{U}}:\mathcal{F}_{Lie}Alg \to \mathcal{F}_{Ass}Alg$ is left adjoint to the functor $l:\mathcal{F}_{Ass}Alg \to \mathcal{F}_{Lie}Alg$. 
\end{prop}

\begin{proof}
We need to show that there is a canonical isomorphism
\[\mathcal{F}_{Ass}Alg(\overline{\mathcal{U}}\mathcal{G}, A)\cong \mathcal{F}_{Lie}Alg(\mathcal{G},l(A)).\]
First assume that $\mathcal{G}$ and $A$ are finite rank $\gf$-modules.  Choose a positive integer $N$ for which $A^N=0$.  Then we have the following isomorphisms:
\begin{equation} \label{frank}
\begin{split}
\mathcal{F}_{Ass}Alg(\overline{\mathcal{U}}\mathcal{G}, A) & \cong \mathcal{F}_{Ass}Alg(\overline{\mathcal{U}}\mathcal{G}/(\overline{\mathcal{U}}\mathcal{G})^N, A), \\
& \cong \Hom_{\gf-alg}(I(\mathcal{G}),A), \\
& \cong \mathcal{F}_{Lie}Alg(\mathcal{G},l(A)).
\end{split}
\end{equation}

Now let $\mathcal{G}=\inlim{\alpha}\mathcal{G}_\alpha$, $A=\inlim{\beta}{A}_\beta$ where $\mathcal{G}_\alpha, A_\beta$ are nilpotent Lie $\gf$-algebras and nilpotent associative $\gf$-algebras respectively.  We obtain using \eqref{frank}:
\begin{align*}
\mathcal{F}_{Ass}Alg(\overline{\mathcal{U}}\mathcal{G}, A) & \cong \inlim{\alpha}\dilim{\beta}\mathcal{F}_{Ass}Alg(\overline{\mathcal{U}}\mathcal{G}_\alpha, A_\beta), \\
& \cong \inlim{\alpha}\dilim{\beta}\mathcal{F}_{Lie}Alg(\mathcal{G}_\alpha,l(A_\beta)), \\
& \cong \mathcal{F}_{Lie}Alg(\mathcal{G},l(A)).
\end{align*}
\end{proof}

\begin{example}
Let $\mathcal{G}:=\widehat{L}V$, the pro-free Lie algebra on a profinite $\gf$-module $V$. Then $\widehat{\mathcal{U}}\mathcal{G}=\widehat{T}V$, the pro-free associative $\gf$-algebra on $V$. 
\end{example}

\begin{rem} \label{rem_ccmult}
The above example is especially important for us. It is well known that for $V\in Mod_{\gf}$ the algebra $T(V)$ is in fact a Hopf algebra. The cocommutative coproduct $\Delta: T(V)\rightarrow T(V)\otimes T(V)$ (sometimes called the \emph{shuffle coproduct}) is defined uniquely by the requirement $\Delta$ be a $\gf$-algebra homomorphism and that $\Delta(v)=v\otimes 1+1\otimes v$ for all $v\in V$. Then $L(V)$ is naturally identified with the Lie subalgebra of primitive elements in $T(V)$.

These facts have formal analogues: let $V$ now be a profinite $\gf$-module. Then $\widehat{T}(V)$ has the structure of a (formal) Hopf algebra with the coproduct $\Delta: \widehat{T}(V)\rightarrow \widehat{T}(V)\cotimes \widehat{T}(V)$ which is uniquely specified by the requirement that $\Delta$ be a continuous homomorphism of formal $\gf$-algebras and that  $\Delta(v)=v\otimes 1+1\otimes v$ for all $v\in V$. Then the Lie algebra $\widehat{L}(V)$ is naturally identified with the Lie subalgebra of primitive elements in $T(V)$.
\end{rem}

Next, we need to define the notion of a \emph{formal module} over a formal algebra.

\begin{defi}
\
\begin{enumerate}
\item 
Let $A$ be a formal $\gf$-algebra (commutative or associative) and $V$ be a profinite $\gf$-module. Then $V$ is said to be a left \emph{formal $A$-module} if there is a map of $\gf$-modules $A\cotimes V\rightarrow V; a\cotimes v\mapsto av$ such that the usual associativity axiom holds: $(ab)v=a(bv)$ for all $a,b\in A, v\in V$. Similarly, $V$ is a right formal $A$-module if there is a map of $\gf$-modules $V\cotimes A\rightarrow V; v\cotimes a\mapsto av$ such that $v(ab)=(va)b$ for all $a,b\in A, v\in V$. Finally, $V$ is a formal $A$-bimodule if it is both a right and left formal $A$-module and if $(av)b=a(vb)$ for   $a,b\in A, v\in V$.
\item Let $\mathcal{G}$ be a formal Lie $\gf$-algebra and $V$ be a profinite $\gf$-module. Then $V$ is said to be a  \emph{formal $\mathcal{G}$-module} if there is a map of $\gf$-modules $\mathcal{G}\cotimes V\rightarrow V; g\cotimes v\mapsto gv$ such that
$[g,h]v=g(hv)-h(gv)$.
\end{enumerate}
\end{defi}

Note that Proposition \ref{adjoint} implies that the structure of a formal $\mathcal{G}$-module on $V$ is equivalent to the structure of a formal $\widehat{\mathcal{U}}\mathcal{G}$-module on $V$. 

\end{document}